\documentclass[12pt]{amsart}
\usepackage{amssymb}

\newtheorem{theorem}{Theorem}

\newtheorem{proposition}[theorem]{Proposition}
\newtheorem{conjecture}[theorem]{Conjecture}
\newtheorem{corollary}[theorem]{Corollary}
\newtheorem{remark}[theorem]{Remark}
\newtheorem{definition}{Definition}
\newtheorem{lemma}[theorem]{Lemma}

\newcommand{\Z}{{\mathbb Z}}

\newcommand{\R}{{\mathbb R}}
\newcommand{\C}{{\mathbb C}}

\newcommand{\RP}{{\mathbb RP}}

\begin{document}

{

\title[Complements of discriminants of parabolic singularities]{Complements of discriminants of real parabolic function singularities. II}

\author{V.A.~Vassiliev}
\thanks{This work was supported by the Absorption Center in Science of the Ministry of Immigration and Absorption of the State of Israel}
\address{Weizmann Institute of Science, Rehovot, Israel}
\email{vavassiliev@gmail.com}

\begin{abstract}
We list all connected components of sets of non-dis\-cri\-mi\-nant functions near 
 all {\em parabolic} function singularities (which are the second most important family of singularity classes of smooth functions after {\em simple} singularities). Thus, we prove (and improve in one particular case) all the corresponding conjectures from the previous work \cite{para} with the same title. As an application, we enumerate all {\em local Petrovskii lacunas} near arbitrary parabolic singularities of wavefronts of hyperbolic PDEs. We also show that the complements of the discriminant varieties of the versal deformations of $X_9^{\pm}$ and $P_8^1$ singularities have nontrivial one-dimensional homology groups, in contrast to all simple singularities.

These results are applications of a general method for investigating and separating non-singular perturbations of real function singularities. An important part of this method is a computer program that formalizes local Picard--Lefschetz theory and surgeries of Morse functions.
\end{abstract}

\maketitle 

\section{Introduction}

\subsection{Motivation and terminology} The {\em discriminant variety} of a family of smooth functions that depend on parameters is the set of parameters for which the zero-level sets of the corresponding functions are singular. These varieties appear in many problems of PDE theory and integral geometry, often under the name {\em ``wavefronts''}, see e.g. \cite{ABG}, \cite{AVGZ}, \cite{Kluwer}, \cite{Gord}, \cite{Leray}, \cite{Vas86}, \cite{APLT}, \cite{kashin}. 

Many functions studied in these theories are defined by integral transforms. Their arguments are the parameters of certain function families (called {\em generating families}), and their values are defined by the integrals of certain differential forms along the cycles in the zero-level sets of the corresponding functions. These integral functions are regular outside the discriminant varieties of the generating families. The analytic behavior of such a function near a discriminant point depends on the singularity type of the corresponding function and on the connected component of the complement of the discriminant variety, from which we approach this point. Thus, enumerating these components for different types of discriminant points is one of the first steps in studying the analytical properties of the integral transforms.

\begin{definition}[see \cite{AVGZ}] \rm
A point $a \in \R^n$ is a {\em critical point} of a smooth function $\R^n \to \R$ if all first partial derivatives of this function vanish at $a$. 
A {\it function singularity} is a germ of a regular analytic function $(\R^n,a) \to (\R,0)$ at its critical point with critical value $0$. Two function singularities at points $a$ and $b$ are {\em equivalent} if they can be transformed to each other via a local diffeomorphism $(\R^n, a) \to (\R^n,b)$ (i.e., they have the same expression in appropriate local coordinate systems centered at $a$ and $b$). 
Two function singularities, $f(x_1, \dots, x_n) $ and $\tilde f (\tilde x_1, \dots, \tilde x_m)$, are {\em stably equivalent} if they become equivalent after the summation with non-degenerate quadratic forms in additional variables, that is, $$f(x_1, \dots, x_n) + Q_2(y_1, \dots, y_{N-n}) \sim \tilde f(\tilde x_1, \dots, \tilde x_m) + \tilde Q_2(\tilde y_1, \dots, \tilde y_{N-m})$$ 
for some $N$ and quadratic polynomials $Q_2,$ $\tilde Q_2$ of maximal ranks.
The critical points of analytic functions $\R^n \to \R$ are equivalent if the function singularities obtained from their germs at these points by subtracting the critical values are equivalent.
An $l$-parametric {\em deformation} of a function singularity $f$ is a family $\{f_\lambda\}$ of analytic functions, where $\lambda$ runs through a neighborhood of the origin point $0 \in \R^l$, such that $f_0 \equiv f$ and the function $F(x, \lambda) \equiv f_\lambda(x)$ is regular analytic in a neighborhood of the origin in $\R^{n+l}$. The {\em discriminant} variety of the deformation $F$ is the set of parameter values $\lambda \in \R^l$ such that the corresponding functions $f_\lambda$ have real critical points with zero critical value near the origin in $\R^n$ (or, equivalently, the zero-level sets of the analytic continuations of $f_\lambda$ to $\C^n$ are singular at these points).
\end{definition}

For each isolated function singularity, studying all its deformations reduces to studying any its sufficiently large (so-called {\em versal}, see \cite{AVGZ}) deformation. The complements of discriminant varieties of all versal deformations of stably equivalent singularities are homotopically equivalent, in particular there is a one-to-one correspondence between the connected components of these complements. Therefore, for any stable equivalence class of singularities it is sufficient to list these components for only one versal deformation of one function of this class.

The list of equivalence classes of function singularities begins with {\em simple singularities,} which were identified by V.~Arnold, see \cite{AVGZ}. All of these singularities are stably equivalent to singularities of functions in two variables, and the entire geometric study of their deformations and discriminants reduces to studying such functions. The components of discriminant complements for simple singularities were listed earlier, see \S~\ref{hist}. 

 \begin{table}
\caption{Normal forms of simple, $X_9$ and $J_{10}$ singularities in two variables and of $P_8$ singularities in three variables}
\label{t1}
\begin{center}
\begin{tabular}{|l|l|l|}
\hline
Notation & Normal form & Restriction \\ 
\hline
$A_{2k-1}$ & $\pm x^{2k} \pm y^2 $ & $k \ge 1$ \\ 
$A_{2k}$ & $x^{2k+1} \pm y^2$ & $k \ge 1$ \\
[3pt]
$D_{k}^{\pm}$ & $x^2y \/\pm y^{k-1} $ & $k \ge 4$ \\ [3pt]
$E_6^{\pm}$ & $x^3 \pm y^4 $ & \cr
$E_7$ & $x^3 + x y^3 $ & \cr
$E_8$ & $x^3 + y^5$ & \cr
\hline
\hline 
$X_9^{\pm}$ & $\pm(x^4 + 2A x^2 y^2 + y^4)$ & $|A|<1$ \\
$X_9^1$ & $ x y (x^2 + 2A x y + y^2)$ & $|A|<1$ \\
$X_9^2$ & $x^4 + 2 A x^2y^2 + y^4$ & $A<-1$ \\
\hline 
$J_{10}^1$ & $(x^2+y^2)(x-\gamma y^2)$ & $\gamma \in (-\infty, \infty)$ \\
$J_{10}^3$ & $(x^2-y^2)(x-\gamma y^2)$ & $ \gamma \in (-1,1)$ \\
\hline
$P_8^1$ & $x^3+y^3+z^3 -3A x y z$ & $A<1$ \\
$P_8^2$ & $x^3+y^3+z^3 -3A x y z$ & $A>1$ \\
\hline \end{tabular} \end{center}
\end{table}

 In this work we list all analogous components for the next natural collection of equivalence classes in the hierarchy of function singularities, the {\em parabolic} singularities, which were also distinguished by V.~Arnold. The list of {\em complex} parabolic singularities up to the stable equivalence consists of three classes $P_8$, $X_9$, and $J_{10}$. The classification of their real forms includes more classes: $P_8^1$, $P_8^2$, $X_9^+$, $X_9^-$, $X_9^1$, $X_9^2$, $J_{10}^1$, and $J_{10}^3$. Singularities of classes $X_9$ and $J_{10}$ are stably equivalent to the singularities of functions in two variables, while the $P_8$ singularities are stably equivalent to the non-degenerate cubic forms in three variables. Any real function singularity in two variables of class $X_9^*$ or $J_{10}^*$, as well as any singularity in three variables of one of classes $P_8^*$, has the normal form indicated in Table \ref{t1} in appropriate local coordinates centered at the critical point. Additionally, this table provides the normal forms of the simple singularities, $A_k, D_k, $ and $E_k$.
\smallskip

\begin{table}[h]
\caption{Numbers of rigid isotopy and topological classes of non-discriminant perturbations of parabolic singularities}
\label{t0}
\begin{center}
\begin{tabular}{|l|c|c|c|c|c|c|c|}
\hline
Singularity class & $X_9^{\pm}$ & $X_9^1$ & $X_9^2$ & $J_{10}^1$ & $J_{10}^3$ & $P_8^1$ & $P_8^2$\\
\hline
Rigid isotopy & 7 & 14 & 52 & 13 & 33 & 7 & 15 \\
\hline
Topology & 6 & 12 & 50 & 10 & 31 & & \\
\hline
\end{tabular}
\end{center}
\end{table}

 O.~Viro \cite{viro2} found all {\em topological types} of non-dis\-cri\-mi\-nant perturbations of the $X_9$ and $J_{10}$ singularities. By definition, two small non-dis\-cri\-mi\-nant perturbations $f_\lambda: \R^2 \to \R$ of the same isolated function singularity have the same topological type if there is an isotopy between their zero-level sets in a neighborhood of the critical point that preserves the behavior of these sets at the boundary of the neighborhood. Clearly, the topological type is an invariant of the local connected components of the discriminant complements.
Our ``rigid isotopy'' classification of non-discriminant perturbations, i.e., the classification whose classes are these components, is only slightly stricter than the topological classification, see Table \ref{t0}. Counting the components of the discriminant complements of $P_8$ singularities is related to 
(though is different from) the isotopy classification of nonsingular cubic surfaces in $\R^3$, which is described in detail in \cite{DK}. This classification consists of five and ten classes for $P_8^1$ and $P_8^2$ singularities, respectively.

\subsection{Structure and interpretation of results}

Each class of parabolic singularities is {\em unimodal}, that is, it is swept out by a one-parametric family of orbits of the natural action of the group of local diffeomorphisms $(\R^n, 0) \to (\R^n, 0)$. In each case, we consider the {\em canonical versal deformation} of the corresponding class. This is a space of polynomials that is it diffeomorphic to $\R^9$ (respectively, $\R^{10}$, $\R^8$) in the case of singularities of classes $X_9^*$, (respectively, $J_{10}^*$, $P_8^*$), see formulas (\ref{vers0}), (\ref{vers1}), (\ref{vers3}), (\ref{versP}). The polynomials with parabolic singularities of the corresponding class form a subset that is diffeomorphic to an open interval in this space. All orbits are represented by some points of these intervals. Our main results are the lists of connected components, into which the discriminant varieties divide the parameter spaces of these deformations, see Theorems \ref{teoX0}, \ref{theX91}, \ref{thmX2}, \ref{thmj1}, \ref{thmj3}, \ref{mthmp81}, and \ref{mthmp2}. However, each of these results can be interpreted as a list of all local components of the discriminant complements in arbitrarily small neighborhoods of an arbitrary particular parabolic singularity of the corresponding class. Indeed, one-parametric groups of homogeneous or weighted homogeneous dilations act on the parameter spaces of canonical versal deformations. These groups are naturally isomorphic to the multiplicative group $\R^1_+$ of positive numbers. For $P_8$ singularities the element $\{t\} \in \R^1_+$ maps each of degree three polynomials $f_\lambda(x, y, z)$ to
$t^{-3}f_\lambda (t x, t y, t z)$, in particular it preserves the principal homogeneous part of the polynomials. For $X_9$ singularities it acts as $ f_\lambda (x,y) \to t^{-4}f_\lambda(t x, t y)$, and for $J_{10}$ as $ f_\lambda(x, y) \to t^{-6}f_\lambda(t^2 x, t y)$. All of these actions preserve the discriminants and their complements. Parabolic singularities are defined by the homogeneous or weighted homogeneous polynomials of top degree and are the only invariant points of these actions. Each point of the discriminant complement is reduced by this action into an arbitrarily small neighborhood of some parabolic singularity. Additionally, according to \cite{LoEll}, the discriminant variety is equisingular along the one-dimensional sets of parabolic singularities in parameter spaces of versal deformations. Therefore, each component of the discriminant complement, which is represented in an arbitrarily small neighborhood of some particular parabolic singularity, is also represented in the neighborhoods of all other parabolic singularities of the same class. 

\subsection{On the method} In \cite{para}, an invariant of the connected components of the complements of the discriminant varieties of function singularities is described, complete lists of possible values of this invariant for all parabolic singularities are found, and all of these values for singularities of classes $X_9$ and $J_{10}$ are realized by concrete polynomials. The remaining problems are to list all components with each value of the invariant, and also to realize them for the $P_8$ singularities. In \cite{para}, a conjectural answer to this problem was proposed. It provides the lists
 consisting of 7, 15, 7, 14, 52, 11, and 33 connected components for the singularities of classes $P_8^1$, $P_8^2$, $X_9^{\pm}$, $X_9^1$, $X_9^2$, $J_{10}^1$, and $J_{10}^3$ respectively. Below, we prove that almost all of these lists are indeed complete. However, some two sets of perturbations of the $J_{10}^1$ singularities indicated in \cite{para} consist of pairs of different connected components. Thus, the number 11 is replaced in Table \ref{t0} by 13.

\subsection{On the discriminant complements for simple singularities} 
\label{hist} E.~Looijenga \cite{Lo} implicitly listed the local connected components of the complements of the discriminants of simple singularities, although the numbers of these components remained unknown. O.~Viro \cite{viro2} listed all possible topological types of non-discriminant functions representing these components. In \cite{isr}, it was shown that each topological type on this list corresponds to exactly one connected component of the discriminant complement. For simple singularities with Milnor numbers at most six this fact also follows from the results of \cite{sed}.
Unfortunately, I only read the work \cite{viro2}, which contains an enumeration of the topological types of perturbations of simple and parabolic singularities in two variables, after the publication of \cite{para} and \cite{isr}. 

\subsection{1-cohomology of discriminant complements}
E.~Looijenga \cite{Lo} pro\-ved that all local connected components of discriminant complements of real simple singularities are homotopically trivial. Conversely, we demonstrate that at least three such components of parabolic singularities of each of the classes $X_9^+$ and $X_9^-$, as well as one component of the $P_8^1$ class, have nontrivial one-dimensional homology groups; see Theorems \ref{t1hom} and \ref{t2hom}.

\subsection{An application to hyperbolic PDEs.}
In the theory of hyperbolic PDEs, {\em local Petrovskii lacunas} are the local connected components of wavefront complements where wave functions behave regularly (see, for example, \cite{ABG}, \cite{Gord}, \cite{Leray}, \cite{Petr}, \cite{APLT}). Near the non-singular points of wavefronts (corresponding to the Morse singularities of generating functions) the local lacunas were enumerated in \cite{Dav} and \cite{Bor}, see also \cite{ABG}. Near the simplest non-Morse points of classes $A_2$ and $A_3$, this was done in \cite{Gord}, and near all simple singularities in \cite{Vas86}. In \cite{kashin}, a list of Petrovskii lacunas at parabolic singular points of wavefronts was presented, and it was conjectured that this list is complete. We prove this conjecture for almost all parabolic classes but discover one additional local lacuna in the $P_8^2$ case, see \S~\ref{lacu}. 

\subsection*{A notation} The symbol $\Box$ indicates either the conclusion of a proof or the absence of a proof, as in the case of immediate corollaries or statements supplied with references to the works where they are proven.

\section{Results for the $X_9^*$ classes}

\subsection{} Real parabolic function singularities of type $X_9$ form four one-parametric families of stable equivalence classes $X_9^+$, $X_9^-$, $X_9^1$, and $X_9^2$, see Table \ref{t1}.

The {\em standard versal deformations} of these singularities in all four cases have the same form: they are families of polynomials
\begin{equation} f_A(x) + \lambda_1 + \lambda_2 x + \lambda_3 y + \lambda_4 x^2 + \lambda_5 x y + \lambda_6 y^2 + \lambda_7 x^2 y + \lambda_8 x y^2 ,
\label{vers0}
\end{equation}
depending on nine parameters $\lambda_1, \dots, \lambda_8,$ and $A$, where $f_A$ is the corresponding normal form from Table \ref{t1}. According to versality theory (see \cite{AVGZ}, \S~8), studying all other deformations of singularities of stable equivalence classes $X_9^*$ can be reduced to studying these families. The entire study of $X_9^-$ singularities and their deformations can be reduced to that of $X_9^+$ via the multiplication by $-1$.

\begin{definition} \rm
A polynomial $f: ({\mathbb C}^2, {\mathbb R}^2) \to ({\mathbb C}, {\mathbb R})$ of the form (\ref{vers0}) is {\em generic} if it has only Morse critical points in ${\mathbb C}^2$, all of its critical values are different and {\em not equal to 0}, and the complex zero-level set of the principal homogeneous part of $f$ consists of four different lines in ${\mathbb C}^2$. Denote by $\Phi_+$, $\Phi_1$, and $\Phi_2$ the spaces of all polynomials of the form (\ref{vers0}) with principal parts of classes $X_9^+$, $X_9^1$, and $X_9^2$, respectively.
\end{definition}

The parameters $\{\lambda_1, \dots, \lambda_8, A\}$ identify the spaces $\Phi_+$ and $\Phi_1$ with $\R^8 \times (-1,1) $, and the space $\Phi_2$ with $\R^8 \times (-\infty, -1) $. The discriminant subvarieties divide the spaces $\Phi_*$ into connected components. Enumerating these components is our goal in this section.

\subsection{Virtual functions and virtual components}
\label{virtx}
In \S 2 of \cite{para} we introduced the {\em set-valued invariant} of connected components of the discriminant complements. Let us briefly review this construction in the case of $X_9$ singularities.

A {\em virtual function associated with a generic real polynomial} $f$ is a collection of 
topological characteristics of the complex zero-level set of this polynomial. Specifically, these data of a generic polynomial $f(x,y)$ of some class $X_9^*$ consists of 
\begin{enumerate}
\item
the $9 \times 9$ \ matrix of the intersection indices of canonically ordered and oriented {\em vanishing cycles} forming the basis of the two-dimensional homology group of
the submanifold $V_f \subset \C^3$ given by the equation $f(x, y)+z^2=0$, 
\item
the string of intersection indices of these nine cycles with the set of real points in $V_f$, 
\item
the Morse indices of real critical points of $f$, and 
\item the numbers of positive and negative critical values of $f$. 
\end{enumerate}

\begin{remark} \rm
\label{rem111}
Besides the function $f$, a virtual function depends on the collection of paths in $\C^1$ defining the base of vanishing cycles in the homology of the level set. If the function has no more than one pair of non-real critical values, then there is only one canonical system of paths, and hence only one associated virtual function.

If all critical points of a polynomial are real, then the intersection indices from item (2) are easily determined by the remaining elements (1), (3), and (4) of a virtual function and hence may not be included in its description, see Theorem V.1.3 in \cite{APLT}.
\end{remark}

The (abstract) {\em virtual functions} of class $\Phi_*$ include all virtual functions associated with generic real polynomials of this class, as well as all analogous collections of data that can be obtained from them via standard flips that model elementary surgeries of real functions. The {\em formal graph} of any type $X_9^*$ is the graph, whose vertices correspond to all the virtual functions of the corresponding class $\Phi_*$, and whose edges correspond to the standard flips. The {\em virtual components} of the formal graph are the connected subgraphs, into which the formal graph splits when all edges modeling the surgeries that cross the discriminant are removed. The {\em set-valued invariant} of a generic polynomial $f$ of class $\Phi_*$ is the virtual component of the formal graph containing the virtual functions associated with this polynomial. This virtual component is denoted by $S(f)$. \smallskip

These objects have the following obvious property.

\begin{proposition}
\label{protrivi}
If two generic polynomials, $f$ and $\tilde f$, 
belong to the same connected component of the complement of the discriminant, then the virtual components $S(f)$ and $S(\tilde f)$ are the same. \hfill $\Box$
\end{proposition}

The formal graphs of all simple and parabolic singularities are finite.

Enumerating the virtual components of such a graph is a combinatorial problem. An (available) computer program solving it is described in \cite{APLT} and \cite{para}. The work \cite{para} provides the results of these computations for all parabolic singularities, as well as the immediate consequences of these computations in the terms of real functions. The remaining problem is to find all real connected components associated with each virtual component.

The results of \cite{Jaw2} imply that in the case of parabolic singularities all standard flips of virtual functions can be realized by real surgeries of generic real functions. In particular, the following statement holds.

\begin{proposition}[see \cite{vasX9}, Proposition 2]
\label{propmain}
For any class of polynomials $\Phi_+$, $\Phi_1$, or \ $\Phi_2$,

a$)$ every virtual function of this class is associated with a generic real polynomial of the same class;

b$)$ for any generic polynomial $f$ of this class, every virtual function from the virtual component $S(f)$ is associated with some generic real polynomial $\tilde f$ from the same connected component of the complement of the discriminant as $f$. \hfill $\Box$
\end{proposition}

\begin{remark} \rm
In the statement of Proposition 2 of \cite{vasX9}, the functions $f$ and $ \tilde f$ belong to the broader class of all degree-four Morse polynomials with non-degenerate principal homogeneous parts. However, the entire construction proving this proposition in \cite{vasX9} is performed within the spaces of the versal deformations (\ref{vers0}). Thus, it also proves the above Proposition \ref{propmain}.
\end{remark}

These constructions and statements allow us to prove
upper bounds of the numbers of isotopy classes of generic Morse perturbations associated with particular virtual functions. These bounds are given in Propositions \ref{estx+}, \ref{estx1}, and \ref{estx2} for classes $X_9^+$, $X_9^1$ and $X_9^2$ separately. Consequently, we obtain estimates of the numbers of connected components of the real discriminant complement which are associated with each virtual component. Specifically, in the cases $X_9^+$ and $X_9^2$, all real connected components of discriminant complements associated with the same virtual component can be obtained from each other by the action of a transformation group of order eight. In the case $X_9^1$ the analogous group is of order four. Then, for each virtual component in the list from \cite{para}, we compare all real components associated with this virtual component and determine which of these components are identical. This comparison provides the results of the present work concerning the $X_9$ singularities. The strategy for the $J_{10}$ and $P_8$ singularities is exactly the same.

\begin{remark} \rm
In \cite{vasX9}, \cite{vasP8}, and \cite{vasJ10}, we use the same strategy with a different definition of virtual components. In these works, we classify the generic perturbations of function singularities up to isotopy within the space of Morse functions instead of the space of non-discriminant functions. Thus, the virtual components of the formal graph are determined there by removing the edges that change the Morse type.
\end{remark}

\begin{definition} \rm 
\label{defind}
For any generic polynomial $f: \R^n \to \R$, \ the number \ $\mbox{Ind}(f)$ \ equals the number of its real critical points with a negative critical value and an even Morse index minus the number of critical points with a negative critical value and an odd Morse index.
\end{definition}

 It is easy to see that this number is an invariant of the connected components of the discriminant complements, and also of associated virtual components.

\subsection{Enumeration of components for $X_9^+$ singularities}
\label{sec90}
In this subsection, we prove the following statement, conjectured in \cite{para}.

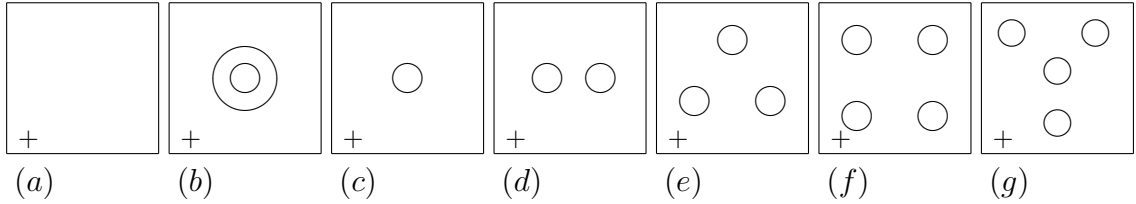
\begin{figure}
\unitlength 0.90mm
\linethickness{0.4pt}
\begin{picture}(20.00,25.00)
\put(0,5){\line(1,0){20}}
\put(0,5){\line(0,1){20}}
\put(0,25){\line(1,0){20}}
\put(20,5){\line(0,1){20}}
\put(1.5,6){\footnotesize $+$}
 \put(1,0){$(a)$}
\end{picture}
\begin{picture}(20.00,25.00)
\put(0,5){\line(1,0){20}}
\put(0,5){\line(0,1){20}}
\put(0,25){\line(1,0){20}}
\put(20,5){\line(0,1){20}}
\put(10,15){\circle{4}}
\put(10,15){\circle{8}}
 \put(1,0){$(b)$}
\put(1.5,6){\footnotesize $+$}
\end{picture}
\begin{picture}(20.00,25.00)
\put(0,5){\line(1,0){20}}
\put(0,5){\line(0,1){20}}
\put(0,25){\line(1,0){20}}
\put(20,5){\line(0,1){20}}
\put(10,15){\circle{4}}
 \put(1,0){$(c)$}
\put(1.5,6){\footnotesize $+$}
\end{picture}
\begin{picture}(20.00,25.00)
\put(0,5){\line(1,0){20}}
\put(0,5){\line(0,1){20}}
\put(0,25){\line(1,0){20}}
\put(20,5){\line(0,1){20}}
\put(7,15){\circle{4}}
\put(14,15){\circle{4}}
\put(1.5,6){\footnotesize $+$}
 \put(1,0){$(d)$}
\end{picture}
\begin{picture}(20.00,25.00)
\put(0,5){\line(1,0){20}}
\put(0,5){\line(0,1){20}}
\put(0,25){\line(1,0){20}}
\put(20,5){\line(0,1){20}}
\put(5,12){\circle{4}}
\put(15,12){\circle{4}}
\put(10,20){\circle{4}}
 \put(1,0){$(e)$}
\put(1.5,6){\footnotesize $+$}
\end{picture}
\begin{picture}(20.00,25.00)
\put(0,5){\line(1,0){20}}
\put(0,5){\line(0,1){20}}
\put(0,25){\line(1,0){20}}
\put(20,5){\line(0,1){20}}
\put(5,10){\circle{4}}
\put(5,20){\circle{4}}
\put(15,10){\circle{4}}
\put(15,20){\circle{4}}
\put(1.5,6){\footnotesize $+$}
 \put(1,0){$(f)$}
\end{picture}
\begin{picture}(20.00,25.00)
\put(0,5){\line(1,0){20}}
\put(0,5){\line(0,1){20}}
\put(0,25){\line(1,0){20}}
\put(20,5){\line(0,1){20}}
\put(10,16){\circle{3.5}}
\put(4,21){\circle{3.5}}
\put(10,9){\circle{3.5}}
\put(15,21){\circle{3.5}}
\put(1.5,6){\footnotesize $+$}
 \put(1,0){$(g)$}
\end{picture}
\caption{Perturbations of $+X_9$ singularity}
\label{x9+}
\end{figure}

\begin{theorem}
\label{teoX0}
There are exactly seven connected components in the space of non-discriminant polynomials of class $\Phi_+$. The shapes of the zero-level sets of some polynomials representing these components are shown in Fig.~\ref{x9+}.
\end{theorem}

\begin{remark} \rm
Pictures $(f)$ and $(g)$ of Fig.~\ref{x9+} are topologically equivalent. However, Bezout's theorem arguments show that the polynomials that realize them cannot be connected within the space of non-discriminant polynomials of class $\Phi_+$. Figure \ref{x9+} essentially coincides with the well-known list of affine rigid isotopy classes of bounded degree four smooth algebraic plane curves. In our case, however, the isotopies are allowed only within the space $\Phi_+$ which is much smaller than the space of all degree four polynomials. In other cases, $\Phi_1$ and $\Phi_2$, the analogous restriction leads to an increase of classification lists.
\end{remark}

The starting point of the proof of Theorem \ref{teoX0} is the following proposition.

\begin{proposition}[see \cite{para}]
\label{proX0}
There are exactly seven virtual components of class $\Phi_+$. The shapes of the zero-level sets of some polynomials of this class that are associated with virtual functions representing these virtual components are shown in Fig.~\ref{x9+}. \hfill $\Box$
\end{proposition}

\begin{definition} \rm
Denote by $G_0$ the group of linear transformations of $\R^2$ generated by reflections $(x,y) \mapsto (x, -y)$ and $(x,y) \mapsto (y,x)$ 
\end{definition}

Clearly, this group is isomorphic to the group of symmetries of the square, in particular is of order eight. This group acts on the space $\Phi_+$ (each transformation $\Gamma: \R^2 \to \R^2$ moves the polynomial $f$ to $f \circ \Gamma$). This action preserves the set of generic polynomials and the set-valued invariant. 

\begin{proposition}
\label{estx+}
If two generic polynomials, $f$ and $\tilde f$, of class $\Phi_+$ are associated with the same virtual function, then the connected components of the space of generic polynomials of class $\Phi_+$ containing $f$ and $\tilde f$ are mapped to each other by an element of the group $ G_0$.
\end{proposition}

\begin{definition}\rm
For any polynomial $f: \R^2 \to \R$ of degree four with principal homogeneous part $x^4 + 2A x^2y^2 + y^4$, denote by ${\bf T}(f)$ the unique polynomial of the form (\ref{vers0}) with the same principal part, which is obtained from $f$ by a parallel translation of $\R^2$. (This is the translation by the vector $\frac{1}{4}(a, b),$ where $a$ and $b$ are the coefficients of $f$ at the monomials $x^3$ and $y^3$).
\end{definition}

\begin{lemma}[see \cite{vasX9}, Lemma 4]
\label{lemx+}
A polynomial of the form $($\ref{vers0}$)$ with a positive principal homogeneous part has a critical point of class $E_6$ with critical value 0 and a critical point of class $A_3$ with critical value $-\frac{27}{256}$ if and only if it is one of the four polynomials 
\begin{equation}
\label{sta+}
\hspace{2.7cm}{\bf T}(x^4 \pm x^3 + y^4) \quad \mbox{ or } \quad {\bf T}(x^4 + y^4 \pm y^3) \ . \hspace{2.7cm} \Box
\end{equation}
\end{lemma}

\noindent
{\it Proof of Proposition \ref{estx+}.} 
\label{proofestx+} This proof almost coincides with the proof of Proposition 3 in \S 3.5 of \cite{vasX9}. In the next paragraph, we present the outline of this proof; we refer to \cite{vasX9} for all common details. 

Using the {\em Lyashko--Looijenga covering} (in the form described in \cite{Jaw2}), we can assume that the polynomials $f$ and $\tilde f$ have the same set of critical values and the same set of paths in $\C^1$ specifying the common associated virtual function. 
Let $I: [0,1] \to \Phi_+$ be a generic path such that $I(0)=f$, $I(1)$ is the polynomial 
${\bf T}(x^4 + x^3 + y^4)$, and the path $I$ has in its interior points only finitely many transversal intersections with the variety of non-generic polynomials in $\Phi_+$ (that is, with the union of the discriminant and the set of polynomials with fewer than nine critical values). Since the virtual function of $\tilde f$ is the same as that of $f$, there exists another path $\tilde I: [0,1] \to \Phi_+$ such that $\tilde I(0) = \tilde f$, the polynomials $I(t)$ and $\tilde I(t)$ have the same set of critical values for any $t \in [0,1]$, and the polynomials of the path $\tilde I$ undergo the surgeries of the same type as $I$ at the same instants $t$. (For this existence statement, the parabolicity of $X_9$ singularities is essential. The proof of it is based on Proposition 2 of \cite{Jaw2}.) In particular, the point $\tilde I(1)$ can only be one of the four points (\ref{sta+}). First, suppose that it is the same point $I(1)$. A neighborhood of this point in the space $\Phi_+$ has only one nontrivial smooth automorphism that preserves the associated virtual functions and commutes with the Lyashko--Looijenga map. This automorphism is induced by the reflection of $\R^2$ in the line $y=0$, see \cite{vasX9}. Therefore, there are only two possibilities: the final segments of the paths $I$ and $\tilde I$ (for $t \in (1-\varepsilon, 1]$) are either the same or are mapped to each other by this automorphism. In both cases, this condition can be continued to the entire paths $I$ and $\tilde I$ and proves that either $\tilde f = f$ or these two functions are mapped to each other by this automorphism, which belongs to the group $G_0$. If the point $\tilde I(1) $ is different from $I(1)$, we first apply the element of the group $G_0$, moving these points to each other, to the entire path $\tilde I$ (in particular, to its starting point $\tilde f$). This leads us to the situation considered previously. \hfill $\Box$
\medskip

\noindent
{\it Proof of Theorem \ref{teoX0}.} By Proposition \ref{estx+}, any virtual component mentioned in Proposition \ref{proX0} is associated with at most eight real connected components. Let us prove that in fact each of these seven virtual components is associated with a single real connected component (that is, these components are $G_0$-invariant). 

Pictures $(a)$, \ $(c)$, \ $(b)$, and \ $(f)$ \ of Fig.~\ref{x9+} are realized, respectively, by 
polynomials $$x^4 + y^4 + t, $$
$$x^4+y^4 - t, $$
$$x^4 + y^4 - t (x^2+y^2) + \alpha t^2, \qquad \alpha \in (0, 1/4) , $$
and
$$x^4+y^4-t(x^2+y^2) + \alpha t^2, \quad \alpha \in (1/4, 1/2) ,$$
with arbitrary $t >0$. All of these polynomials are invariant under the entire group $G_0$.

Picture \ $(d)$ \ can be represented by the polynomial
\begin{equation}
x^4+2Ax^2y^2 +y^4 + 2 x^2 - 2 y^2+ 1/2 
\label{u1}
\end{equation} 
with arbitrary $A \in (-1,1)$.
This polynomial is invariant under the transformations $(x,y) \mapsto (x,-y)$ and $(x,y) \mapsto (-x,y)$. Thus, its $G_0$-orbit consists only of it and the polynomial 
\begin{equation}
\label{u2}
x^4+2Ax^2y^2 +y^4 - 2 x^2 + 2 y^2+ 1/2 \ .
\end{equation}
The polynomials (\ref{u1}) and (\ref{u2}) are connected by the family of polynomials 
\begin{equation}
\label{fam1}
x^4+2Ax^2y^2 + y^4 + 2 (x \cos \varphi + y \sin \varphi )^2 - 2 (x \sin \varphi - y \cos \varphi )^2+ 1/2,
\end{equation}
that depend on the parameter $\varphi \in [0, \pi/2]$.

\begin{lemma}
\label{lem10}
If $A<1$ is sufficiently close to $1$, then all polynomials of the family $(\ref{fam1})$ are non-discriminant, in particular their zero-level sets are all smooth and isotopic to each other. 
\end{lemma} 

\noindent
{\it Proof.} If we take $A=1$ in (\ref{fam1}), then all these polynomials will be obtained from each other by the rotations of $\R^2$ by angles $\varphi \in [0,\pi/2]$. Each of these polynomials has only three real Morse critical points: one at the origin and two on the unit circle. After replacing $A=1$ by $1-\varepsilon$ with sufficiently small $\varepsilon$, the real zero-level set deforms only slightly in the domain $\{x^2+y^2 \leq 2\}.$ In particular, it maintains its topological structure in this domain for any $\varphi$. It also remains empty outside this domain. \hfill $\Box$
\medskip

Thus, the connected component of the space of non-discriminant polynomials containing (\ref{u1}) is also invariant under the transformation $(x,y) \mapsto (y,x)$.
\medskip

The topological type of the curve shown in picture \ $(e)$ \ can be realized by any of the four polynomials 
\begin{equation}
\label{uu1}
{\bf T}(x^4+2Ax^2y^2+y^4 \pm(x^3-3x y^2)+ \varkappa) 
\end{equation} 
and
\begin{equation}
{\bf T}(x^4+2Ax^2y^2+ y^4\pm (y^3-3x^2y) + \varkappa) 
\end{equation}
with arbitrary $A \in (-1,1)$ and sufficiently small $\varkappa >0$ (where the restriction on $\varkappa$ depends on $A$).
These four polynomials form an orbit of the group $G_0$. All of them belong to the one-parametric family of polynomials 
\begin{equation}
\label{fam2}
{\bf T}(x^4+ 2Ax^2y^2+ y^4 + \mbox{Re}(e^{i \tau} (x+ i y)^3)) + \varkappa), 
\end{equation} 
$\tau \in [0,2\pi]$. For $A=1$, all these polynomials are obtained from each other by the rotations of $\R^2$. They are bounded from below, have a critical point of class $D_4^-$ with critical value $\varkappa$ at the origin, three real minimum points with critical value $-\frac{27}{256}+\varkappa$, and no other critical points. Therefore, for $\varkappa\in \left(0, \frac{27}{256}\right)$ all these polynomials realize the topological picture \ $(e)$. 
Analogously to the proof of Lemma \ref{lem10}, taking $A<1$ very close to 1 will not change these topological pictures for any $\tau$. So, the connected component of the discriminant complement containing these polynomials is invariant under the action of the group $G_0$.
\medskip

The picture \ $(g)$ \ can be realized by either of the four polynomials 
\begin{equation}
\label{weq}
 {\bf T}(x^4 + (2-\varepsilon) x^2y^2 + y^4 \pm 10 x(x^2-3y^2) + 15 (x^2 + y^2) - 5 ) \ ,
\end{equation} 
\begin{equation}
\label{weq8}
 {\bf T}(x^4 + (2-\varepsilon)x^2y^2+ y^4 \pm 10 y(y^2-3x^2) + 15 (x^2 + y^2) - 5 )
\end{equation}
with sufficiently small $\varepsilon>0$.
 These four polynomials form an orbit of the group $G_0$. 

The family of polynomials 
\begin{equation}
\label{fam3}
{\bf T}(x^4+(2-\varepsilon)x^2y^2+ y^4 + 10 \mbox{Re}\left(e^{i\tau} (x+i y)^3\right) + 15 (x^2+y^2) - 5), 
\end{equation}
$\tau \in [0, 2\pi],$ contains all four polynomials (\ref{weq}), (\ref{weq8}) for $\tau = 0,$ $\frac{\pi}{2}$, $\pi$, and $\frac{3\pi}{2}$. 
All elements of this family lie in the same connected component of the complement of the discriminant. This completes the proof of Theorem \ref{teoX0}. \hfill $\Box$
\medskip

\begin{remark} \rm
To realize the pictures ($d$), \ ($e$), \ and \ ($g$) \ by arbitrarily small perturbations of the $X_9^+$ singularities, we can multiply by $t^{4-k}$ all monomials of degree $k$ in the formulas (\ref{u1}), (\ref{uu1}), and (\ref{weq}).
\end{remark}

\subsection{The 1-cohomology of discriminant complements of $X_9^+$ singularities} Recall the notation $B(\R^2, k)$ for the {\em unordered configuration space}, whose points are arbitrary $k$-element subsets of the plane. For any $k \geq 2$, $H^1(B(\R^2,k)) \simeq \Z$.
The basic element of this group is induced from the generator of the group $H^1(\C^1 \setminus \{0\}) \sim \Z $ by the map $B(\R^2, k) \equiv B(\C^1,k) \to \C^1 \setminus \{0\}$ sending any collection of \ $k$ \ distinct point $z_1, \dots, z_k$ to the product of all \ $k(k-1)$ \ complex numbers $z_i - z_j$, $i \neq j$. Let $U$ be a connected component of the discriminant complement in $\Phi_+$ that consists of polynomials whose zero-level sets are unions of $k \geq 2$ non-nested ovals. Then a map $H^*(B(\R^2,k)) \to H^*(U) $ is defined as follows. For any polynomial $f \in U,$ define $\beta (f)$ as the space of all choices of a point inside any of its ovals. These spaces $\beta(f)$ form a locally trivial fiber bundle over the component $U$. The fibers of this fiber bundle are contractible, hence it has a cross-section, which is a continuous map $U \to B(\R^2, k)$. Clearly, these maps defined by different cross-sections are homotopic and define the same homomorphism of cohomology rings. 

Define the class $W \in H^1(U)$ as the class obtained by this map from the generator of the group $H^1(B(\R^2,k))$.

\begin{theorem}[cf. \cite{VS}, \S~1.8.3]
\label{t1hom}
The class $W$ takes unit values on some 1-cycles in each of the connected components containing the polynomials shown in pictures $($d$)$, \ $($e$)$, \ and \ $($g$)$ \ of Fig.~\ref{x9+}. 
\end{theorem}

\noindent
{\it Proof.} These cycles are the families (\ref{fam1}), (\ref{fam2}), and (\ref{fam3}). \hfill $\Box$

\subsection{The components for the $X_9^1$ singularities}
\label{sec91}

Define the group $G_1$ as the group of linear transformations of $\R^2$ generated by the symmetries $(x,y) \mapsto (y, x)$ and $(x,y) \mapsto (-x, -y)$. Clearly, this group is isomorphic to $\Z_2^2$. It acts on the space $\Phi_1$ of polynomials (\ref{vers0}) with principal part 
\begin{equation}
\label{equas}
x y (x^2 + 2A x y + y^2), \quad |A|<1,
\end{equation} 
 and preserves the set of generic polynomials. Polynomials obtained from each other by this action have equal associated virtual functions and so are not separated by the set-valued invariant.
\medskip

Additionally, the involution 
\begin{equation}
\label{invol}
f(x,y) \mapsto -f(x,-y)
\end{equation}
acts on the space $\Phi_1$. It preserves the discriminant variety and permutes the connected components of its complement. This involution naturally extends to the set of associated virtual components: each virtual component which is the value of the set-valued invariant of the polynomial $f(x,y)$ goes to the value of the polynomial $-f(x,-y)$. 
This involution has no fixed elements: indeed, the sum of the values of the invariant \ Ind \ (see Definition \ref{defind}) for any two components (real or virtual) related by the involution is always equal to $-1$.

 For any polynomial $f$ with a principal homogeneous part of the form (\ref{equas}),
its zero-level set goes to the infinity along four branches: two vertical and two horizontal. If the polynomial is non-discriminant, then each vertical branch is matched with the horizontal branch belonging to the same component of the zero-level set. The involution (\ref{invol}) acts on the zero-level sets as the reflection in the line $\{y=0\}$. In particular, it always alters the matching of the vertical and horizontal branches.

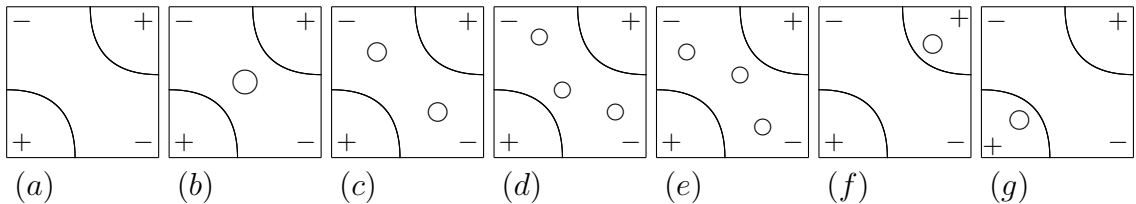
\begin{figure}
\unitlength 0.90mm
\linethickness{0.4pt}
\begin{picture}(20.00,25.00)
\put(0,5){\line(1,0){20}}
\put(0,5){\line(0,1){20}}
\put(0,25){\line(1,0){20}}
\put(20,5){\line(0,1){20}}
\bezier{300}(11,25)(11,16)(20,16)
\bezier{300}(9,5)(9,14)(0,14)
\put(2,7){\makebox(0,0)[cc]{\tiny $+$}}
\put(2,23){\makebox(0,0)[cc]{\tiny $-$}}
\put(18,23){\makebox(0,0)[cc]{\tiny $+$}}
\put(1,0){$(a)$}
\end{picture}
\begin{picture}(20.00,25.00)
\put(0,5){\line(1,0){20}}
\put(0,5){\line(0,1){20}}
\put(0,25){\line(1,0){20}}
\put(20,5){\line(0,1){20}}
\bezier{300}(11,25)(11,16)(20,16)
\bezier{300}(9,5)(9,14)(0,14)
\put(10,15){\circle{3}}
\put(2,7){\makebox(0,0)[cc]{\tiny $+$}}
\put(2,23){\makebox(0,0)[cc]{\tiny $-$}}
\put(18,23){\makebox(0,0)[cc]{\tiny $+$}}
\put(1,0){$(b)$}
\end{picture}
\begin{picture}(20.00,25.00)
\put(0,5){\line(1,0){20}}
\put(0,5){\line(0,1){20}}
\put(0,25){\line(1,0){20}}
\put(20,5){\line(0,1){20}}
\bezier{300}(11,25)(11,16)(20,16)
\bezier{300}(9,5)(9,14)(0,14)
\put(6,19){\circle{2.5}}
\put(14,11){\circle{2.5}}
\put(2,7){\makebox(0,0)[cc]{\tiny $+$}}
\put(2,23){\makebox(0,0)[cc]{\tiny $-$}}
\put(18,23){\makebox(0,0)[cc]{\tiny $+$}}
\put(1,0){$(c)$}
\end{picture}
\begin{picture}(20.00,25.00)
\put(0,5){\line(1,0){20}}
\put(0,5){\line(0,1){20}}
\put(0,25){\line(1,0){20}}
\put(20,5){\line(0,1){20}}
\bezier{300}(11,25)(11,16)(20,16)
\bezier{300}(9,5)(9,14)(0,14)
\put(9,14){\circle{2}}
\put(16,11){\circle{2}}
\put(6,21){\circle{2}}
\put(2,7){\makebox(0,0)[cc]{\tiny $+$}}
\put(2,23){\makebox(0,0)[cc]{\tiny $-$}}
\put(18,23){\makebox(0,0)[cc]{\tiny $+$}}
\put(1,0){$(d)$}
\end{picture}
\begin{picture}(20.00,25.00)
\put(0,5){\line(1,0){20}}
\put(0,5){\line(0,1){20}}
\put(0,25){\line(1,0){20}}
\put(20,5){\line(0,1){20}}
\bezier{300}(11,25)(11,16)(20,16)
\bezier{300}(9,5)(9,14)(0,14)
\put(11,16){\circle{2}}
\put(14,9){\circle{2}}
\put(4,19){\circle{2}}
\put(2,7){\makebox(0,0)[cc]{\tiny $+$}}
\put(2,23){\makebox(0,0)[cc]{\tiny $-$}}
\put(18,23){\makebox(0,0)[cc]{\tiny $+$}}
\put(1,0){$(e)$}
\end{picture}
\begin{picture}(20.00,25.00)
\put(0,5){\line(1,0){20}}
\put(0,5){\line(0,1){20}}
\put(0,25){\line(1,0){20}}
\put(20,5){\line(0,1){20}}
\bezier{300}(11,25)(11,16)(20,16)
\bezier{300}(9,5)(9,14)(0,14)
\put(15,20){\circle{2.5}}
\put(2,7){\makebox(0,0)[cc]{\tiny $+$}}
\put(2,23){\makebox(0,0)[cc]{\tiny $-$}}
\put(18.5,23.5){\makebox(0,0)[cc]{\tiny $+$}}
\put(1,0){$(f)$}
\end{picture}
\begin{picture}(20.00,25.00)
\put(0,5){\line(1,0){20}}
\put(0,5){\line(0,1){20}}
\put(0,25){\line(1,0){20}}
\put(20,5){\line(0,1){20}}
\bezier{300}(11,25)(11,16)(20,16)
\bezier{300}(9,5)(9,14)(0,14)
\put(5,10){\circle{2.5}}
\put(1.5,6.5){\makebox(0,0)[cc]{\tiny $+$}}
\put(2,23){\makebox(0,0)[cc]{\tiny $-$}}
\put(18,23){\makebox(0,0)[cc]{\tiny $+$}}
\put(1,0){$(g)$}
\end{picture}
\caption{Perturbations of $X_9^1$ singularity}
\label{x91}
\end{figure}

\begin{theorem}
\label{theX91}
There are exactly 14 connected components of the complement of the discriminant in the space of polynomials $($\ref{vers0}$)$ of class $\Phi_1$. The shapes of the zero-level sets of polynomials representing some seven of them are shown in Fig.~\ref{x91}, and the remaining seven can be obtained from these by the involution $($\ref{invol}$)$.
\end{theorem}

The proof of this theorem is based on the following two propositions.

\begin{proposition}[see \cite{para}]
\label{proX91}
There are exactly ten virtual components of type $X_9^1$. They split into five pairs of virtual components related to each other by the involution $($\ref{invol}$)$. The topological types of zero-level sets of polynomials associated with elements 
of some five of the virtual components $($one from each pair$)$ are shown in pictures \ $($a$)$, \ $($b$)$, \ $($c$)$, \ $($d$)$, \ and \ $($f$)$ \ of Fig.~\ref{x91}. \hfill $\Box$
\end{proposition}

\begin{proposition}
\label{estx1}
If two generic polynomials $f$ and $\tilde f$ of class $\Phi_1$ are associated with the same virtual function, then their connected components in the space of generic polynomials are moved to each other by an element of the group $ G_1$.
\end{proposition}

\begin{definition}\rm
For any polynomial $f: \R^2 \to \R$ of degree four with the principal homogeneous part 
(\ref{equas}), denote by ${\bf \tilde T}(f)$ the unique polynomial of the form (\ref{vers0}) with the same principal part, which is obtained from $f$ by a parallel translation of the plane. (This is the translation by the vector $ ( \alpha, \beta),$ where $\alpha$ and $\beta$ are the coefficients of $f$ at the monomials $y^3$ and $x^3$).
\end{definition}

\begin{lemma}
\label{lemx1}
A polynomial of the form $($\ref{vers0}$)$ with the principal part $($\ref{equas}$)$ has a critical point of class $E_7$ with zero critical value and a critical point of class $A_2$ with critical value $3^{\frac{15}{2}}2^{-8}$ if and only if it is one of four polynomials 
\begin{equation}
\label{ea1}
{\bf \tilde T}(x^3y + \sqrt{3} x^2y^2+ x y^3 \pm x^3) \ \mbox{ \rm or } \
 {\bf \tilde T}(x^3y + \sqrt{3}x^2y^2 + x y^3 \pm y^3)\ .
\end{equation}
\end{lemma}

\noindent
{\it Proof.} Let $f$ be a polynomial that satisfies the conditions of this lemma. By Lemma 5 of \cite{vasX9}, the $j$-invariant of its principal homogeneous part is equal to 0. By the formula for $j$-invariants of quartic forms (see \cite{Mukai}, \S 1.3), $A = \pm \frac{\sqrt{3}}{2}$ in (\ref{equas}). We can assume that the coordinates $(x, y)$ are centered at the $E_7$ point. In these coordinates, the monomials of degree $\leq 2$ of the polynomial $f$ vanish, and the homogeneous part of degree three is proportional to the third degree of a linear function that also divides the degree four homogeneous part of $f$. Thus, the degree three homogeneous part has the form $\alpha x^3$ or $\alpha y^3$, $\alpha \neq 0$. The critical value of the polynomial at the $A_2$ \ point is equal to $\pm 3^{\frac{15}{2}}2^{-8} \alpha^2$, where $\pm$ is the sign of $A$. This value is equal to $3^{\frac{15}{2}}2^{-8}$ if and only if this sign is $+$ and $\alpha$ equals $1$ or $-1$.
\hfill $\Box$
\medskip

\noindent
{\it Proof of Proposition \ref{estx1}.} The deduction of this proof from Lemma \ref{lemx1} 
essentially repeats the proof of Proposition \ref{estx+} with the following simplification. In this case, the neighborhoods of the point ${\bf \tilde T}(x^3y+\sqrt{3}x^2y^2+x y^3+x^3)$
 in the parameter space of the deformation (\ref{vers0}) have no nontrivial automorphisms that commute with the Lyashko--Looijenga covering (see \cite{vasX9}, \S~4.5). Therefore, the statement of Proposition \ref{estx1} follows immediately from the fact that all four polynomials (\ref{ea1}) belong to the same orbit of the group $G_1$. \hfill $\Box$
\medskip

\noindent
{\it Proof of Theorem \ref{theX91}.} Proposition \ref{estx1} implies that any virtual component of the formal graph of type $X_9^1$ is associated with at most four real connected components mapped to each other by the group $G_1$. Below we prove that six of the virtual components described in Proposition \ref{proX91} are associated with single real connected component each, while each of the four remaining virtual components is associated with a pair of real connected components that are mapped into each other by the central symmetry in $\R^2$. 

Pictures \ $(a)$, \ $(b)$, \ and \ $(c)$ \ of Fig.~\ref{x91} are realized by the polynomials
$$x y (x^2 + y^2) -t, $$
$$ x y(x^2+y^2) - t (x^2+y^2) + \alpha t^2, \quad \alpha \in (0,1/2), $$
and
$$x y (x^2+y^2 - t) - \alpha t^2, \quad \alpha \in (0, 1/8), $$
for any $t>0$.
All of these polynomials are invariant under the action of the group $G_1$. The polynomials obtained from them by the involution (\ref{invol}) and representing other three connected components of the discriminant complement are also invariant under this action.

Picture \ $(d)$ \ (respectively, \ $(e)$) is realized by the polynomials 
\begin{equation}
{\bf \tilde T}\left((x^2+y^2-t^2)(x+ 0.9 t)(y+ 0.9 t)+ \alpha t^4 \right), \quad \alpha \in (0, 0.003), 
\label{eq73}
\end{equation}
with positive (respectively, negative) $t$. These two pictures are topologically equivalent. However, it easily follows from Bezout's theorem that these polynomials represent different connected components of the discriminant complement
in the space $\Phi_1,$ see \cite{vasX9}.

Each polynomial (\ref{eq73}) is invariant under the transformation $(x,y) \mapsto (y,x)$. The transformation $(x,y) \mapsto (-x, -y)$ maps each such polynomial to the analogous polynomial with the opposite value of parameter \ $t$. In particular, it permutes the connected components of the discriminant complement that are described by pictures \ $(d)$ \ and \ $(e)$. 

Pictures \ $(f)$ \ and \ $(g)$ \ of Fig.~\ref{x91} are realized by polynomials
\begin{equation}
\label{tau2}
 x y \left((x-2t)^2 + (y-2t)^2 - t^2\right) - t^4 =0
\end{equation}
with positive and negative values of \ $t$, respectively. Each of these polynomials is invariant under the transformation $(x,y) \mapsto (y,x)$. The transformation $(x,y) \mapsto (-x, -y)$ maps each polynomial (\ref{tau2}) to the analogous polynomial with the opposite value of the parameter \ $t$, thus it permutes two connected components represented by these polynomials.

The same considerations hold for components obtained from these by the involution (\ref{invol}). This concludes the proof of Theorem \ref{theX91}. \hfill $\Box$

\subsection{Components for the $X_9^2$ singularities}
\label{sec92}

The group $G_0$ (see \S~\ref{sec90}) acts on the space $\Phi_2$ of all polynomials (\ref{vers0}) with $f_A= x^4 + 2A x^2y^2 + y^4,$ $A<-1$, preserving the coefficient $A$ and permuting the connected components of the discriminant complement. All these components belonging to the same orbit of this action are associated with the same virtual component of the formal graph of type $X_9^2$. 

In addition, the operation 
\begin{equation}
\label{ze}
f(x,y) \mapsto -f\left(\frac{x+y}{\sqrt{2}}, \frac{x-y}{\sqrt{2}}\right)
\end{equation}
acts on the space $\Phi_2$. This operation always alters the real connected components of the discriminant complement, as well as the associated virtual components of type $\Phi_2$. Indeed, it maps any generic polynomial with value \ $I$ \ of the invariant \ Ind \ to a polynomial with value \ $-3-I$. This operation is of order eight. All even iterations of it act as rotations by angles of $\frac{\pi k}{2}$ and thus belong to the group $G_0$. The entire set of connected components of the discriminant complement in the space $\Phi_2$ splits into the orbits of this operation. 

\unitlength 0.75mm
\linethickness{0.4pt}
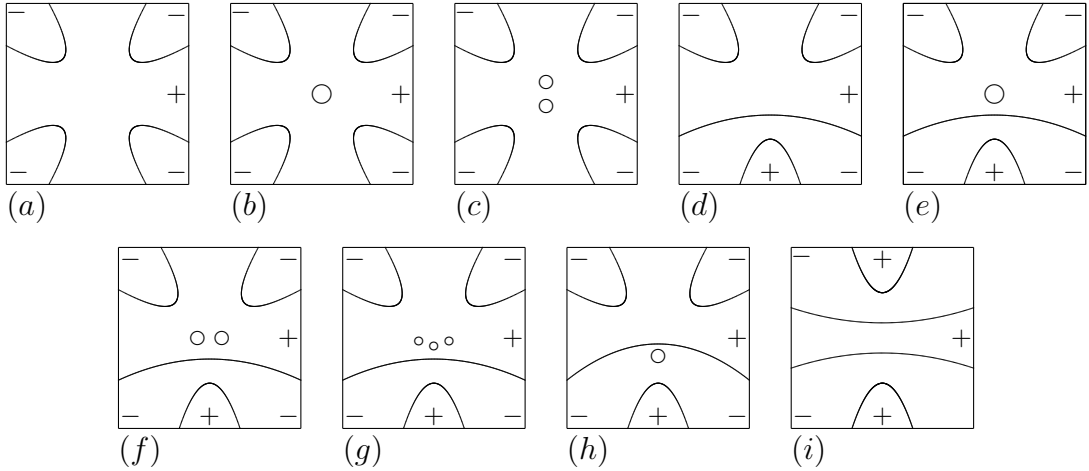
\begin{figure}
\begin{center}
\begin{picture}(30.00,35.00)
\put(0,5){\line(1,0){30}}
\put(0,5){\line(0,1){30}}
\put(0,35){\line(1,0){30}}
\put(30,5){\line(0,1){30}}
\bezier{200}(0,12)(15,20)(7,5)
\bezier{200}(30,12)(15,20)(23,5)
\bezier{200}(0,28)(15,20)(7,35)
\bezier{200}(30,28)(15,20)(23,35)
\put(2,7){\makebox(0,0)[cc]{\tiny $-$}}
\put(1.7,33.5){\makebox(0,0)[cc]{\tiny $-$}}
\put(28,7){\makebox(0,0)[cc]{\tiny $-$}}
\put(28,33){\makebox(0,0)[cc]{\tiny $-$}}
\put(28,20){\makebox(0,0)[cc]{\tiny $+$}}
 \put(0,0){$(a)$}
\end{picture} \quad 
\begin{picture}(30.00,35.00)
\put(0,5){\line(1,0){30}}
\put(0,5){\line(0,1){30}}
\put(0,35){\line(1,0){30}}
\put(30,5){\line(0,1){30}}
\bezier{200}(0,12)(15,20)(7,5)
\bezier{200}(30,12)(15,20)(23,5)
\bezier{200}(0,28)(15,20)(7,35)
\bezier{200}(30,28)(15,20)(23,35)
\put(2,7){\makebox(0,0)[cc]{\tiny $-$}}
\put(1.7,33.5){\makebox(0,0)[cc]{\tiny $-$}}
\put(28,7){\makebox(0,0)[cc]{\tiny $-$}}
\put(28,33){\makebox(0,0)[cc]{\tiny $-$}}
\put(28,20){\makebox(0,0)[cc]{\tiny $+$}}
\put(15,20){\circle{3}}
 \put(0,0){$(b)$}
\end{picture} \quad 
\begin{picture}(30.00,35.00)
\put(0,5){\line(1,0){30}}
\put(0,5){\line(0,1){30}}
\put(0,35){\line(1,0){30}}
\put(30,5){\line(0,1){30}}
\bezier{200}(0,12)(15,20)(7,5)
\bezier{200}(30,12)(15,20)(23,5)
\bezier{200}(0,28)(15,20)(7,35)
\bezier{200}(30,28)(15,20)(23,35)
\put(2,7){\makebox(0,0)[cc]{\tiny $-$}}
\put(1.7,33.5){\makebox(0,0)[cc]{\tiny $-$}}
\put(28,7){\makebox(0,0)[cc]{\tiny $-$}}
\put(28,33){\makebox(0,0)[cc]{\tiny $-$}}
\put(28,20){\makebox(0,0)[cc]{\tiny $+$}}
\put(15,22){\circle{2}}
\put(15,18){\circle{2}}
 \put(0,0){$(c)$}
\end{picture} \quad
\begin{picture}(30.00,40.00)
\put(0,5){\line(1,0){30}}
\put(0,5){\line(0,1){30}}
\put(0,35){\line(1,0){30}}
\put(30,5){\line(0,1){30}}
\bezier{200}(0,28)(15,20)(7,35)
\bezier{200}(30,28)(15,20)(23,35)
\bezier{200}(10,5)(15,20)(20,5)
\bezier{300}(0,13)(15,20)(30,13)
\put(2,33){\makebox(0,0)[cc]{\tiny $-$}}
\put(28,33){\makebox(0,0)[cc]{\tiny $-$}}
\put(2,7){\makebox(0,0)[cc]{\tiny $-$}}
\put(28,7){\makebox(0,0)[cc]{\tiny $-$}}
\put(15,7){\makebox(0,0)[cc]{\tiny $+$}}
\put(28,20){\makebox(0,0)[cc]{\tiny $+$}}
 \put(0,0){$(d)$ }
\end{picture} \quad 
\begin{picture}(30.00,35.00)
\put(0,5){\line(1,0){30}}
\put(0,5){\line(0,1){30}}
\put(0,35){\line(1,0){30}}
\put(30,5){\line(0,1){30}}
\put(0,5){\line(1,0){30}}
\put(0,5){\line(0,1){30}}
\put(0,35){\line(1,0){30}}
\put(30,5){\line(0,1){30}}
\bezier{200}(0,28)(15,20)(7,35)
\bezier{200}(30,28)(15,20)(23,35)
\bezier{200}(10,5)(15,20)(20,5)
\bezier{300}(0,13)(15,20)(30,13)
\put(2,33){\makebox(0,0)[cc]{\tiny $-$}}
\put(28,33){\makebox(0,0)[cc]{\tiny $-$}}
\put(2,7){\makebox(0,0)[cc]{\tiny $-$}}
\put(28,7){\makebox(0,0)[cc]{\tiny $-$}}
\put(15,7){\makebox(0,0)[cc]{\tiny $+$}}
\put(28,20){\makebox(0,0)[cc]{\tiny $+$}}
\put(15,20){\circle{3}}
 \put(0,0){$(e)$}
\end{picture} 

\begin{picture}(30.00,35.00)
\put(0,5){\line(1,0){30}}
\put(0,5){\line(0,1){30}}
\put(0,35){\line(1,0){30}}
\put(30,5){\line(0,1){30}}
\bezier{200}(0,28)(15,20)(7,35)
\bezier{200}(30,28)(15,20)(23,35)
\bezier{200}(10,5)(15,20)(20,5)
\bezier{300}(0,13)(15,20)(30,13)
\put(2,33){\makebox(0,0)[cc]{\tiny $-$}}
\put(28,33){\makebox(0,0)[cc]{\tiny $-$}}
\put(2,7){\makebox(0,0)[cc]{\tiny $-$}}
\put(28,7){\makebox(0,0)[cc]{\tiny $-$}}
\put(15,7){\makebox(0,0)[cc]{\tiny $+$}}
\put(28,20){\makebox(0,0)[cc]{\tiny $+$}}
\put(13,20){\circle{2}}
\put(17,20){\circle{2}}
 \put(0,0){$(f)$}
\end{picture} \quad
\begin{picture}(30.00,40.00)
\put(0,5){\line(1,0){30}}
\put(0,5){\line(0,1){30}}
\put(0,35){\line(1,0){30}}
\put(30,5){\line(0,1){30}}
\bezier{200}(0,28)(15,20)(7,35)
\bezier{200}(30,28)(15,20)(23,35)
\bezier{200}(10,5)(15,20)(20,5)
\bezier{300}(0,13)(15,20)(30,13)
\put(2,33){\makebox(0,0)[cc]{\tiny $-$}}
\put(28,33){\makebox(0,0)[cc]{\tiny $-$}}
\put(2,7){\makebox(0,0)[cc]{\tiny $-$}}
\put(28,7){\makebox(0,0)[cc]{\tiny $-$}}
\put(15,7){\makebox(0,0)[cc]{\tiny $+$}}
\put(28,20){\makebox(0,0)[cc]{\tiny $+$}}
\put(12.5,19.5){\circle{1.2}}
\put(15,18.7){\circle{1.2}}
\put(17.5,19.5){\circle{1.2}}
 \put(0,0){$(g)$}
\end{picture} \quad 
\begin{picture}(30.00,37.00)
\put(0,5){\line(1,0){30}}
\put(0,5){\line(0,1){30}}
\put(0,35){\line(1,0){30}}
\put(30,5){\line(0,1){30}}
\bezier{200}(0,28)(15,20)(7,35)
\bezier{200}(30,28)(15,20)(23,35)
\bezier{200}(10,5)(15,20)(20,5)
\bezier{300}(0,13)(15,25)(30,13)
\put(2,33){\makebox(0,0)[cc]{\tiny $-$}}
\put(28,33){\makebox(0,0)[cc]{\tiny $-$}}
\put(2,7){\makebox(0,0)[cc]{\tiny $-$}}
\put(28,7){\makebox(0,0)[cc]{\tiny $-$}}
\put(15,7){\makebox(0,0)[cc]{\tiny $+$}}
\put(28,20){\makebox(0,0)[cc]{\tiny $+$}}
\put(15,17){\circle{2}}
 \put(0,0){$(h)$}
\end{picture} \quad 
\begin{picture}(30.00,37.00)
\put(0,5){\line(1,0){30}}
\put(0,5){\line(0,1){30}}
\put(0,35){\line(1,0){30}}
\put(30,5){\line(0,1){30}}
\bezier{200}(0,15)(15,20)(30,15)
\bezier{200}(0,25)(15,20)(30,25)
\bezier{300}(10,5)(15,20)(20,5)
\bezier{300}(10,35)(15,20)(20,35)
\put(28,20){\makebox(0,0)[cc]{\tiny $+$}}
\put(15,7){\makebox(0,0)[cc]{\tiny $+$}}
\put(15,33){\makebox(0,0)[cc]{\tiny $+$}}
\put(2,7){\makebox(0,0)[cc]{\tiny $-$}}
\put(1.7,33.5){\makebox(0,0)[cc]{\tiny $-$}}
 \put(0,0){$(i)$}
\end{picture}
\caption{Non-discriminant perturbations of $X_9^2$}
\label{X92}
\end{center}
\end{figure}

\begin{proposition}[see \cite{para}]
\label{eee}
There are exactly 18 virtual components of class $\Phi_2$. The shapes of the zero-level sets of the polynomials representing nine of these virtual components are shown in Fig.~\ref{X92} $($see also Fig.~33 in \cite{viro2}$)$, and the polynomials representing the remaining nine virtual components are obtained from these by the operation $($\ref{ze}$)$.
\end{proposition}

\begin{theorem}
\label{thmX2}
There are exactly 52 connected components of the discriminant complement in the space $\Phi_2$. 

Specifically, the virtual components represented by pictures \ $(a)$ \ and \ $(b)$ \ of Fig.~\ref{X92} are associated with single real connected components of the discriminant complement;
the virtual components represented by pictures \ $(c)$ \ and \ $(i)$ \ are associated with pairs of such components $($whose pictures are obtained from each other by rotating the plane by the angle $\frac{\pi}{2})$; each of the virtual components represented by pictures \ $(d)$, \ $(e)$, \ $(f)$, \ $(g)$, \ and \ $(h)$ \ is associated with four real connected components of the discriminant complement, whose pictures are obtained from each other by rotating the plane by the angles \ $\frac{\pi}{2}$, \ $\pi$, \ and \ $\frac{3\pi}{2}$. \ The same is true for the components obtained from these nine by the operation $($\ref{ze}$)$.
\end{theorem} 

\begin{remark} \rm
Picture ($c$) and the result of its rotation by the angle \ $\frac{\pi}{2}$ \ are topologically equivalent, but the connected components of the discriminant complement realized by polynomials with these shapes are different.
\end{remark}

Theorem \ref{thmX2} was conjectured in \cite{para}, and the difference of all the connected components mentioned there was proven. It remains to prove that there are no additional components, that is, the virtual components listed in Proposition \ref{eee} are associated with no more real components than indicated in this theorem.

\begin{proposition}
\label{estx2}
If two generic polynomials $f$ and $\tilde f$ of class $\Phi_2$ are associated with the same virtual function, then their connected components of the space of generic polynomials of the form $($\ref{vers0}$)$ are mapped to each other by an element of the group $ G_0$. 
\end{proposition}

\begin{lemma}
\label{prost8}
There are exactly eight polynomials $\R^2 \to \R$ of class $\Phi_2$ that have a critical point of class $D_6$ with critical value 0 and a critical point of class $A_3$ with critical value $\frac{-1}{324}$. These eight polynomials form an orbit of the natural action of the group $G_0$ on $\Phi_2$.
\end{lemma}

\begin{lemma}
\label{le833}
Each polynomial of degree four on $\R^2$ having a critical point of class $D_6$ with zero critical value and a critical point of class $A_3$ can be reduced by an affine transformation of $\R^2$ to the form 
\begin{equation}
\label{eqqt}
\pm \left(t \left(\tilde x^2 \tilde y + \frac{1}{3}\tilde x^3\right) + \frac{3}{4}\tilde x^4 + 3 \tilde x^3 \tilde y + \frac{3}{2}\tilde x^2 \tilde y^2-\tilde x \tilde y^3\right) \ , \quad t > 0 \ .
\end{equation}
In particular, the j-invariant of its principal homogeneous part is equal to 5.
\end{lemma}

\noindent
{\it Proof of Lemma \ref{le833}.}
This lemma follows from Proposition 22 of \cite{vasX9} by the substitution $\tilde x = |a|^{\frac{3}{4}} x$, $\tilde y= \mbox{sign}(a) |a|^{\frac{-1}{4}} y$, $t=|a|^{\frac{-5}{4}}$, and the sign $\pm$ in (\ref{eqqt}) equal to $\mbox{sign}(a)$. \hfill $\Box$ 

\begin{corollary}
If a polynomial of degree four in $\R^2$ has a critical point of class $D_6$ with zero critical value and a point of class $A_3$ with critical value $\frac{-1}{324}$, then it has the form 
\begin{equation}
\label{eqqt7}
\tilde x^2 \tilde y + \frac{1}{3}\tilde x^3 + \frac{3}{4}\tilde x^4 + 3 \tilde x^3 \tilde y + \frac{3}{2}\tilde x^2 \tilde y^2- \tilde x \tilde y^3
\end{equation}
in some affine coordinates. 
\end{corollary}

\noindent
{\it Proof.} Of all the polynomials (\ref{eqqt}), only the polynomial (\ref{eqqt7}) has critical value $\frac{-1}{324}$. \hfill $\Box$

\begin{corollary}
If a polynomial $f(x, y)$ of degree four with principal homogeneous part
\begin{equation}
\label{fo}
x^4 + 2A x^2y^2 +y^4, \quad A < -1, 
\end{equation}
 has critical points of classes $D_6$ and $A_3$, then $A$ is a root of the equation 
\begin{equation}
\label{135}
\frac{(A^2+3)^3}{27 (A^2-1)^2}=5.
\end{equation}
\end{corollary}

\noindent
{\it Proof.} This follows from the well-known formula for $j$-invariant of quartic forms, see e.g. \cite{Mukai}. \hfill $\Box$
\medskip

\noindent
{\it Proof of Lemma \ref{prost8}.}
There are only two roots of equation (\ref{135}) in the domain $\{A<-1\}$: they are equal approximately to \ $-1.395$ \ and \ $-11.118$. In the first (respectively, in the second) case, any linear transformation of $\R^2$ that moves the zero lines of the principal homogeneous part of (\ref{eqqt7}) to those of (\ref{fo}) sends the domain of positive values of the former quartic form to the domain of positive (respectively, negative) values of the latter. Therefore, only the polynomials of class $\Phi_2$ with $A \approx -1.395$ can have the form (\ref{eqqt7}) in some affine coordinates $(\tilde x, \tilde y)$. The line $\{\tilde x=0\}$ should be parallel to one of the four lines of zeros of the principal part (\ref{fo}). For any choice of this line, of an orientation of this line, and of a point $O \in \R^2$, there is exactly one pair $(f, \{\tilde x, \tilde y\})$ consisting of a polynomial function of degree four with principal homogeneous part (\ref{fo}) and an affine coordinate system $\{\tilde x, \tilde y\}$ centered at the point $O$ such that the direction $\frac{\partial}{\partial \tilde y}$ in these coordinates coincides with the direction of the chosen oriented line, and the polynomial $f$ has the form (\ref{eqqt7}) in these coordinates. However, for any choice of an oriented line of zeros of the polynomial (\ref{fo}) with $A$ equal to the root $\approx -1.395$ of the equation (\ref{135}) only one choice of the point $O$ provides in this way a function $f$ satisfying the condition $ f \in \Phi_2$ (i.e., the vanishing of the monomials $x^3$ and $y^3$). Thus, this function is uniquely determined by a choice of an oriented line of zeros of (\ref{fo}). All eight polynomials $f$ defined by these choices form an orbit of the group $G_0$. \hfill $\Box$ 

\medskip
The rest of the proof of Proposition \ref{estx2} repeats the proof of Proposition \ref{estx+} with the simplification indicated in the proof of Proposition \ref{estx1}. Let $f_0$ be one of eight polynomials mentioned in in Lemma \ref{prost8}. The space $\Phi_2$ does not have automorphisms that commute with the Lyashko--Looijenga map and map this polynomial to itself, see the proof of Proposition 
3 in \S~5.4 of \cite{vasX9}. Let $f$ and $\tilde f$ be two generic polynomials of class $\Phi_2$ with the same associated virtual function and sets of critical values and paths. Let $I$ be a generic path in $\Phi_2$ connecting $f$ with $f_0$, and $\tilde I$ be a path starting at $\tilde f$ and repeating all the critical values and surgeries of the path $I$. Then the endpoint $\tilde f_0$ of this path is well-defined and is also one of these eight points, and the entire path $\tilde I$ is mapped to $I$ by the element of the group $G_0$ that maps the point $\tilde f_0$ to $f_0$. 
\hfill $\Box$
\medskip

\noindent
{\it Proof of Theorem \ref{thmX2}.}
Pictures \ ($a$) \ and \ ($b$) \ of Fig.~\ref{X92} are realized, respectively, by the
polynomials (\ref{x2a}) and (\ref{x2b}) with arbitrary $t>0$. 
\begin{equation}
\label{x2a}
x^4-6x^2y^2+y^4+ t 
\end{equation}
\begin{equation}
\label{x2b}
x^4-6x^2y^2+y^4 -2t(x^2+y^2)+ \frac{1}{2}t^2
\end{equation}
These polynomials are invariant under the entire group $G_0$. Therefore, their virtual components are associated with single real connected components of the discriminant complement.
 
Pictures \ ($c$) \ and \ ($i$) \ are realized respectively by polynomials (\ref{x2c}) and (\ref{x2i}) with arbitrary $t>0$. 
\begin{equation}
\label{x2c}
(x^2-3t)^2-6(x^2-3t)(y^2-t/2)+ (y^2-t/2)^2+t^4
\end{equation}
\begin{equation}
\label{x2i}
x^4-6x^2(y^2+t)+(y^2+t)^2-t^2/2 
\end{equation}
These polynomials are invariant under the reflections $(x,y) \mapsto (x,-y)$ and $(x,y) \mapsto (-x,y),$ which generate a subgroup of index two in $G_0$. Thus, in each of these two cases there are no more than two connected components of the discriminant complement associated with the corresponding virtual component.

Pictures \ ($d$), \ ($e$), \ ($f$), \ ($g$), \ and \ ($h$) \ are realized by the equations (\ref{x2d}), (\ref{x2e}), (\ref{x2f}), (\ref{x2g}), and (\ref{x2h}) for sufficiently small $\varepsilon>0$ and arbitrary $t>0$. 
\begin{equation}
\label{x2d}
{\bf T}\left(x^4-6x^2y^2+y^4+t x^3+\varepsilon t^3 x\right) 
\end{equation}
\begin{equation}
\label{x2e}
{\bf T}\left(x^4-6x^2y^2+y^4+t x^3- \varepsilon t^3 x\right)
\end{equation}
\begin{equation}
\label{x2f}
{\bf T}\left(x^4-6x^2y^2+t(x^3-\varepsilon^4 t^2 x) +(y^2-\varepsilon^3t^2)^2\right)
\end{equation}
\begin{equation}
\label{x2g}
{\bf T}\left(x^4-6x^2y^2+y^4+ t x^3-6\varepsilon t x y^2 - 2 \varepsilon^4 t^3 x+ 8 \varepsilon^3 t^2 y^2\right)
\end{equation}
\begin{equation}
\label{x2h}
{\bf T}\left((x+0.7t)(x+1.3t)(x+6t)(x+7t)-6x^2y^2+y^4\right)
\end{equation}
These five polynomials are all constructed as follows.
First we take the polynomial 
$$x^4-6x^2y^2+y^4 + x^3,$$ whose zero-level set is topologically equivalent to the picture (d), but has a singular point of class $E_6$ at its center. Then, we perturb this singular point by adding very small terms of lower degrees. These terms do not change the topology of the curve far away from this singular point, but provide the standard dissolutions of the $E_6$ singularity, see Figure 24 in \cite{viro2}. All the resulting polynomials (\ref{x2d})--(\ref{x2h}) 
are invariant under the reflection in the axis $y=0$. Thus, in each of these five cases, the virtual component is associated with no more than four real connected components of the discriminant complement. These four components differ for the topological reasons, see \cite{para}.

For the remaining nine virtual components, the convenient realizations of the corresponding real connected components can be obtained from the polynomials (\ref{x2a})--(\ref{x2h}) via the operation (\ref{ze}). \hfill $\Box$

\section{Components of discriminant complements for $J_{10}$ singularities}

\subsection{}
The canonical versal deformations of the singularities of the classes $J_{10}^1$ and $J_{10}^3$ consist, respectively, of polynomials
\begin{multline}
(x-A y^2)(x^2+y^4)+ \lambda_1 + \lambda_2 y + \lambda_3 y^2 + \lambda_4 y^3 + \lambda_5 y^4 + \\
\lambda_6 x +
\lambda_7 x y + \lambda_8 x y^2 + \lambda_9 (x^2+3y^4)y , \qquad A \in (-\infty, +\infty), \label{vers1} \end{multline}
and
\begin{multline}
(x-A y^2)(x^2-y^4)+ \lambda_1 + \lambda_2 y + \lambda_3 y^2 + \lambda_4 y^3 + \lambda_5 y^4 + \\
\lambda_6 x + \lambda_7 x y + \lambda_8 x y^2 + \lambda_9 (x^2-3y^4)y , \qquad A \in (-1,1), \label{vers3}
\end{multline} 
 with ten parameters $\lambda_1, \dots, \lambda_9,$ and $A$, cf. \cite{Jaw2}, \cite{vasJ10}.

\begin{definition} \rm
\label{genJ}
A polynomial $f: ({\mathbb C}^2, {\mathbb R}^2) \to ({\mathbb C}, {\mathbb R})$ of the form (\ref{vers1}) or (\ref{vers3}) is {\em generic} if it has only Morse critical points in ${\mathbb C}^2$, all of its critical values are different and not equal to 0, and the complex zero-level set of the principal weighted homogeneous part of $f$ (of degree six with weights $\deg x =2, \deg y =1$)
consists of three different complex parabolas or two parabolas and one line in ${\mathbb C}^2$.

Denote by $\Omega_1$ and $\Omega_3$ the spaces of all polynomials of the form (\ref{vers1}) and (\ref{vers3}), respectively. The virtual functions associated with generic polynomials of these types, the formal graphs of types $J_{10}^1$ and $J_{10}^3$, and their virtual components are defined in the same way as for the case $X_9$, except that the sets of vanishing cycles consist of ten elements.
\end{definition}

The following analog of Propositions \ref{protrivi} and \ref{propmain} follows from the same considerations as these statements. 

\begin{proposition}
\label{protriviJ}
For any class of polynomials $\Omega_1$ or \ $\Omega_3$,
1$)$ if two generic polynomials, $f$ and $\tilde f$, of this class belong to the same connected component of the complement of the discriminant, then the virtual components $S(f)$ and $S(\tilde f)$ are the same. 

2$)$ every virtual function of this class is associated with a generic real polynomial of the same class;

3$)$ for any generic polynomial $f$ of this class, every virtual function from the virtual component $S(f)$ is associated with some generic real polynomial $\tilde f$ from the same connected component of the complement of the discriminant as $f$. \hfill $\Box$ \end{proposition}

\begin{proposition}[see \cite{vasJ10}, Proposition 9]
\label{estj1}
If two generic polynomials $f$, $\tilde f$ of class $\Omega_1$ or $\Omega_3$ are associated with the same virtual function, then their connected components of the space of generic polynomials either are the same or are mapped to each other by the involution 
\begin{equation}
\label{invol4}
f(x,y) \leftrightarrow f(x,-y) \ .
\end{equation}
In particular, any virtual component of class $\Omega_1$ or $\Omega_3$ is associated with at most two real connected components of the discriminant complement in $\Omega_*$. \hfill $\Box$
\end{proposition}

In the following pictures (\ref{J101})--(\ref{J38}), it is assumed that the $x$ coordinate axis is directed upwards, and the horizontal line is the axis $\{x=0\}$.

\subsection{Components for $J_{10}^1$ singularities }

\begin{figure}
\unitlength 1mm
\begin{picture}(20.00,25.00)
\put(0,5){\line(1,0){20}}
\put(0,5){\line(0,1){20}}
\put(0,25){\line(1,0){20}}
\put(20,5){\line(0,1){20}}
\put(0,15){\line(1,0){20}}
\put(8.5,22){\footnotesize $+$}
\put(8.5,6){\footnotesize $-$}
 \put(1.5,0){$(a)$}
\end{picture} \ \
\begin{picture}(20.00,25.00)
\put(0,5){\line(1,0){20}}
\put(0,5){\line(0,1){20}}
\put(0,25){\line(1,0){20}}
\put(20,5){\line(0,1){20}}
\put(0,15){\line(1,0){20}}
\put(8.5,22){\footnotesize $+$}
\put(15.5,6){\footnotesize $-$}
\put(10,10){\circle{3}}
 \put(1.5,0){$(b)$}
\end{picture} \ \
\begin{picture}(20.00,25.00)
\put(0,5){\line(1,0){20}}
\put(0,5){\line(0,1){20}}
\put(0,25){\line(1,0){20}}
\put(20,5){\line(0,1){20}}
\put(0,15){\line(1,0){20}}
\put(7,10){\circle{3}}
\put(13,10){\circle{3}}
\put(8.5,22){\footnotesize $+$}
\put(8.5,6){\footnotesize $-$}
 \put(1.5,0){$(c)$}
\end{picture} \ \
\begin{picture}(20.00,25.00)
\put(0,5){\line(1,0){20}}
\put(0,5){\line(0,1){20}}
\put(0,25){\line(1,0){20}}
\put(20,5){\line(0,1){20}}
\put(0,15){\line(1,0){20}}
\put(5,8.5){\circle{2.5}}
\put(15,8.5){\circle{2.5}}
\put(10,11){\circle{2.5}}
\put(8.5,22){\footnotesize $+$}
\put(8.5,6){\footnotesize $-$}
 \put(1.5,0){$(d)$}
\end{picture} \ \
\begin{picture}(20.00,25.00)
\put(0,5){\line(1,0){20}}
\put(0,5){\line(0,1){20}}
\put(0,25){\line(1,0){20}}
\put(20,5){\line(0,1){20}}
\put(0,15){\line(1,0){20}}
\put(3,8){\circle{2.5}}
\put(17,8){\circle{2.5}}
\put(7,10){\circle{2.5}}
\put(13,10){\circle{2.5}}
\put(8.5,22){\footnotesize $+$}
\put(8.5,6){\footnotesize $-$}
 \put(1.5,0){$(e)$}
\end{picture} \ \
\begin{picture}(20.00,27.00)
\put(0,5){\line(1,0){20}}
\put(0,5){\line(0,1){20}}
\put(0,25){\line(1,0){20}}
\put(20,5){\line(0,1){20}}
\put(0,15){\line(1,0){20}}
\put(6,10){\circle{3}}
\put(14,19){\circle{3}}
\put(15.5,22){\footnotesize $+$}
\put(15.5,6){\footnotesize $-$}
 \put(1.5,0){$(f)$}
\end{picture}

\begin{picture}(20.00,25.00)
\put(0,0){$\phantom{A}$}
\end{picture} \ \
\begin{picture}(20.00,25.00)
\put(0,5){\line(1,0){20}}
\put(0,5){\line(0,1){20}}
\put(0,25){\line(1,0){20}}
\put(20,5){\line(0,1){20}}
\put(0,15){\line(1,0){20}}
\put(15.5,22){\footnotesize $+$}
\put(8.5,6){\footnotesize $-$}
\put(10,20){\circle{3}}
 \put(1.5,0){($\tilde b$)}
\end{picture} \ \ 
\begin{picture}(20.00,25.00)
\put(0,5){\line(1,0){20}}
\put(0,5){\line(0,1){20}}
\put(0,25){\line(1,0){20}}
\put(20,5){\line(0,1){20}}
\put(0,15){\line(1,0){20}}
\put(7,20){\circle{3}}
\put(13,20){\circle{3}}
\put(8.5,22){\footnotesize $+$}
\put(8.5,6){\footnotesize $-$}
 \put(1.5,0){($\tilde c$)}
\end{picture} \ \
\begin{picture}(20.00,25.00)
\put(0,5){\line(1,0){20}}
\put(0,5){\line(0,1){20}}
\put(0,25){\line(1,0){20}}
\put(20,5){\line(0,1){20}}
\put(0,15){\line(1,0){20}}
\put(5,21.5){\circle{2.5}}
\put(15,21.5){\circle{2.5}}
\put(10,19){\circle{2.5}}
\put(8.5,22){\footnotesize $+$}
\put(8.5,6){\footnotesize $-$}
 \put(1.5,0){($\tilde d$)}
\end{picture} \ \ 
\begin{picture}(20.00,29.00)
\put(0,5){\line(1,0){20}}
\put(0,5){\line(0,1){20}}
\put(0,25){\line(1,0){20}}
\put(20,5){\line(0,1){20}}
\put(0,15){\line(1,0){20}}
\put(3,22){\circle{2.5}}
\put(17,22){\circle{2.5}}
\put(7,20){\circle{2.5}}
\put(13,20){\circle{2.5}}
\put(8.5,22){\footnotesize $+$}
\put(8.5,6){\footnotesize $-$}
 \put(1.5,0){($\tilde e$)}
\end{picture} \ \
\begin{picture}(20.00,27.00)
\put(0,0){$\phantom{A}$}
\end{picture}
\caption{Non-discriminant perturbations of $J_{10}^1$ singularities}
\label{J101}
\end{figure}
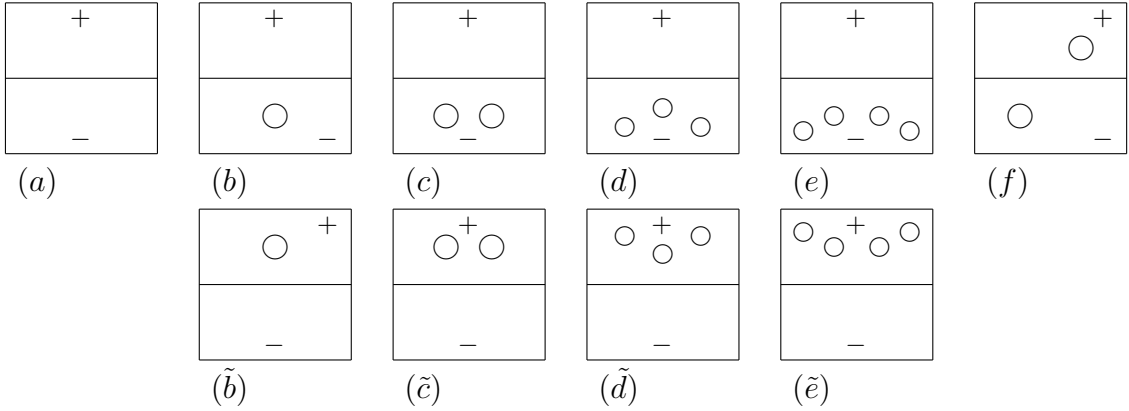

\begin{proposition}[see \S 4.2 in \cite{viro2}]
There are exactly ten topological types of zero-level sets of non-discriminant perturbations of $J_{10}^1$ singularities, see Fig.~\ref{J101}. \hfill $\Box$
\end{proposition}

\begin{proposition}[see \cite{para}]
\label{paraj101}
There are exactly ten virtual components of the formal graph of type $J_{10}^1$. These virtual components are realized by polynomials, the topological shapes of whose zero-level sets are shown on Fig.~\ref{J101}. \hfill $\Box$
\end{proposition}

\begin{theorem}
\label{thmj1}
There are exactly thirteen connected components of the discriminant complement in the space of all polynomials of class $\Omega_1$. 
Each of the virtual components realized by polynomials with topological types shown in pictures $(a)$, $(b)$, $(\tilde b)$, $(c)$, $(\tilde c)$, $(d)$, and $(\tilde d)$ of Fig.~\ref{J101} is associated with a single real connected component.
Each of the virtual components/topological types represented in pictures $(e)$, $(\tilde e),$ and $(f)$ is associated with two different components that are related by the involution $($\ref{invol4}$)$.
\end{theorem}

\noindent
{\it Proof.} The topological types shown in pictures \ ($a$), \ ($b$), \ ($c$), \ and \ ($d$) \ are represented respectively by the following polynomials: 
$$x(x^2 +y^4+t), $$
$$(x-2t) ((x+t)^2+ y^4-t^2),$$
$$(x-2t)((x+t)^2 + y^4-3t y^2 +2 t^2),$$
and
$$(x- y^2 + 2 t ) (x^2-2t x+ y^4-5t^2) - t^3/10 $$ 
with arbitrary $t>0$. These polynomials are all invariant under the involution (\ref{invol4}). Therefore, the corresponding virtual components are associated with single connected components of the discriminant complement.

The realizations of pictures \ $(\tilde b)$, \ $(\tilde c)$, \ and \ $(\tilde d)$ \ can be obtained from these realizations of pictures \ $(b),$ \ $(c)$, \ and \ $(d)$ \ by the involution
\begin{equation}
\label{invol94}
f(x,y) \mapsto -f(-x, y).
\end{equation}
All of them are also invariant under the involution (\ref{invol4}).

No polynomial $f(x,y)$ of class $\Omega_1$ that realizes the picture \ $(f)$ \ can belong to the same connected component as its image $f(x, -y)$ under the involution (\ref{invol4}). Indeed, the polynomials (\ref{vers1}) are of degree three in the variable \ $x$, so any vertical line $\{y=\mbox{const}\}$ can intersect the curve $\{f=0\}$ at most three times. Therefore, the ovals placed on the different sides of the infinite component of the zero-level set cannot pass one above the other. On the other hand, the involution (\ref{invol94}) moves the component of picture $(f)$ to the same component as the involution (\ref{invol4}).

Every polynomial $f(x,y) $ of class $\Omega_1$ that realizes the picture \ $(e)$ \ belongs to a different connected component than its mirror image $f(x, -y)$. Indeed, the four ovals defined by this polynomial are canonically ordered from the right to the left, and this ordering is continuous along the connected components of the discriminant complement. 
Let us take an arbitrary point inside each of these four ovals. There is a single curve with the equation \ $x = a y^3+b y^2 + c y + d$ \ that passes through these four points. The coefficient $a$ of this equation cannot be zero. Otherwise, the restriction of the polynomial $f$ to this curve would be a polynomial in the coordinate $y$ of degree at most six with at least eight roots. Therefore, the sign of the coefficient $a$ does not depend on the choice of the four points inside the ovals and remains constant along the connected component of the discriminant complement. Any polynomial $f(x,y)$ in the considered component and its mirror image $f(x,-y)$ define opposite signs.
 The same is true for polynomials that realize the picture \ $(\tilde e)$. \hfill $\Box$

\subsection{Components for $J_{10}^3$ singularities}

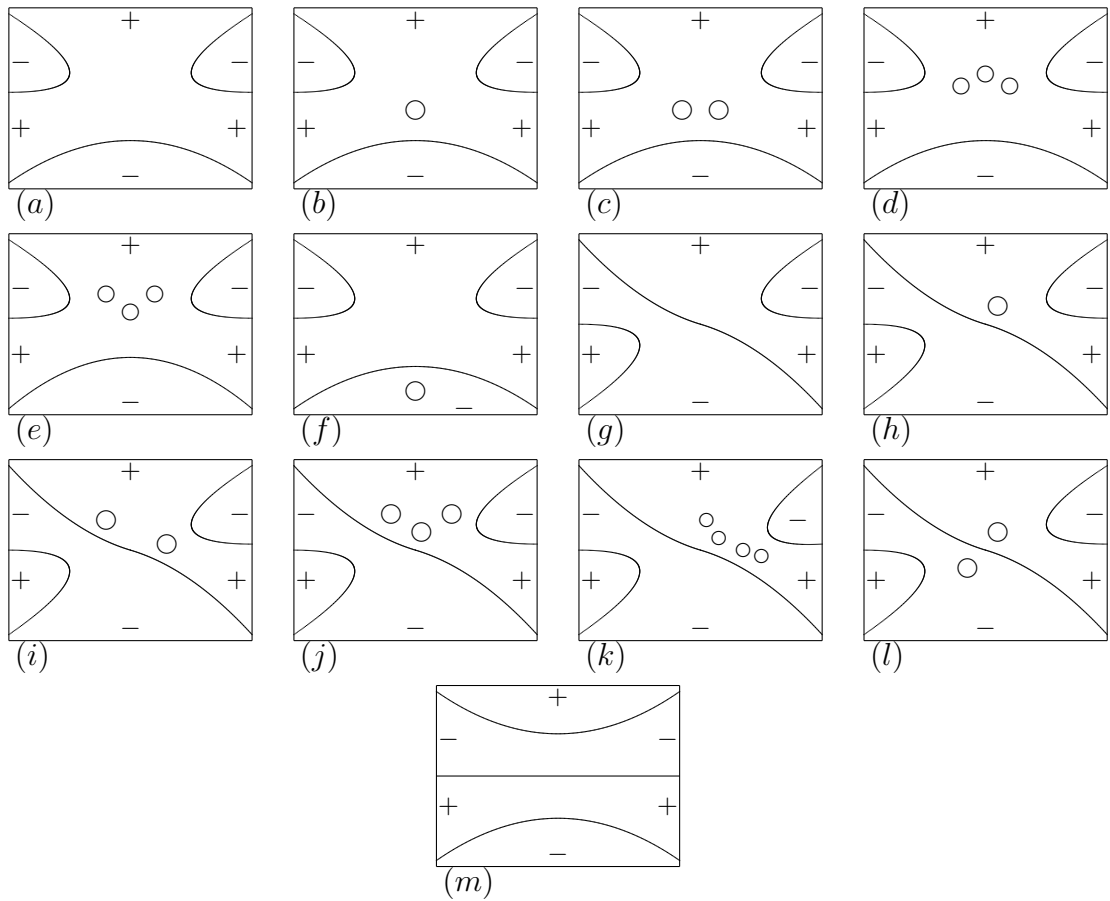
\begin{figure}
\unitlength 0.75 mm
\begin{center}
\begin{picture}(40.00,35.00)
\put(0,5){\line(1,0){40}}
\put(0,5){\line(0,1){30}}
\put(0,35){\line(1,0){40}}
\put(40,5){\line(0,1){30}}
\bezier{300}(0,6)(20,20)(40,6)
\bezier{200}(0,21)(20,21)(0,34)
\bezier{200}(40,21)(20,21)(40,34)
\put(2,26){\makebox(0,0)[cc]{\tiny $-$}} 
\put(20,7){\makebox(0,0)[cc]{\tiny $-$}} 
\put(38,26){\makebox(0,0)[cc]{\tiny $-$}}
\put(2,15){\makebox(0,0)[cc]{\tiny $+$}}
\put(37.5,15){\makebox(0,0)[cc]{\tiny $+$}}
\put(20,33){\makebox(0,0)[cc]{\tiny $+$}}
\put(1,1){$(a)$}
\end{picture} \quad
\begin{picture}(40.00,37.00)
\put(0,5){\line(1,0){40}}
\put(0,5){\line(0,1){30}}
\put(0,35){\line(1,0){40}}
\put(40,5){\line(0,1){30}}
\bezier{300}(0,6)(20,20)(40,6)
\bezier{200}(0,21)(20,21)(0,34)
\bezier{200}(40,21)(20,21)(40,34)
\put(20,18){\circle{3}}
\put(2,26){\makebox(0,0)[cc]{\tiny $-$}} 
\put(20,7){\makebox(0,0)[cc]{\tiny $-$}} 
\put(38,26){\makebox(0,0)[cc]{\tiny $-$}}
\put(2,15){\makebox(0,0)[cc]{\tiny $+$}}
\put(37.5,15){\makebox(0,0)[cc]{\tiny $+$}}
\put(20,33){\makebox(0,0)[cc]{\tiny $+$}}
\put(1,1){$(b)$}
\end{picture} \quad
\begin{picture}(40.00,35.00)
\put(0,5){\line(1,0){40}}
\put(0,5){\line(0,1){30}}
\put(0,35){\line(1,0){40}}
\put(40,5){\line(0,1){30}}
\bezier{300}(0,6)(20,20)(40,6)
\bezier{200}(0,21)(20,21)(0,34)
\bezier{200}(40,21)(20,21)(40,34)
\put(17,18){\circle{3}}
\put(23,18){\circle{3}}
\put(2,26){\makebox(0,0)[cc]{\tiny $-$}} 
\put(20,7){\makebox(0,0)[cc]{\tiny $-$}} 
\put(38,26){\makebox(0,0)[cc]{\tiny $-$}}
\put(2,15){\makebox(0,0)[cc]{\tiny $+$}}
\put(37.5,15){\makebox(0,0)[cc]{\tiny $+$}}
\put(20,33){\makebox(0,0)[cc]{\tiny $+$}}
\put(1,1){$(c)$} 
\end{picture} \quad
\begin{picture}(40.00,35.00)
\put(0,5){\line(1,0){40}}
\put(0,5){\line(0,1){30}}
\put(0,35){\line(1,0){40}}
\put(40,5){\line(0,1){30}}
\bezier{300}(0,6)(20,20)(40,6)
\bezier{200}(0,21)(20,21)(0,34)
\bezier{200}(40,21)(20,21)(40,34)
\put(16,22){\circle{2.5}}
\put(24,22){\circle{2.5}}
\put(20,24){\circle{2.5}}
\put(2,26){\makebox(0,0)[cc]{\tiny $-$}} 
\put(20,7){\makebox(0,0)[cc]{\tiny $-$}} 
\put(38,26){\makebox(0,0)[cc]{\tiny $-$}}
\put(2,15){\makebox(0,0)[cc]{\tiny $+$}}
\put(37.5,15){\makebox(0,0)[cc]{\tiny $+$}}
\put(20,33){\makebox(0,0)[cc]{\tiny $+$}}
\put(1,1){$(d)$}
\end{picture} \quad
\begin{picture}(40.00,37.00)
\put(0,5){\line(1,0){40}}
\put(0,5){\line(0,1){30}}
\put(0,35){\line(1,0){40}}
\put(40,5){\line(0,1){30}}
\put(16,25){\circle{2.5}}
\put(24,25){\circle{2.5}}
\put(20,22){\circle{2.5}}
\bezier{300}(0,6)(20,23)(40,6)
\bezier{200}(0,21)(20,21)(0,34)
\bezier{200}(40,21)(20,21)(40,34)
\put(2,26){\makebox(0,0)[cc]{\tiny $-$}} 
\put(20,7){\makebox(0,0)[cc]{\tiny $-$}} 
\put(38,26){\makebox(0,0)[cc]{\tiny $-$}}
\put(2,15){\makebox(0,0)[cc]{\tiny $+$}}
\put(37.5,15){\makebox(0,0)[cc]{\tiny $+$}}
\put(20,33){\makebox(0,0)[cc]{\tiny $+$}}
\put(1,1){$(e)$} 
\end{picture}
\quad
\begin{picture}(40.00,37.00)
\put(0,5){\line(1,0){40}}
\put(0,5){\line(0,1){30}}
\put(0,35){\line(1,0){40}}
\put(40,5){\line(0,1){30}}
\bezier{300}(0,6)(20,20)(40,6)
\bezier{200}(0,21)(20,21)(0,34)
\bezier{200}(40,21)(20,21)(40,34)
\put(20,9){\circle{3}}
\put(2,26){\makebox(0,0)[cc]{\tiny $-$}} 
\put(28,6){\makebox(0,0)[cc]{\tiny $-$}} 
\put(38,26){\makebox(0,0)[cc]{\tiny $-$}}
\put(2,15){\makebox(0,0)[cc]{\tiny $+$}}
\put(37.5,15){\makebox(0,0)[cc]{\tiny $+$}}
\put(20,33){\makebox(0,0)[cc]{\tiny $+$}}
\put(1,1){$(f)$}
\end{picture} \quad
\begin{picture}(40.00,35.00)
\put(0,5){\line(1,0){40}}
\put(0,5){\line(0,1){30}}
\put(0,35){\line(1,0){40}}
\put(40,5){\line(0,1){30}}
\bezier{200}(0,20)(20,20)(0,6)
\bezier{200}(40,21)(20,21)(40,34)
\bezier{200}(0,34)(10,23)(20,20)
\bezier{200}(20,20)(30,17)(40,6)
\put(2,26){\makebox(0,0)[cc]{\tiny $-$}} 
\put(20,7){\makebox(0,0)[cc]{\tiny $-$}} 
\put(38,26){\makebox(0,0)[cc]{\tiny $-$}}
\put(2,15){\makebox(0,0)[cc]{\tiny $+$}}
\put(37.5,15){\makebox(0,0)[cc]{\tiny $+$}}
\put(20,33){\makebox(0,0)[cc]{\tiny $+$}}
\put(1,1){$(g)$}
\end{picture} \quad
\begin{picture}(40.00,35.00)
\put(0,5){\line(1,0){40}}
\put(0,5){\line(0,1){30}}
\put(0,35){\line(1,0){40}}
\put(40,5){\line(0,1){30}}
\bezier{200}(0,20)(20,20)(0,6)
\bezier{200}(40,21)(20,21)(40,34)
\bezier{200}(0,34)(10,23)(20,20)
\bezier{200}(20,20)(30,17)(40,6)
\put(22,23){\circle{3}}
\put(2,26){\makebox(0,0)[cc]{\tiny $-$}} 
\put(20,7){\makebox(0,0)[cc]{\tiny $-$}} 
\put(38,26){\makebox(0,0)[cc]{\tiny $-$}}
\put(2,15){\makebox(0,0)[cc]{\tiny $+$}}
\put(37.5,15){\makebox(0,0)[cc]{\tiny $+$}}
\put(20,33){\makebox(0,0)[cc]{\tiny $+$}}
\put(1,1){$(h)$}
\end{picture} \quad
\begin{picture}(40.00,37.00)
\put(0,5){\line(1,0){40}}
\put(0,5){\line(0,1){30}}
\put(0,35){\line(1,0){40}}
\put(40,5){\line(0,1){30}}
\bezier{200}(0,20)(20,20)(0,6)
\bezier{200}(40,21)(20,21)(40,34)
\bezier{200}(0,34)(10,23)(20,20)
\bezier{200}(20,20)(30,17)(40,6)
\put(16,25){\circle{3}}
\put(26,21){\circle{3}}
\put(2,26){\makebox(0,0)[cc]{\tiny $-$}} 
\put(20,7){\makebox(0,0)[cc]{\tiny $-$}} 
\put(38,26){\makebox(0,0)[cc]{\tiny $-$}}
\put(2,15){\makebox(0,0)[cc]{\tiny $+$}}
\put(37.5,15){\makebox(0,0)[cc]{\tiny $+$}}
\put(20,33){\makebox(0,0)[cc]{\tiny $+$}}
\put(1,1){$(i)$}
\end{picture} \quad
\begin{picture}(40.00,37.00)
\put(0,5){\line(1,0){40}}
\put(0,5){\line(0,1){30}}
\put(0,35){\line(1,0){40}}
\put(40,5){\line(0,1){30}}
\bezier{200}(0,20)(20,20)(0,6)
\bezier{200}(40,21)(20,21)(40,34)
\bezier{200}(0,34)(10,23)(20,20)
\bezier{200}(20,20)(30,17)(40,6)
\put(16,26){\circle{3}}
\put(26,26){\circle{3}}
\put(21,23){\circle{3}}
\put(2,26){\makebox(0,0)[cc]{\tiny $-$}} 
\put(20,7){\makebox(0,0)[cc]{\tiny $-$}} 
\put(38,26){\makebox(0,0)[cc]{\tiny $-$}}
\put(2,15){\makebox(0,0)[cc]{\tiny $+$}}
\put(37.5,15){\makebox(0,0)[cc]{\tiny $+$}}
\put(20,33){\makebox(0,0)[cc]{\tiny $+$}}
\put(1,1){$(j)$}
\end{picture} \quad
\begin{picture}(40.00,37.00)
\put(0,5){\line(1,0){40}}
\put(0,5){\line(0,1){30}}
\put(0,35){\line(1,0){40}}
\put(40,5){\line(0,1){30}}
\bezier{200}(0,20)(20,20)(0,6)
\bezier{200}(40,21)(22,20)(40,34)
\bezier{200}(0,34)(10,23)(20,20)
\bezier{200}(20,20)(30,17)(40,6)
\put(21,25){\circle{2}}
\put(23,22){\circle{2}}
\put(27,20){\circle{2}}
\put(30,19){\circle{2}}
\put(2,26){\makebox(0,0)[cc]{\tiny $-$}} 
\put(20,7){\makebox(0,0)[cc]{\tiny $-$}} 
\put(36,25){\makebox(0,0)[cc]{\tiny $-$}}
\put(2,15){\makebox(0,0)[cc]{\tiny $+$}}
\put(37.5,15){\makebox(0,0)[cc]{\tiny $+$}}
\put(20,33){\makebox(0,0)[cc]{\tiny $+$}}
\put(1,1){$(k)$}
\end{picture} \quad
\begin{picture}(40.00,35.00)
\put(0,5){\line(1,0){40}}
\put(0,5){\line(0,1){30}}
\put(0,35){\line(1,0){40}}
\put(40,5){\line(0,1){30}}
\bezier{200}(0,20)(20,20)(0,6)
\bezier{200}(40,21)(20,21)(40,34)
\bezier{200}(0,34)(10,23)(20,20)
\bezier{200}(20,20)(30,17)(40,6)
\put(22,23){\circle{3}}
\put(17,17){\circle{3}}
\put(2,26){\makebox(0,0)[cc]{\tiny $-$}} 
\put(20,7){\makebox(0,0)[cc]{\tiny $-$}} 
\put(38,26){\makebox(0,0)[cc]{\tiny $-$}}
\put(2,15){\makebox(0,0)[cc]{\tiny $+$}}
\put(37.5,15){\makebox(0,0)[cc]{\tiny $+$}}
\put(20,33){\makebox(0,0)[cc]{\tiny $+$}}
\put(1,1){$(l)$}
\end{picture} \quad 
\begin{picture}(40.00,37.00)
\put(0,5){\line(1,0){40}}
\put(0,5){\line(0,1){30}}
\put(0,35){\line(1,0){40}}
\put(40,5){\line(0,1){30}}
\bezier{300}(0,6)(20,20)(40,6)
\bezier{300}(0,34)(20,20)(40,34)
\put(0,20){\line(1,0){40}}
\put(2,26){\makebox(0,0)[cc]{\tiny $-$}} 
\put(20,7){\makebox(0,0)[cc]{\tiny $-$}} 
\put(38,26){\makebox(0,0)[cc]{\tiny $-$}}
\put(2,15){\makebox(0,0)[cc]{\tiny $+$}}
\put(38,15){\makebox(0,0)[cc]{\tiny $+$}}
\put(20,33){\makebox(0,0)[cc]{\tiny $+$}}
\put(1,1){$(m)$}
\end{picture}
\caption{Perturbations of $J_{10}^3$ singularities}
\label{J103}
\end{center}
\end{figure}

\begin{proposition}[see \cite{para}]
\label{virtj103}
There are exactly 23 virtual components of the formal graph of type $J_{10}^3$. The shapes of the zero-level sets of some polynomials $f(x,y)$ representing thirteen of these components are shown in Fig.~\ref{J103}. In all these cases except those shown in pictures \ $(g)$, \ $(l)$, \ and \ $(m)$, \ the polynomials $f(x,y)$ and $-f(-x,y)$, which are related by the involution $($\ref{invol94}$)$, represent different virtual components. \hfill $\Box$
\end{proposition}

\begin{remark}\rm
Twelve topological types shown in these pictures are given in Fig.~29 of \cite{viro2}.
Pictures \ $(d)$ \ and \ $(e)$ \ of Fig.~\ref{J103} are topologically equivalent, however they represent different virtual components. 
\end{remark}

We will denote by \ $(\tilde a)$, $(\tilde b)$,\dots the pictures obtained from the corresponding pictures \ $(a)$, \ $(b)$,\dots \ of Fig.~\ref{J103} by the transformation (\ref{invol94}), that is, by the reflection in a horizontal line and simultaneous change of all signs. 

\begin{theorem}
\label{thmj3}
There are exactly 33 connected components of the complement of the discriminant variety of the versal deformation $($\ref{vers3}$)$. Each virtual component represented by one of ten pictures 
\begin{equation}
\label{listex}
(g), \ (h), \ (\tilde h), \ (i), \ (\tilde i), \ (j), \ (\tilde j), \ (k), \ (\tilde k), \ (l) 
\end{equation}
is associated with two different connected components of the discriminant complement that are obtained from each other by the involution $(\ref{invol4})$. All other virtual components are associated with single real connected components that are invariant under this involution.
\end{theorem}

\noindent
{\it Proof.} The difference of all these 33 connected components follows from the topological considerations and Bezout's theorem, see \cite{para}. By Proposition \ref{estj1}, it remains to be shown that, in all cases not covered by the list (\ref{listex}), the involution (\ref{invol4}) does not change the real connected component. 

The following polynomials (\ref{in9}), (\ref{in8}), (\ref{in7}), (\ref{in6}), are all invariant under this involution, so their virtual components are represented by single connected components of the complements of the discriminants. 

The polynomial 
\begin{equation}
\label{in9}
(x-t)^3-(x-t) y^4+3(x-t)^2+ 5(x-t) y^2-4(x-t) - C t^3
\end{equation}
with arbitrary $t>0$ realizes the picture \ $(a)$ \ for $C< 6-\frac{14}{3}\sqrt{\frac{7}{3}}$, the picture $(b)$ for 
$C \in \left(6-\frac{14}{3}\sqrt{\frac{7}{3}}, -\frac{1}{2}\right)$, and the picture \ $(d)$ \ for $C \in \left(-\frac{1}{2}, 0\right)$.

Picture \ $(m)$ \ is represented by the polynomials 
\begin{equation}
\label{in8}
x(x - y^2 -t)(x + y^2 +t), \quad t>0. 
\end{equation}

The polynomials 
\begin{equation}
\label{in7}
(x-t)^3-(x-t)y^4+3t(x-t)^2+3t(x-t)y^2+C t^3
\end{equation}
realize picture \ $(c)$ \ for any $C \in (0, 1/2)$. 

The polynomials 
\begin{equation}
\label{in6}
(x-t)^3 - (x-t) y^4 + 3t(x-t)^2 + 5t(x-t) y^2 - 4t^2(x-t) + \frac{t^2}{10}y^2 - 0.35 t^3
\end{equation}
 realize picture \ $(f)$.

\begin{figure} 
\unitlength 0.5mm
\begin{picture}(100,45)
\put(0,30){\line(1,0){100}}
\bezier{500}(12,0)(34,80)(56,0)
\bezier{500}(44,0)(66,80)(88,0)
\put(30.5,33){\tiny $-$}
\put(64.5,33){\tiny $-$}
\put(47.5,25){\tiny $-$}
\end{picture} \qquad \qquad \quad
\begin{picture}(80,45)
\bezier{300}(0,30)(20,30)(18,27)
\bezier{300}(18,27)(14,20)(10,0)
\bezier{300}(100,30)(80,30)(82,27)
\bezier{300}(82,27)(86,20)(90,0)
\bezier{200}(43,0)(50,25)(57,0)
\bezier{150}(23,36)(30,44)(37,36)
\bezier{150}(23,36)(20,32)(30,32)
\bezier{150}(30,32)(40,32)(37,36)
\bezier{150}(63,36)(70,44)(77,36)
\bezier{150}(63,36)(60,32)(70,32)
\bezier{150}(70,32)(80,32)(77,36)
\bezier{70}(52,23)(50,18)(48,23)
\bezier{70}(52,23)(55,27)(50,27)
\bezier{70}(48,23)(45,27)(50,27)
\put(27.5,34){\tiny $-$}
\put(67.5,34){\tiny $-$}
\put(47.5,22.9){\tiny $-$}
\end{picture}
\caption{Realization of perturbation $(e)$}
\label{J38}
\end{figure}

To realize the picture \ $(e)$, \ we start with the polynomial
\begin{equation}
x(x+(y-2)^2-2)(x+(y+2)^2-2) , \label{J8d}
\end{equation}
whose zero-level set is shown in Fig.~\ref{J38}. Then, we add the term $\varepsilon y^6$ with sufficiently small $\varepsilon >0$ so that the principal weighted homogeneous part will be of class $J_{10}^3$. The obtained polynomial will then keep eight real critical points (three minima and five saddlepoints), which are visible in this figure, and do not significantly change their critical values. It will also have two critical points in the complex domain. Then we add a small positive constant to the obtained polynomial in such a way that all saddlepoints will have positive critical values and all minima will stay in the negative domain. The zero-level set of the resulting polynomial will have the desired shape shown in picture $(e)$ of Fig.~\ref{J103}, see Fig.~\ref{J38} right. Then we apply a coordinate change of the form $\tilde x = x + \alpha y^2, \tilde y = \beta y$ in such a way that the principal weighted homogeneous part of the resulting polynomial will have the canonical form, see (\ref{vers3}). Finally, we apply a parallel translation of the $x$ coordinate to kill the monomial $x^2$. The resulting polynomial realizes the picture \ $(e)$. It is invariant under the involution (\ref{invol4}), in particular it has the form (\ref{vers3}) with $\lambda_2=\lambda_4=\lambda_7=\lambda_9=0$. 
 \hfill $\Box$ 

\section{The case of $P_8$ singularities}

\subsection{} The two classes $P_8^1$ and $P_8^2$ of function singularities in three variables have the canonical normal form 
\begin{equation}
\label{nfP}
x^3+y^3+z^3 -3 A x y z, \quad A \neq 1,
\end{equation}
with $A<1 $ for $P_8^1$ singularities and $A>1$ for $P_8^2$, see \cite{BM}. The zero-level sets of these normal forms of class $P_8^1$ (respectively, $P_8^2$) define cubic curves with one (respectively, two) connected component in the projective plane.
In the case of two components, one of them is a convex oval. The point $(1:1:1)$ lies in the disc bounded by this oval, because both the curve and this point are invariant under the action of the group $S(3)$ of coordinate permutations in $\R^3$.
Thus, each polynomial of class $P_8^1$ defines a conical surface in $\R^3$ that is homeomorphic to a plane and has a single singular point at the origin. The zero-level set of a polynomial of class $P_8^2$ consists of a surface like this one and a surface that is homeomorphic to the standard quadratic cone.

The standard versal deformations of $P_8^1$ and $P_8^2$ singularities consist of the polynomials
\begin{equation}
\label{versP}
f_{\lambda,A} = f_A + \lambda_1+ \lambda_2 x + \lambda_3 y + \lambda_4 z + \lambda_5 x y + \lambda_6 x z + \lambda_7 y z ,
\end{equation}
where $f_A$ is the corresponding normal form (\ref{nfP}).

\begin{definition} \rm
A polynomial $f_{\lambda,A}: ({\mathbb C}^3, {\mathbb R}^3) \to ({\mathbb C}, {\mathbb R})$ of the form (\ref{versP}) is {\em generic} if it has only Morse critical points in ${\mathbb C}^3$, its critical values are all distinct and not equal to 0, and the complex zero-level set of its principal homogeneous part $f_A$ defines a smooth curve in ${\mathbb C}P^2$.

Let $\Xi_1$ and $\Xi_2$ denote the parameter spaces $\R^7 \times (-\infty,1) $ and $\R^7 \times (1, \infty)$ of the deformations (\ref{versP}) of the functions (\ref{nfP}) of classes $P_8^1$ and $P_8^2$, respectively.
\end{definition}

\subsection{The homological invariant of real components}

\label{hominv}

Let $f$ be a polynomial in the class $\Xi_1$ or $\Xi_2$. Let $T_-$ (respectively, $T_+$) be an arbitrary negative (respectively, positive) constant which is lower (respectively, higher) than all critical values of $f$. 

Consider three subsets in $\R^3$: \ $W_{0} = f^{-1} ((-\infty, 0])$, \ 
$W_{-} = f^{-1}((-\infty, T_-]),$ \ and \ $W_{+} = f^{-1}((-\infty, T_+])$. 
If $f \in \Xi_1$ (respectively, $f \in \Xi_2$), then the space $W_-$ is contractible (respectively, consists of two components, one homotopy equivalent to a point and the other to a circle). In particular, the reduced homology groups \ $\tilde H_i(W_{-})$ \ for $f \in \Xi_2$ are equal to $\Z$ for $i=0$ or 1 and are trivial for all other $i$. 

\begin{definition} \rm
The {\em homology index} $\mbox{HI}( f)$ of a non-discriminant polynomial $f$ of class $\Xi_1$ or $\Xi_2$ is the triple of Betti numbers $(b_0, b_1, b_2)$ of the relative homology group \ $H_*(W_{0}, W_{-})$. 
\end{definition}

It is convenient to consider this homology group in the context of the exact homological sequence of the triple \begin{equation} 
\label{triple}
(W_{+}, W_{0}, W_{-}). 
\end{equation} 
The groups \ $H_*(W_{+}, W_{-})$ \ participating in this exact sequence are the same for all $f$ from the same space $\Xi_i$, because the pairs of sets \ $(W_{+}, W_{-})$ \ corresponding to different functions $f $ form a trivializable fiber bundle over these spaces. In the case $\Xi_1$ these groups are trivial. For $\Xi_2$, they are equal to $\Z$ for $i=1$ or $2$ and are trivial for all other $i$. 

Clearly, the homology index is an invariant of the connected components of discriminant complement, and the invariant \ $\mbox{\rm Ind}$ \ is equal to the alternating sum of three elements of the homology index. 

The values that the invariant $\mbox{\rm Ind}$ can take on polynomials of classes $\Xi_1$ and $\Xi_2$ are described in Propositions \ref{prop1} and \ref{proP2} below. These values are $1, 0, -1, -2,$ and $-3$. The closure in $\RP^3$ of the set of zeros of any non-discriminant perturbation $f$ of a singularity (\ref{nfP}) with \ $\mbox{\rm Ind}(f)$ \ equal to 0 $($respectively, $-1$, $-2$, $-3$, or 1$)$ is homeomorphic to the projective plane $($respectively, the projective plane with one, two, or three handles, or the disjoint union of the projective plane and the sphere$)$.

\subsection{The symmetry group and the relation between the real and virtual components}

The virtual functions, formal graphs, and virtual components of the $P_8^*$ singularities are defined almost identically to those of the $X_9$ singularities, see \S~\ref{virtx}, with only the following alterations: the frames of vanishing cycles in the groups 
$H_2(f^{-1}(0))$ for non-discriminant polynomials $f \in \Xi_1$ or $\Xi_2$ consist of eight elements, and the data sets of virtual functions only include the {\em parities} of the Morse indices of critical points, and not their integer values. We also have the following analog of Proposition \ref{propmain}.

\begin{proposition}[see \cite{vasP8}, Theorem 10]
\label{propmainZ}
 For any generic Morse polynomial $f$ of type $\Xi_1$ or $\Xi_2$, a virtual Morse function $\varphi$ associated with it, and a virtual Morse function $\tilde \varphi$ connected with $\varphi$ by an edge of the formal graph, there exists a generic Morse polynomial $\tilde f$ and a path in the space $\Xi_1$ or \ $\Xi_2$ \ connecting $f$ and $\tilde f$ and containing only one non-generic point at which it experiences a
standard surgery of the same type as the edge $[\varphi, \tilde \varphi]$.
\end{proposition}

\begin{corollary}
\label{propmainP}
\label{cormain}
1. Any virtual function of the class \ $\Xi_1$ or $\Xi_2$ is associated with a real polynomial of the same type. 

2. If two virtual functions of the class \ $\Xi_1$ \ $($respectively, \ $\Xi_2)$ belong to the same virtual component, then they are associated with two generic polynomials of type $\Xi_1$ $($respectively, $\Xi_2)$ that belong to the same connected component of the space of generic polynomials. \hfill $\Box$
\end{corollary}

The group $S(3)$ of coordinate permutations in $\R^3$ acts on the spaces \ $\Xi_1$ \ and \ $\Xi_2$ preserving the normal form (\ref{nfP}) and the discriminant variety. This action is generated by the following transformations:
\begin{equation}
\label{symP}
f (x,y,z) \mapsto f(y,x, z) \quad \mbox{and} \quad f(x,y,z) \mapsto f(y, z, x) .
\end{equation}
The first transformation is an involution, and the second is of order three.

Consequently, this group acts on the set of the connected components of the discriminant complements in the spaces $\Xi_*$. 
The transformations (\ref{symP}) preserve all topological characteristics of the functions. Therefore, all these connected components within the same orbit of this action are associated with the same virtual component of the formal graph of type $P_8^1$ or $P_8^2$.

Additionally, the involution 
\begin{equation}
\label{invoP}
f(x,y,z) \mapsto -f(-x,-y,-z)
\end{equation}
also acts on the spaces $\Xi_1$ and $ \Xi_2$.
 It defines an involution on both sets of real and virtual components.
This involution preserves the invariant \ $\mbox{Ind}$. 

If $f$ and $\tilde f$ are two polynomials related by the involution (\ref{invoP}), then the corresponding exact sequences (\ref{triple}) are related by a kind of Poincar\'e duality: the group \ $H_*(W_{+}, W_{0})$ \ for the polynomial $f$ is naturally dual to \ $H_*(W_{0}, W_{-})$ \ for $\tilde f$, and vice versa.

\begin{proposition}
\label{estp1}
If two generic polynomials, $f$ and $\tilde f$, of class $\Xi_1$ or $\Xi_2$ are associated with the same virtual function, then the connected components of the space of generic polynomials of this class containing $f$ and $\tilde f$ are mapped to each other by an element of the group $ S(3)$. In particular, in this case the connected components of the discriminant complement containing these polynomials are also mapped to each other by an element of the group $ S(3)$.
\end{proposition}

\begin{definition} \rm
For any polynomial $f: \R^3 \to \R$ of degree three with principal homogeneous part of the form (\ref{nfP}), ${\bf T}(f)$ is the unique polynomial of the form (\ref{versP}) with the same principal homogeneous part, which is obtained from $f$ by a parallel translation in $\R^3$. Specifically, it is the translation by the vector $\frac{1}{3}(a,b,c)$, where $a, b,$ and $c$ are the coefficients of $f$ at the monomials $x^2$, $y^2$, and $z^2$.
\end{definition}

\begin{lemma}
\label{le83}
The set of all polynomials of the form $(\ref{versP})$ having a critical point of class $E_6$ with zero critical value and a critical point of class $A_2$ with critical value $\frac{64}{3}$ and signature $($1,1$)$ of the second differential at the latter point consists of the polynomial 
\begin{equation} \label{ean}
{\bf T}\left(x^3 + y^3 + z^3 + 6 x y z + (x +y -2z)^2\right)
\end{equation}
and two other polynomials obtained from this one by permuting the coordinates.
\end{lemma}

\noindent
{\it Proof.} According to Lemma 38 of \cite{vasP8}, there are only six polynomials of the form (\ref{versP}) that have a critical point of class $E_6$ with zero critical value and a critical point of class $A_2$ with critical value $\frac{64}{3}$. Three of these polynomials (mentioned in formula (27) of \cite{vasP8} with a certain value of the parameter $\beta$) have signature (0,2) of the second differential at the $A_2$ point. The remaining three polynomials have the form (\ref{ean}) up to the permutations of coordinates, see formula (28) in \cite{vasP8}. In particular, they belong to one orbit of the group $S(3)$.
\hfill $\Box$

\begin{lemma}
\label{le38}
The set of all polynomials of of the form $($\ref{versP}$)$ and class $\Xi_2$ that have a critical point of class $E_6$ with zero critical value, and also two real critical points of class $A_1$ with the same critical value $\frac{64}{3}+\frac{128}{3\sqrt{3}}$, consists of three polynomials: the polynomial
\begin{equation}
\label{cur2}
{\bf T}\left(x^3 + y^3 + z^3 - 3(1+\sqrt{3})x y z + (x + y+(1+\sqrt{3})z)^2 \right), \ 
\end{equation}
and two polynomials obtained from $(\ref{cur2})$ by permuting the coordinates.
\end{lemma}

\noindent
{\it Proof.} This lemma follows almost immediately from statement 2 of Lemma 47 in \cite{vasP8} (the coefficients $\varkappa$ and $c$ in formula (45) of \cite{vasP8} are fixed by our conditions on the critical values). \hfill $\Box$ \medskip

\noindent
{\it Proof of Proposition \ref{estp1}.} The deduction of this proposition for the case $\Xi_1$ (respectively, $\Xi_2$) from Lemma \ref{le83} (respectively, \ref{le38}) 
 repeats the deduction of Proposition \ref{estx+} from Lemma \ref{lemx+}. Specifically, the neighborhood of the polynomial (\ref{ean}) or (\ref{cur2}) in $\Xi_2$ has only one non-trivial automorphism that commutes with the Lyashko--Looijenga map: the involution $f(x,y,z) \leftrightarrow f(y,x,z)$. \hfill $\Box$

\subsection{Components of discriminant complements of $P_8^1$ singularities}

\begin{proposition}[see \cite{para}, Proposition 7]
\label{prop1}
1. There are exactly 6503 virtual functions of the class $P_8^1$.

2. There are exactly seven virtual components of the formal graph of type $P_8^1$: two with $\mbox{\rm Ind}=1$, two with $\mbox{\rm Ind}=-3$, and one virtual component with each value $0$, $-1$, and $-2$ of the {\rm Ind} invariant. 

3. The involution $(\ref{invoP})$ acts non-trivially on each pair of virtual components with $\mbox{\rm Ind}=1$ or $\mbox{\rm Ind} = -3$. \hfill $\Box$
\end{proposition}

\begin{theorem}
\label{mthmp81}
Each virtual component of type $P_8^1$ is associated with one and only one connected component of the space of non-discriminant polynomials in class $\Xi_1$.
In particular, there are exactly seven such connected components.
\end{theorem}

\noindent
{\it Proof.}
The invariant $ \mbox{Ind}$ of the polynomial
\begin{equation}
x^3 + y^3 +z^3 - 3 (x + y +z) + C 
\end{equation}
is equal to $0 $ \ for \ $ C< -6$ and $C > 6$, \ to \ $1$ for \ $C \in (-6,-2)$ \ and \ $C \in (2, 6)$, \ and to \ $-2$ \ for \ $C \in (-2,2)$. 
These polynomials are invariant under the action of the group $S(3)$ of permutations of coordinates. Therefore, according to Proposition \ref{estp1}, all virtual components represented by these polynomials
are associated with single real connected components of the set of non-discriminant polynomials. 

The two real components with $\mbox{Ind}=1$ corresponding to $C \in (-6,-2)$ and $C \in (2,6)$ are separated by the homology invariant of \S~\ref{hominv}. For the first component, this invariant is equal to $(1,0,0)$, and for the second it is equal to $(0,0,1)$. According to the last sentence of the previous paragraph, the associated virtual components are also different.

The polynomials 
\begin{equation}
{\bf T}\left(x^3+y^3+z^3+3 (x + y + z)^2 - \frac{1}{2}(x^2+y^2+z^2)+\frac{7}{36}\right)
\end{equation}
and 
\begin{equation}
{\bf T}\left(x^3+y^3+z^3-3 (x + y + z)^2 + \frac{1}{2} (x^2+y^2+z^2)-\frac{7}{36}\right)
\end{equation} 
have $\mbox{Ind}=-3$. These polynomials are invariant under the entire group $S(3)$. Therefore, the virtual components associated with them are associated only with the real components of discriminant complements that contain these polynomials.
 These two polynomials are related to each other by the involution (\ref{invoP}). Therefore, according to the final statement of Proposition \ref{prop1}, their associated virtual components are different. 
Thus, Theorem \ref{mthmp81} is proved for all virtual components with \ $\mbox{Ind}$ \ equal to $1, 0, -2,$ and $-3$. Let us prove it for the only remaining virtual component with $\mbox{Ind}=-1$. 

Let $L \sim S^1$ be the space of all linear functions $\R^3 \to \R$ of unit Euclidean norm that vanish on the vector $(1,1,1)$. 
For any $ l \in L$, the restriction of the polynomial \ $x^3+y^3+z^3$ \ to the plane \ $\ker l$ \ is a non-degenerate cubic form with only one line of zeros. Therefore, in some linear coordinates $(\tilde p, \tilde q)$ in this plane this restriction is equal to $\tilde p^3+\tilde p \tilde q^2$. Let $p_l$ and $q_l$ be the linear functions $\R^3 \to \R$ that are constant along the lines perpendicular to the plane \ $\ker l$ \ and whose restrictions to this plane coincide with $\tilde p$ and $\tilde q$, respectively. 
The coordinate $p_l$ is uniquely defined by the plane \ $\ker l$, while $q_l$ is defined up to a sign, so that the monomial $q^2_l$ is also uniquely defined by this plane.

\begin{lemma}
\label{le338}
For each linear function $l \in L$, 

1$)$ the polynomial 
\begin{equation}
\label{Pl} x^3+y^3+z^3 - l^2
\end{equation}
has only one real critical point, namely 
a critical point of class $D_4^+$ at the origin. In particular, this polynomial has four critical points $($counted with multiplicities$)$ in the complex domain.
 
2$)$ For any sufficiently small $\varepsilon>0$, the polynomial \begin{equation}
x^3+y^3+z^3 - l^2 - \varepsilon p_l^2
\label{Ql}
\end{equation} 
has only two critical points in $\R^3$:
a critical point of type $A_3$ at the origin, at which this polynomial is equal to 
\begin{equation}
\label{nfzet}
\xi^4 -\eta^2-\zeta^2
\end{equation}
 in some local coordinates, and a Morse critical point with Morse index $1$ and a negative critical value. 
\end{lemma}

\noindent
{\it Proof.} 1. The sum of the three partial derivatives of the polynomial (\ref{Pl})
is identically equal to $3(x^2+y^2+z^2),$ so the origin is the only critical point of (\ref{Pl}) in the real domain.
In coordinates $l, p_l, q_l$, this polynomial has the form 
\begin{equation}
\label{nco}
-l^2 + p_l^3 + p_l q_l^2 + l \psi_2(p_l, q_l) + l^2 \psi_1(p_l, q_l, l),
\end{equation}
where $\psi_2$ and $\psi_1$ are some homogeneous polynomials of degrees 2 and 1, respectively.

The lowest weighted homogeneous part of the polynomial (\ref{nco}) in the coordinates $l, p_l,$ and $ q_l$ 
with weights $\deg l=3$ and $ \deg p_l = \deg q_l =2$
is equal to $-l^2+p_l^3+p_l q_l^2.$ Therefore, by the techniques of \cite{AVGZ}, \S 12, this polynomial can be reduced to the normal form $-\tilde z^2+\tilde x^3+\tilde x \tilde y^2$ of $D_4^+$ class in some local coordinates $\tilde x, \tilde y,$ and $\tilde z$.
 
2. Let $\alpha$ be the coefficient at the monomial $q_l^2$ of the polynomial $\psi_2$ in (\ref{nco}). Then the lowest weighted homogeneous part of the polynomial (\ref{Ql}) with weights $\deg p_l=2, \deg q_l=1, $ and $\deg l=2$ is equal to $-l^2 -\varepsilon p_l^2 + p_l q_l^2 + \alpha l q_l^2$.
Substituting $\tilde l = l-\frac{\alpha q_l^2}{2}$ and $ P_l =p_l - \frac{q_l^2}{2\varepsilon}$
turns this polynomial into the polynomial with the lowest weighted homogeneous part 
$$-\tilde l^2 -\varepsilon P_l^2 + \left(\frac{\alpha^2}{4} + \frac{1}{4\varepsilon}\right) q_l^4 \ .$$ Again, using the techniques of \S 12 of \cite{AVGZ}, this function can be reduced to the standard form (\ref{nfzet}) of $A_3$ singularities via a local change of coordinates. If $\varepsilon$ is sufficiently small, then adding $-\varepsilon p_l^2$ to 
(\ref{Pl}) does not return the non-real critical points of (\ref{Pl}) to the real domain. By the index considerations, the obtained polynomial $F(l)$ has one real Morse critical point in addition to the $A_3$ point at the origin. This Morse critical point depends continuously on $l$, therefore its Morse index is the same for all $l \in L$. It is easy to calculate that for $l \equiv x-y$ it is equal to 1. 

The restriction of the polynomial (\ref{Ql}) to the line connecting the origin and this Morse point is a degree three polynomial with vanishing derivative at these points; its critical point at the origin is a local maximum. Therefore, its critical value at the Morse critical point is negative. \hfill $\Box$
\medskip

Denote by $C(l)$ the critical value of the polynomial (\ref{Ql}) at its Morse critical point. Let $F(l)$ be the polynomial obtained from (\ref{Ql}) by subtracting the constant $C(l)/2$. Then the polynomial $F(l)$ is non-discriminant and has $\mbox{Ind}=-1$. 

The latter polynomial is canonically determined by the linear form $l$
and depends continuously on it. The union of polynomials $F(l) $ over all $l \in L$ is connected and invariant under the action of the group $S(3)$. Therefore, the virtual component with $\mbox{Ind}=-1$ is also associated with only one component of the discriminant complement. This completes the proof of Theorem \ref{mthmp81}.
\hfill $\Box$

\subsection{One-dimensional cohomology of the discriminant com\-ple\-ment of $P_8^1$ singularities}

\begin{theorem}
\label{t2hom}
The first cohomology group of the connected component of the discriminant com\-ple\-ment in $\Xi_1$ that consists of polynomials with $\mbox{Ind} = -1$ is non-trivial.
\end{theorem}

\noindent
{\it Proof.} The zero-level sets of the polynomials in this component are homeomorphic to the plane with a handle. Indeed, the topology of these sets remains constant along the component. For any polynomial $F(l)$ from the previous proof this set is obtained from the very low level sets (homeomorphic to a plane) by passing through a Morse critical point of Morse index 1. 

Thus, the one-dimensional homology groups of all zero-level sets of class $\Xi_1$ polynomials with $\mbox{Ind}=-1$ are isomorphic to $\Z^2$. These homology groups form a local system over the space of all such polynomials, i.e., over the connected component of the discriminant complement under consideration. In particular, the fundamental group of this component acts on these groups by monodromy. For each $l \in L$, one basic homology class of the zero-level set of the polynomial $F(l)$ is defined by the intersection of this set and the plane $\ker l$. When $l \in L$ runs a half of the circle $L$, the planes $\ker l$ \ run a closed cycle in $\RP^2$, and the polynomials $ F(l)$ run a closed cycle in the discriminant complement. The monodromy over the latter cycle brings our basic 1-cycle to itself with the opposite orientation. Thus, the half-circle $L/\Z_2$ defines non-trivial elements of the fundamental group and of the mod 2 homology group of the component. 
 \hfill $\Box$

\subsection{Components of the discriminant com\-ple\-ment of the $P_8^2$ singularities}

\begin{definition} \rm
The invariant $\mbox{Card}$ of a virtual component of the formal graph is the number of vertices in that virtual component. 
\end{definition}

The involution (\ref{invoP}) of virtual components of classes $\Xi_1$ and $\Xi_2$ preserves this invariant.

\begin{table}
\caption{Statistics for $P_8^2$ singularity}
\label{p82stat}
\noindent
{\footnotesize
\begin{tabular}{|c|cccc|ccc|ccc|cc|ccc|c|}
\hline
No. & 1 & 2 & 3 & 4 & 5 & 6 & 7 & 8 & 9 & 10 & 11 & 12 & 13 & 14 & 15 \\
\hline
Ind & \multicolumn{4}{c|}{$-3$} & \multicolumn{3}{c|}{$-2$} & \multicolumn{3}{c|}{$-1$} & \multicolumn{2}{c|}{0} & \multicolumn{3}{c|}{1} \\
\hline
Card & 258 & 156 & 60 & 60 & 1216 & 336 & 336 & 1318 & 844 & 844 & 1648 & 1648 & 262 & 94 & 94 \\
\hline
\end{tabular}
}
\end{table}

\begin{proposition}[see \cite{para}, Proposition 8]
\label{proP2}
1. There are exactly 9174 distinct virtual functions of class $\Xi_2$.

2. These virtual functions are distributed among fifteen virtual components. The values of invariants \ $\mbox{\rm Ind}$ \ and \ $ \mbox{\rm Card}$ \ of these components are shown in the second and the third rows of Table \ref{p82stat}. 

3. All pairs of different virtual components of class $\Xi_2$ with the same \ $\mbox{\rm Card}$ \ value $($i.e., the pairs of components 3 and 4, 6 and 7, 9 and 10, 11 and 12, and 14 and 15 in Table \ref{p82stat}$)$ can be realized by generic polynomials connected by involution $($\ref{invoP}$)$.
\end{proposition}

\begin{theorem}
\label{mthmp2}
Each virtual component of the formal graph of type $P_8^2$ is associated with only one real connected component of the discriminant complement. In particular, there are exactly fifteen such connected components.
\end{theorem}

Proof of this theorem takes the rest of this section.
\smallskip

All of the polynomials given in the formulas (\ref{ci}), (\ref{ci3}), (\ref{ci4}), and 
(\ref{e258}) below are invariant under the $S(3)$-action. Thus, if $f$ is any of these polynomials and $C$ is a noncritical value of $f$, then Proposition \ref{estp1} implies that the virtual component associated with the polynomial $f-C$ is associated with only one real component of the discriminant complement.
 We will show that all virtual components with $\mbox{Card}$ equal to 258, 336, 1648, 262, or 94 can be realized in this way. Hence, the assertion of Theorem \ref{mthmp2} is true for these seven components.
\medskip

\noindent
{\bf Components with $\mbox{\rm Card}=1648$.} 
For each polynomial $f$ of class $\Xi_2$, the polynomials $f-C$ with sufficiently large $|C|$ and different signs of $C$ realize two different virtual components with $\mbox{Ind}=0$. The homology indices of these components are equal to $(0,0,0)$ and $(0,1,1)$. 
For example, let $f$ be the $S(3)$-invariant polynomial 
\begin{equation}
\label{ci}
 x^3+y^3+z^3 -4x y z .
\end{equation}
The polynomial $f-C$ with arbitrary $C \neq 0$ then represents one or the other of these two components depending on the sign of $C$. These polynomials are not Morse, but all their sufficiently small generic perturbations have $\mbox{Card}=1648$.
\medskip

\noindent
{\bf Component with $\mbox{\rm Card}=262$.} The polynomial 
\begin{equation}
\label{ci3}
x^3+y^3+z^3-4xyz +(x+y+z)
\end{equation}
has only two real Morse critical points: one with Morse index 2 and critical value $-2$, and the other with Morse index 1 and critical value 2. Thus, it has
$\mbox{Ind}=1$ and homology index $(0,0,1)$. This polynomial is invariant under the involution (\ref{invoP}). According to the last statement of Proposition \ref{proP2}, the associated virtual component cannot have $\mbox{Card}=94$. Therefore, this is the only remaining component with $\mbox{Ind}=1$, which has \ $\mbox{Card}=262$. 
\medskip

\noindent
{\bf Components with $\mbox{\rm Card}=94$ and $336$.} 
The polynomial 
\begin{equation}
\label{ci4}
f \equiv x^3+y^3+z^3-6xyz +3(x^2+y^2+z^2)
\end{equation}
has eight real Morse critical points: one local minimum with critical value 0, 
three critical points with Morse index $1$ and critical value $\frac{88}{49}$,
three critical points with Morse index $2$ and critical value $4$, and one critical point with Morse index $1$ and critical value $12$. Therefore, the invariant \ $\mbox{Ind}(f-C)$ \ is equal to \ $0$ \ for \ $C<0$, to \ $1$ \ for \ $C\in \left(0, \frac{88}{49}\right),$ \ to \ $-2$ \ for \ $C \in \left(\frac{88}{49}, 4\right),$ \ to \ 1 \ for \ $C\in \left(4, 12\right),$ \ and to \ 0 \ for \ $C> 12$. 
 
The group \ $H_*(W_{0}, W_{-})$ corresponding to any polynomial \ $f-C$ \ with $C \in \left(0, \frac{88}{49}\right) $ is isomorphic to $\Z$ in dimension 0 and is trivial in all other dimensions. Therefore, the homology index $\mbox{HI}(f- 1)$
is $ (1,0,0)$. The involution (\ref{invoP}) transforms this polynomial into one with $\mbox{HI}=(0,1,2)$. In particular, these two polynomials represent two different components. They are $S(3)$-invariant. The associated virtual components are therefore associated with single real connected components and are also different, though their $\mbox{Card}$ invariants are the same. According to Proposition \ref{proP2} and Table \ref{p82stat}, these are two components with $\mbox{Card}=94$.

According to elementary Morse theory, the homology index of each polynomial $f-C$ with $C \in \left(\frac{88}{49}, 4 \right)$ is equal to $(0,2,0)$ or $(1,3,0)$, while its image under the involution (\ref{invoP}) has the homology index $(0,3,1)$. Therefore, these two polynomials belong to different components with $\mbox{Ind}=-2$ and equal $\mbox{Card} $ invariants. The only possibility for this is two components with $\mbox{Card}=336.$
\medskip

\noindent
{\bf Component with $\mbox{\rm Card}=258$.} The polynomial 
\begin{equation}
\label{e258}
f= x^3+y^3+z^3-12xyz + 3(x^2+y^2+z^2)
\end{equation}
has eight real Morse critical points: one local minimum with critical value 0; 
three critical points with Morse index 1 and critical value \ $\frac{88}{147}$; \
one critical point with Morse index 1 and critical value \ $\frac{4}{3}$; \ and
three critical points with Morse index 2 and critical value \ $4$. \ Thus, \ $\mbox{Ind}(f-C)$ \ is equal to \ $0$ \ for \ $C<0$, \ to \ $1$ \ for \ $C \in \left(0, \frac{88}{147}\right)$, \ to \ $-2$ \ for $C \in \left(\frac{88}{147}, \frac{4}{3} \right)$, \ to \ $-3$ \ for \ $C \in \left(\frac{4}{3}, 4\right)$, \ and to \ 
0 \ for \ $C> 4$. 

\begin{proposition} The polynomials $f-C$, where $f$ is given by $(\ref{e258})$
and \ $C \in \left(\frac{4}{3}, 4\right),$ \ have \ $\mbox{Card}=258.$
\end{proposition}

\noindent
{\it Proof.} By Proposition 40 of \cite{vasP8}, the intersection matrix of the canonical basis of vanishing cycles for any generic polynomial of class $\Xi_2$ with eight real critical points, one of which is an extremum, is expressed by the Coxeter-Dynkin graph shown in 
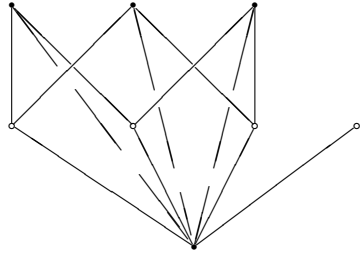
\begin{figure}
\unitlength=0.6mm
\begin{center}
\begin{picture}(118,88)
\put(60,5){\circle*{1.5}} 
\put(0,45){\circle{1.5}} 
\put(40,45){\circle{1.5}} 
\put(80,45){\circle{1.5}} 
\put(0,85){\circle*{1.5}} 
\put(40,85){\circle*{1.5}} 
\put(80,85){\circle*{1.5}} 
\put(113.5,45){\circle{1.5}} 
\put(59.3,5.7){\line(-3,2){58.1}}
\put(59.5,5.8){\line(-1,2){19}}
\put(60.5,5.8){\line(1,2){19}}
\put(0.7,45.7){\line(1,1){38.5}}
\put(40.7,45.7){\line(1,1){38.5}}
\put(60.8,5.6){\line(4,3){51.7}}
\put(0,46){\line(0,1){38.3}}
\put(19.5,65.5){\line(-1,1){18.5}}
\put(20.5,64.5){\line(1,-1){19}}
\put(80,46){\line(0,1){38.5}}
\put(59.5,65.5){\line(-1,1){18.6}}
\put(60.5,64.5){\line(1,-1){19}}
\put(15,65){\line(-3,4){14.5}}
\put(21,57){\line(3,-4){15}}
\put(42,29){\line(3,-4){16}}
\put(47,57){\line(-1,4){6.7}}
\put(48.5,51){\line(1,-4){5}}
\put(55.5,23){\line(1,-4){3.5}}
\put(73,57){\line(1,4){6.7}}
\put(71.5,51){\line(-1,-4){5}}
\put(65,25){\line(-1,-4){4.8}}
\end{picture}
\end{center}
\caption{Coxeter--Dynkin graph for the polynomial (\ref{e258})}
\label{504}
\end{figure}
Fig.~\ref{504}. This graph does not contain edges of multiplicities other than 1 and $-1$. Upon a special request, our program verified that, for each virtual function of $P_8^2$ type from the virtual components with \ $\mbox{Card}=60$ or $156$, \ the intersection index of some two different vanishing cycles of the canonical basis is equal to $-2$. Therefore, the virtual component associated with our polynomial is the only remaining virtual component with $\mbox{Ind}=-3$.
\hfill $\Box$

\begin{remark} \rm 
It is easy to check that the polynomials $f-C$ with $f$ given by (\ref{e258}) and \
$C \in \left(\frac{88}{147}, \frac{4}{3} \right)$ \ or \ $C \in \left(0, \frac{88}{147}\right)$ \ give new realizations of components with \ $\mbox{Card} = 336$ \ and $94$, \ respectively.
\end{remark}

\noindent
{\bf Components with \ $\mbox{Card}=60$.}
\label{vc60}

\begin{proposition}
\label{th43}
There exists a $S(3)$-invariant polynomial $f \in \Xi_2$ that realizes a virtual component with $\mbox{\rm Card}=60.$
\end{proposition}

\noindent
{\it Proof.} According to Proposition \ref{proP2}, two virtual components exist with $\mbox{Card}=60$. Our program shows that one of these components contains a virtual function that models a generic polynomial of class $\Xi_2$ such that

1) all of its critical points are real, 

2) the three highest critical values are positive and the other five critical values are negative, 

\begin{figure}
\unitlength=0.6mm
\begin{center}
\begin{picture}(123,52)
\put(0,5){\circle{1.5}} 
\put(40,5){\circle{1.5}}
\put(120,5){\circle{1.5}}
\put(120,45){\circle*{1.5}}
\put(0,45){\circle*{1.5}}
\put(40,45){\circle*{1.5}}
\put(80,45){\circle*{1.5}}
\put(80,5){\circle{1.5}} 
\put(0,6){\line(0,1){38.5}}
\put(0.8,5.6){\line(2,1){78.5}}
\put(40,26){\line(0,1){18.1}}
\put(40,24.6){\line(0,-1){19}}
\put(40.7,5.7){\line(1,1){38.5}}
\put(119.3,5.7){\line(-1,1){38.5}}
\put(120,6){\line(0,1){38.5}}
\put(41,5){\line(1,0){38.5}}
\put(119,5){\line(-1,0){38.5}}
\bezier{600}(0.7,4.3)(50,-5)(79,4)
\put(80.5,14){\line(0,-1){8}}
\put(79.5,14){\line(0,-1){8}}
\put(80.5,16){\line(0,1){8}}
\put(79.5,16){\line(0,1){8}}
\put(80.5,26){\line(0,1){8}}
\put(79.5,26){\line(0,1){8}}
\put(80.5,36){\line(0,1){8}}
\put(79.5,36){\line(0,1){8}}
\put(0,0){\makebox(0,0)[cc]{\tiny 1}}
\put(35,5){\makebox(0,0)[cc]{\tiny 2}}
\put(80,0){\makebox(0,0)[cc]{\tiny 5}}
\put(120,0){\makebox(0,0)[cc]{\tiny 3}}
\put(0,51){\makebox(0,0)[cc]{\tiny 6}}
\put(40,51){\makebox(0,0)[cc]{\tiny 7}}
\put(80,51){\makebox(0,0)[cc]{\tiny 4}}
\put(120,51){\makebox(0,0)[cc]{\tiny 8}}
\end{picture}
\end{center}
\caption{Coxeter-Dynkin graph of the matrix (\ref{matp8-2})}
\label{945}
\end{figure}

3) the intersection matrix of the canonically oriented vanishing cycles, ordered by increasing critical values, is as follows:
\begin{equation}
 \left|
\begin{array}{cccccccc}
$-2$ \ & 0 & 0 & 1 & 1 & 1 & 0 & 0 \\
0 & $-2$ \ & 0 & 1 & 1 & 0 & 1 & 0 \\
0 & 0 & $-2$ \ & 1 & 1 & 0 & 0 & 1 \\
1 & 1 & 1 & $-2$ \ & $-2$ & 0 & 0 & 0 \\
1 & 1 & 1 & $-2$ & $-2$ \ & 0 & 0 & 0 \\
1 & 0 & 0 & 0 & 0 & $-2$ \ & 0 & 0\\
0 & 1 & 0 & 0 & 0 & 0 & $-2$ \ & 0 \\
0 & 0 & 1 & 0 & 0 & 0 & 0 & $-2$ 
\end{array}
\right|\label{matp8-2}
\end{equation}
and hence its Coxeter-Dynkin graph is shown in Fig.~\ref{945}, 

4) the Morse indices of the points number 1, 2, 3, and 5 with respect to this order are odd, while the Morse indices of points 4, 6, 7, and 8 are even. 

According to Corollary \ref{propmainP}, there exists a real polynomial $f_0 \in \Xi_2$ with these data. 
According to Proposition 41 of \cite{vasP8}, for all functions in class $\Xi_2$ that satisfy these conditions 1)--4), the odd Morse indices of their critical points are all equal to 1, and the even Morse indices are all equal to 2. Therefore, the polynomial $f_0$ has no local minima or maxima.

Let $ s: [0,1]\to \R^8$ be the parameterized line segment such that $s(0)$ is the
 collection of critical values $(c_1, \dots, c_8)$ of the polynomial $f_0$, ordered by increasing value, and $s(1) = (c_-, c_-, c_-, c_4, c_5, c_+, c_+, c_+)$, where $c_-$ is the mean value of $c_1, c_2,$ and $c_3$, and $c_+$ is the mean value of $c_6, c_7,$ and $c_8$. Using the Lyashko--Looijenga covering, whose existence in the spaces $\Xi_2$ is shown in \cite{Jaw2}, one can uniquely construct a path $I: [0,1) \to \Xi_2$ such that $I(0)=f_0$ and, for any $t\in [0,1)$, $s(t) $ is the sequence of critical values of the polynomial $I(t)$. As in the proof of Theorem 10 of \cite{vasP8}, the limit point $I(1) \in \Xi_2$ of this lifted path is well-defined. Indeed, according to Proposition 2 of \cite{Jaw2}, this is guaranteed by the fact that all intersection forms corresponding to the limit positions $c_-$ and $c_+$ of the critical values of the polynomials $I(t)$ are elliptic. (These forms are given respectively by the upper left and the bottom right $3 \times 3$ corners of the matrix (\ref{matp8-2}).)
 The point $I(1)$ is a Morse (albeit non-generic) polynomial of class $\Xi_2$ having three critical points with critical value $c_-$ and three critical points with critical value $c_+$. 

\begin{lemma}
This polynomial $I(1)$ is $S(3)$-invariant.
\end{lemma}

\noindent
{\it Proof.} The critical points of the polynomial $I(1)$ are ordered according to the previous ordering of the critical points of the polynomial $f_0$, from which they were obtained via continuous deformation along the path $I$. This ordering can be uniquely extended to a continuous system of orderings of the critical points of each polynomial $\tilde f$ from a sufficiently small neighborhood $U$ of the point $I(1)$ in $\Xi_2$.
Define the map $\Lambda: U \to \R^8 $, mapping each point $\tilde f \in U$ to the sequence of its correspondingly ordered critical values. 
By the versality property of the family $\Xi_2$, this map is a local diffeomorphism. 
Let $\varepsilon >0$ be a very small number. 
For any permutation $ \sigma \equiv \binom{1\ 2\ 3}{i_1\, i_2\, i_3} $
of critical points of $I(1)$ with critical value $c_-$, there is a unique point $\tilde f_\sigma \in U$, whose critical values at the critical points with numbers \ $i_1, 4, 5,$ and $5+i_1$ \ are the same as those of the polynomial $I(1)$, the values at the points $i_2$ and $5+i_2$ are equal respectively to $c_- +\varepsilon$ and $c_+ + \varepsilon$, and the values at the points $i_3$ and $5+i_3$ are equal respectively to $c_- -\varepsilon$ and $c_+ - \varepsilon$. Different permutations $\sigma$ define six different points $\tilde f_\sigma \in U$. These points are all associated with the same virtual function, because any permutation $\sigma$ of the first three rows and columns in (\ref{matp8-2}) and the simultaneous permutation of the last three rows and columns by the same rule $\sigma$ preserves the matrix. 

As in the proof of Proposition \ref{estp1},
the six points $\tilde f_\sigma \in U$ are all mapped to each other by the action of elements of the group $S(3)$. The order of this group is also six, therefore each element of the group only permutes these six points. In particular, it maps a point from an
arbitrarily small neighborhood $U$ of the polynomial $I(1)$ to another point in the same neighborhood. By continuity, it also maps the polynomial $I(1)$ to itself.
\hfill $\Box$ \medskip

Due to the construction of the polynomial $I(1)$, the $\mbox{Card}$ invariant of arbitrarily small generic perturbations of this polynomial is equal to that of the initial polynomial $f_0$, that is, it is 
equal to 60. This proves Proposition \ref{th43}. \hfill $\Box$ \medskip

By statement 3 of Proposition \ref{proP2}, the involution (\ref{invoP}) moves the polynomial $I(1)$ to a polynomial in a different component of the discriminant complement with $\mbox{Card}=60$. Thus, this component is also represented by an $S(3)$-invariant polynomial. This proves the assertion of Theorem \ref{mthmp2} for both virtual components with $\mbox{Card}=60$.
\medskip

\noindent
{\bf Component with $\mbox{Card}=1216$.}
\label{vc1216}
Denote by $\hat f$ the polynomial \begin{equation}
\label{1216}
x^3+y^3+z^3-6xyz+3(x-y)^2 .
\end{equation}
This polynomial has 
a critical point of class $D_4^-$ with critical value 0 at the origin,
two Morse critical points with Morse index 2 and a common critical value $\omega_1 \approx 3.14$, and two non-real critical points with a common critical value $\omega_2 \approx -11.63.$ 

\begin{lemma}
\label{lem1216}
$\mbox{Card}(\hat f - 1)=1216$.
\end{lemma}

\noindent
{\it Proof.} Due to the versality of the space $\Xi_2$, for any sufficiently small $\varepsilon >0$ we can slightly perturb the polynomial $\hat f$ within the space $\Xi_2$ in such a way that the resulting polynomial $f_\varepsilon$ satisfies the following conditions:
\begin{itemize}
\item[--]
its critical point of class $D_4^-$ splits into one minimum point with critical value $-\varepsilon$, and three critical points with Morse index 1 and critical values $0$, $\varepsilon^2,$ and $-\varepsilon^2$;
\item[--] the critical values at its two critical points with Morse index 2 are equal to $\omega_1$ and $\omega_1+\varepsilon^2$,
\item[--] the critical values at the non-real critical points are equal to $\omega_2+ i \varepsilon^2$ and $\omega_2-i\varepsilon^2$.
\end{itemize}

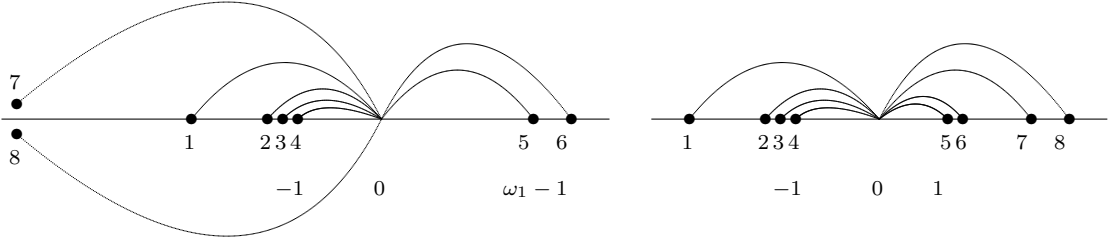
\begin{figure}
\unitlength 0.90mm
\linethickness{0.4pt}
\begin{center}
\begin{picture}(85,40)
\put(0,20){\line(1,0){80}}
\put(25,20){\circle*{1.5}}
\put(35,20){\circle*{1.5}}
\put(73,20){\circle*{1.5}}
\put(75,20){\circle*{1.5}}
\bezier{200}(50,20)(37,35)(25,20)
\bezier{160}(50,20)(42,28)(35,20)
\bezier{150}(50,20)(42,25)(37,20)
\bezier{150}(50,20)(42,23)(39,20)
\bezier{180}(50,20)(61,33)(73,20)
\bezier{200}(50,20)(62,40)(75,20)
\put(37,20){\circle*{1.5}}
\put(39,20){\circle*{1.5}}
\put(2,22){\circle*{1.5}}
\put(2,18){\circle*{1.5}}
\bezier{300}(50,20)(35,50)(2,22)
\bezier{300}(50,20)(35,-10)(2,18)
\put(49,10){\bf {\tiny $0$}}
\put(24,16){{\tiny 1}}
\put(34,16){{\tiny 2}}
\put(36,16){{\tiny 3}}
\put(38,16){{\tiny 4}}
\put(72,16){{\tiny 5}}
\put(66,10){{\tiny $\omega_1-1$}}
\put(34,10){{\tiny $-1$}}
\put(75,16){{\tiny 6}}
\put(1,24){{\tiny 7}}
\put(1,14){{\tiny 8}}
\end{picture} \quad
\begin{picture}(60,40)
\put(0,20){\line(1,0){60}}
\put(5,20){\circle*{1.5}}
\put(15,20){\circle*{1.5}}
\put(53,20){\circle*{1.5}}
\put(55,20){\circle*{1.5}}
\bezier{200}(30,20)(17,35)(5,20)
\bezier{160}(30,20)(22,28)(15,20)
\bezier{150}(30,20)(22,25)(17,20)
\bezier{150}(30,20)(22,23)(19,20)
\bezier{180}(30,20)(41,33)(53,20)
\bezier{200}(30,20)(40,40)(55,20)
\put(17,20){\circle*{1.5}}
\put(19,20){\circle*{1.5}}
\put(39,20){\circle*{1.5}}
\put(41,20){\circle*{1.5}}
\bezier{150}(30,20)(35,24)(39,20)
\bezier{150}(30,20)(36,26)(41,20)
\put(29,10){\bf {\tiny $0$}}
\put(4,16){{\tiny 1}}
\put(14,16){{\tiny 2}}
\put(16,16){{\tiny 3}}
\put(18,16){{\tiny 4}}
\put(14,10){{\tiny $-1$}}
\put(52,16){{\tiny 7}}
\put(55,16){{\tiny 8}}
\put(38,16){{\tiny 5}}
\put(38,10){{\tiny $1$}}
\put(40,16){{\tiny 6}}
\end{picture}
\caption{System of paths for polynomials $\bar f$ (left) and $f!$ (right)}
\label{sp}
\end{center}
\end{figure}

Denote the polynomial $f_\varepsilon - 1$ by $\bar f$. It belongs to the same component of the discriminant complement as the polynomial $\hat f - 1$.
The canonical system of paths in $\C^1$ defining the frame of vanishing cycles in $H_2(\bar f^{-1}(0))$ then looks as shown in Fig.~\ref{sp} (left), where the numbers 1--8 denote the order of basic vanishing cycles corresponding to the critical values.

The virtual function associated with the polynomial $\bar f$ is among the virtual functions of class $P_8^2$ listed by our program. Based on the information on the polynomial $\bar f$, let us try to find it in this list. According to \S~5.2 of \cite{vasP8} there are exactly 1897 virtual functions of class $\Xi_2$ with six real critical points. Only 140 of these virtual functions are associated with real polynomials that have a minimum point. Each of these haves, in addition to the minimum, three critical points with Morse index 1 and two points with Morse index 2. Exactly $20=140/7$ of these virtual functions have four negative and two positive critical values at real points. Of these 20 virtual functions, only 12 have all these positive critical values at the critical points with Morse index 2. Thus, the virtual function associated with $\bar f$ is among these 12. Upon a special request, our program printed these 12 virtual functions. It turns out that for all of them the intersection index $\langle \Delta_7, \Delta_8\rangle$ of two cycles vanishing at non-real points is equal to 1. In particular, this is true for the real polynomial $\bar f$.
Using the Lyashko--Looijenga map in the form of \cite{Jaw2} (see the proof of Theorem 7 in \cite{vasP8}), we can move this polynomial within the discriminant complement in $\Xi_2$ in such a way that all its real critical values stay fixed, but two non-real critical values move along the corresponding paths shown in Fig.~\ref{sp} (left) almost until the point 0, then travel close to the real axis to the point $1 \in \R$ and finally meet at this point. Due to the equality $\langle \Delta_7, \Delta_8\rangle=1,$ the critical points corresponding to these values also meet at a real critical point of class $A_2$. With a slight further perturbation, we obtain a generic polynomial with eight real critical points, whose six critical values are the same as those $\bar f$, and the other two critical points, obtained from the $A_2$ point, have real critical values $1-\zeta$ and $1+\zeta$ with some very small $\zeta>0$. Denote this polynomial by $f!$, see Fig.~\ref{sp} (right). Its critical point corresponding to the critical value $-\varepsilon -1$ and marked with 1 in this figure is a minimum point. The critical points with values marked with 2, 3, and 4 have Morse index 1, and the critical points marked with 7 and 8 have Morse index 2.
According to Lemmas V.3.2 and V.3.3 of \cite{APLT}, the newborn critical point with critical value $1- \zeta$ (marked with 5 in the picture) has an odd Morse index, and the point with critical value $1 + \zeta$ has an even Morse index. According to Proposition 43 of \cite{vasP8}, a generic polynomial of class $\Xi_2$ with eight real critical points cannot have more than one point of local extremum, therefore these odd and even Morse indices are in fact equal to one and two, respectively.

Proposition 43 of \cite{vasP8} implies that the intersection form of any real polynomial of class $\Xi_2$ with eight critical points, one of which is a minimum, is defined by the Coxeter-Dynkin graph shown in Fig.~\ref{504}, in such a way that
\begin{itemize}
\item[--] the lowest vertex in this figure corresponds to the minimum point, the vertices in the middle row correspond to critical points with Morse index 1, and the upper row vertices correspond to points with Morse index 2;
\item[--] for each pair of vertices of the graph that are connected by an edge, the critical value of the critical point corresponding to the higher vertex in this figure is greater than the critical value of the critical point corresponding to the lower vertex.
\end{itemize}

The virtual function associated with the polynomial $f!$ is determined by this graph supplied with the information on the Morse indices of critical points corresponding to its vertices, and with the numbering of these vertices corresponding to the ascending order of critical values, i.e., the correspondence of these vertices with the black points in Fig.~\ref{sp} (right). The graph and the Morse indices are known to us, while concerning the numbering we know the following.
The lowest vertex in Fig.~\ref{504} has number 1, the vertices in the middle row have numbers 2, 3, 4, and 5, and the vertices in the upper row have the numbers 6, 7, and 8. The ``leaf'' vertex of the second row (connected with only one other vertex) cannot correspond to the point 5. Indeed, the intersection index $\langle \Delta_5, \Delta_6 \rangle$ of cycles vanishing in two newborn Morse critical points is equal to 1, see Lemma V.3.2 in \cite{APLT}.
In particular, the positive critical values of the polynomial $f!$ correspond to all vertices in the top row of Fig.~\ref{504} and one vertex of valence three in the middle row. 

 All virtual functions with intersection matrix defined by the graph of Fig.~\ref{504} and the Morse indices and numbering of vertices satisfying these restrictions belong to the same virtual component. Indeed, they can all be obtained from each other by transpositions of orders of vertices that are not connected by edges of the graph and correspond to the critical values of the same sign. Such a transposition is a flip of virtual functions which corresponds to the elementary function surgery when one critical value overtakes another critical value at a distant critical point. This surgery is not related with crossing the discriminant. Our program, when applied to one of the virtual functions from this set, says that this virtual component has 1216 vertices. \hfill $\Box$ 
\medskip

It remains to prove that the connected component of the discriminant complement that contains the polynomial $f!$ is $S(3)$-invariant. Let $L\sim S^1$ again be the set of all linear functions $\R^3 \to \R$ of unit norm that vanish at the vector $(1,1,1)$. 

\begin{lemma}
For each $l \in L$, the polynomial 
\label{lemccd}
\begin{equation}
\label{ccd}
F(l) \equiv x^3+y^3+z^3 - 6x y z + 3 \sqrt{2} l^2
\end{equation}
has a critical point of class $D_4^-$ with zero critical value at the origin. All its other critical values at real critical points are positive.
\end{lemma}

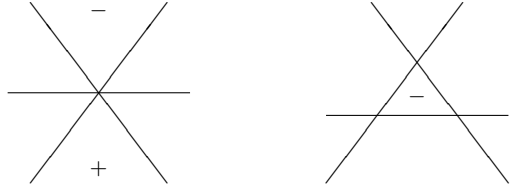
\begin{figure}
\begin{center}
\unitlength 0.6mm
\begin{picture}(50,40)
\put(5,0){\line(3,4){30}}
\put(35,0){\line(-3,4){30}}
\put(0,20){\line(1,0){40}}
\put(18,37){\tiny $+$}
\put(18,2){\tiny $-$}
\end{picture} \qquad \qquad \qquad
\begin{picture}(40,40)
\put(10,0){\line(3,4){30}}
\put(30,0){\line(-3,4){30}}
\put(-4,27){\line(1,0){48}}
\put(18,20){\tiny $-$}
\end{picture}
\caption{$D_4^-$ singularity and its perturbation}
\label{perd4}
\end{center}
\end{figure}

\noindent
{\it Proof.} The restriction of the polynomial (\ref{ccd}) to the plane $\{l=0\}$ is a non-degenerate homogeneous cubic polynomial that vanishes on three real lines, see Fig.~\ref{perd4} (left). The first assertion of the lemma follows from this analogously to the proof of Lemma~\ref{le338}. Let $a \in \R^3$ be an arbitrary critical point of this polynomial outside the plane $\{l=0\}$. The restriction of the polynomial to the line through $0$ and $a$ is a polynomial function of degree three with two critical points $0$ and $a$, that has a local minimum with critical value 0 at the point $0$. \hfill $\Box$ \medskip

The set of all polynomials $F(l)$, $l \in L$, is parametrized by the points of the circle $L/\Z_2 \sim S^1$, because $F(l)=F(-l)$.
Let $c(l)$ be the smallest critical value of the polynomial (\ref{ccd}) at its real critical points different from the origin. These numbers $c(l)$ form a lower semi-continuous function on $L/\Z_2$. By the previous lemma, this function is positive everywhere. Therefore, it is separated from zero by a constant $\delta >0$. The union of the points in $\Xi_2$ corresponding to the polynomials $F(l) -\delta$ over all $l \in L$ does not intersect the discriminant. This union is connected and invariant under the action (\ref{symP}) of the group $S(3)$. By Lemma \ref{lem1216}, for $l \equiv \frac{x-y}{\sqrt{2}}$ the polynomial $F(l)-1$ belongs to the connected component of the discriminant complement with $\mbox{Card}=1216$. The polynomial $F(l)-\delta$ also belongs to this component because all the polynomials $F(l)-t$, where $t \in [\delta,1]$, are non-discriminant. Therefore, the entire component is $S(3)$-invariant. This concludes the proof of Theorem \ref{mthmp2} for the virtual component with $\mbox{Card}=1216$. 
\medskip

\noindent
{\bf Components with $\mbox{Card}=1318$ and $844$.} 
Denote by $\{D_4^-\}$ the stratum of the discriminant variety in $\Xi_2$ that consists of polynomials having a singularity of class $D_4^-$ and no other real critical points with critical value 0. By the versality of the deformation (\ref{versP}), this stratum is a smooth submanifold of codimension four in $\Xi_2$. All polynomials $F(l)$ of the form (\ref{ccd}) belong to this stratum.
 
Consider a tubular neighborhood of this stratum in $\Xi_2$ fibered into its orthogonal slices diffeomorphic to small four-dimensional balls centered at the points of the stratum. Suppose this tubular neighborhood is narrow enough that, for any $l\in L$, the fiber over the point $\{F(l)\}$ of the stratum is the parameter space of a versal deformation of the singularity of the function $F(l)$ at the origin. Then, all the intersections of these fibers with the discriminant 
are locally ambient isotopic to the discriminant set of the model versal deformation 
\begin{equation}
\label{vdD4}
X^3 - 3X Y^2 + Z^2 + \varkappa_1 + \varkappa_2 X + \varkappa_3 Y + \varkappa_4 (X^2+Y^2)
\end{equation}
of the sample polynomial $X^3 -3X Y^2+ Z^2$ of class $D_4^-$ with parameters $\varkappa_i$.

\begin{definition} \rm
\label{defUVW}
For each $l \in L$, denote the fiber of this tubular neighborhood over the point $\{F(l)\}$ by $U(l)$, and define the subset $V(l) \subset U(l)$ as the union of all Morse polynomials $f \in U(l)$ such that 
\begin{itemize}
\item[(a)] $f$ has four real critical points near the origin that are obtained by splitting the critical point of the polynomial $F(l)$: specifically, one minimum and three critical points with Morse index 1,
\item[(b)] $\mbox{Card}(f) = 844$.
\end{itemize}
Similarly, define the subset $W(l) \subset U(l)$ as the union of polynomials $f$ satisfying condition (a) but with condition (b) replaced by $\mbox{Card}(f) = 1318$.
\end{definition}

\begin{remark} \rm
\label{re16}
According to Proposition \ref{proP2} and Lemma \ref{lemccd}, the set $V(l) \cup W(l)$ is the union of all polynomials in $U(l)$ with $\mbox{Ind}=-1$ that satisfy condition (a).
\end{remark}

\begin{lemma}
\label{lem1318}
For each $l \in L$, the subset $V(l)$ is contractible, and the subset $W(l)$ consists of two contractible connected components. 
\end{lemma}

\noindent
{\bf Preliminary remarks.}
By the versality property, the discriminant is equisingular along the smooth stratum $\{D_4^-\}$. Therefore, it is sufficient to prove Lemma \ref{lem1318} for only one particular $l \in L$, say, for $l \equiv \frac{x-y}{\sqrt{2}}$.

The union of subsets $V(l)$ over all $l \in L$ is $S(3)$-invariant. Lemma \ref{lem1318} implies that it is also connected, and therefore the corresponding virtual component with $\mbox{Card}=844$ is associated with only one real connected component of the discriminant complement. Applying the involution (\ref{invoP}) proves the assertion of Theorem \ref{mthmp2} also for the other virtual component with $\mbox{Card}=844.$ The union of subsets $W(l)$ is also $S(3)$-invariant. It is fibered over the circle $L/{\Z_2}$ with the fiber consisting of two contractible components. In Lemma \ref{monod}, we will see that the monodromy of this fibration permutes these components. Therefore, their union over all $l \in L$ is connected, and we again have only one real connected component. 
\medskip

To prove Lemma \ref{lem1318}, we will describe the geometry of the discriminant in the parameter space of the versal deformation (\ref{vdD4}) in detail. The Morse polynomials of the form (\ref{vdD4}) can have two or four critical points depending on the parameters $\varkappa_2, \varkappa_3,$ and $\varkappa_4$. In a neighborhood of the origin of the space $\R^3$ of these parameters, the set of non-Morse polynomials looks as shown in Fig.~\ref{pyramid} (left), see \cite{AVGZ}, \S~21.3.

\begin{figure}
\unitlength 0.8mm
\linethickness{0.4pt}
\begin{center}
\begin{picture}(35.00,50.00)
\put(17.20,22.80){\circle{1.00}}
\bezier{170}(17.30,23.40)(22.00,31.00)(18.00,42.00)
\bezier{250}(17.30,23.40)(24.00,34.00)(26.00,51.00)
\bezier{200}(17.30,23.40)(24.00,34.00)(33.00,41.00)
\bezier{110}(18.00,42.00)(23.00,43.00)(26.00,51.00)
\bezier{110}(33.00,41.00)(28.00,43.00)(26.00,51.00)
\bezier{170}(16.70,22.60)(14.00,16.00)(2.00,11.00)
\bezier{170}(2.00,11.00)(6.00,4.00)(19.00,3.00)
\bezier{160}(19.00,3.00)(8.00,3.00)(0.00,0.00)
\bezier{100}(0.00,0.00)(2.00,3.00)(2.00,11.00)
\bezier{70}(0.00,0.00)(4.00,2.00)(7.50,5.60)
\bezier{200}(19.00,3.00)(12.00,10.00)(16.70,22.60)
\bezier{25}(20.00,42.00)(22.00,42.00)(23.00,41.90)
\bezier{30}(26.00,41.85)(27.00,41.84)(30.50,41.50)
\bezier{60}(8.7,7.3)(12,11)(14,16)
\bezier{150}(20,0)(30,-3)(40,0)
\put(20,0){\vector(-3,1){0.1}}
\put(40,0){\vector(3,1){0.1}}
\end{picture}
\quad
\unitlength 1.2 mm
\begin{picture}(40.00,36.00)
\bezier{300}(20.00,36.00)(20.00,12.00)(0.00,0.00)
\bezier{300}(20.00,36.00)(20.00,12.00)(40.00,0.00)
\bezier{300}(0.00,0.00)(20.00,6.00)(40.00,0.00)
\put(0,0){\line(2,1){26.2}}
\put(40,0){\line(-2,1){26.2}}
\put(20,32){\line(0,-1){29}}
\put(16.5,5){\makebox(0,0)[cc]{\tiny $a b c$}}
\put(23.5,5){\makebox(0,0)[cc]{\tiny $a c b$}}
\put(15,10){\makebox(0,0)[cc]{\tiny $b a c$}}
\put(25,10){\makebox(0,0)[cc]{\tiny $c a b$}}
\put(17.5,14){\makebox(0,0)[cc]{\tiny $b c a$}}
\put(22.5,14){\makebox(0,0)[cc]{\tiny $c b a$}}
\end{picture} \qquad
\begin{picture}(40.00,36.00)
\bezier{300}(20.00,36.00)(20.00,12.00)(0.00,0.00)
\bezier{300}(20.00,36.00)(20.00,12.00)(40.00,0.00)
\bezier{300}(0.00,0.00)(20.00,10.00)(40.00,0.00)
\put(20,10){\line(2,1){6}}
\put(20,10){\line(-2,1){6}}
\put(20,10){\line(0,-1){5}}
\put(15,8){\makebox(0,0)[cc]{\footnotesize $c$}}
\put(25,8){\makebox(0,0)[cc]{\footnotesize $b$}}
\put(20,16){\makebox(0,0)[cc]{\footnotesize $a$}}
\end{picture}
\caption{Caustic and Maxwell set for $D_4^-$ singularity}
\label{pyramid}
\end{center}
\end{figure}
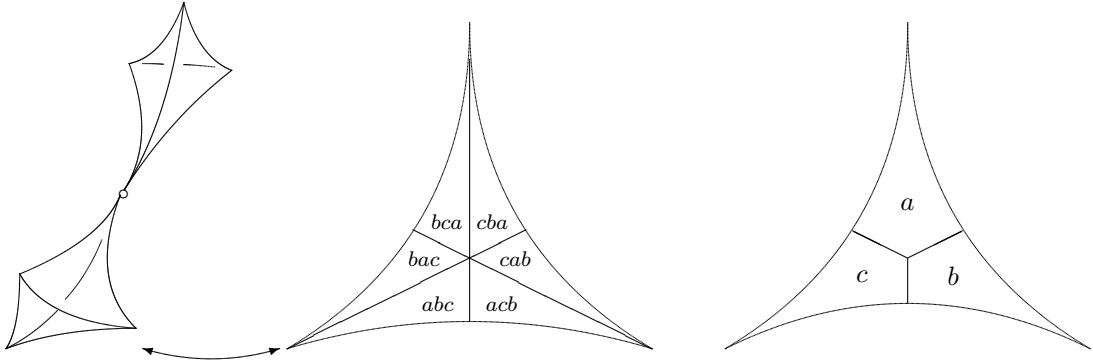
This set divides the space $\R^3$ into three parts. The larger part consists of polynomials with only two real critical points. Two symmetric ``pyramids'' are filled with polynomials (\ref{vdD4}) that have three critical points with Morse index 1 and one critical point with Morse index 0 or 2, depending on the pyramid. We are only interested in the polynomials (\ref{vdD4})
with a minimum point. For each {\em generic} interior point $(\varkappa_2, \varkappa_3, \varkappa_4)$ of the corresponding pyramid, the set of non-discriminant polynomials (\ref{vdD4}) with these values or parameters consists of five intervals. Indeed, the coefficient $\varkappa_1$ of a polynomial from this set should be not equal to minus a critical value of the polynomial 
\begin{equation}
\label{trun}
X^3 - 3X Y^2 + Z^2 + \varkappa_2 X + \varkappa_3 Y + \varkappa_4 (X^2+Y^2),
\end{equation} 
i.e., it can belong to one of the five intervals into which these critical values divide the line $\R$.
For exceptional interior points of the pyramid, the number of intervals is less than five; this occurs when some of critical values of the polynomial are equal. The set of these non-generic polynomials (\ref{trun}) divides the pyramid into six parts. 
 The intersection of the pyramid with a plane orthogonal to its axis looks as the curvilinear triangle shown in Fig.~\ref{pyramid} (center). The non-generic polynomials form three segments in it. In particular, the central point of this triangle is represented by the polynomial of the form (\ref{trun}) with $\varkappa_2=\varkappa_3=0$ and equal critical values at all three points with Morse index 1. These critical values are equal to $\frac{4}{27}\varkappa_4^3$. The orthogonal plane sections of the pyramid are defined by the conditions $\varkappa_4 = C$ for different constants $C>0$. The set of all polynomials (\ref{vdD4}) having a minimum point and three critical points with critical value zero consists of one curve in $\R^4$ satisfying the conditions
\begin{equation}
\label{con7}
\varkappa_2 = \varkappa_3 = 0, \quad \varkappa_1 = - \frac{4}{27} \varkappa_4^3, \quad \varkappa_4>0.
\end{equation}
The zero-level sets of all such polynomials in the plane $\{Z = 0\}$ look as in Fig.~\ref{perd4} (right). 

Choose an arbitrary polynomial inside the pyramid and 
label its three critical points of index 1 as $a, b,$ and $c$. This labeling of critical points can be uniquely extended by continuity to all other polynomials within the pyramid. Each of the six domains, into which the curvilinear triangle in Fig.~\ref{pyramid} (center) is split by the segments, is characterized by the order of the critical values of the corresponding polynomials at the critical points labeled in this way. 

The set of non-discriminant polynomials (\ref{vdD4}) whose parameters $(\varkappa_2, \varkappa_3, \varkappa_4)$ belong to this pyramid consists of nine connected domains: the number of domains with 0, 1, 2, 3, or 4 negative critical values is equal to 1, 1, 3, 3, or 1, respectively. We are only interested in three of these domains, whose points are the polynomials (\ref{vdD4}) which have exactly one positive critical value. The projections of these domains to the space $\R^3$ of coefficients $(\varkappa_2, \varkappa_3, \varkappa_4)$ are three sectors of the pyramid: their intersections with a section of the pyramid are shown in Fig.~\ref{pyramid} (right). Each sector, marked by a letter \ $a, b,$ or $c$ \ in the picture, consists of polynomials (\ref{trun}) whose critical value at the critical point with that marking is the highest.
The corresponding domain in the parameter space $\R^4$ of the family (\ref{vdD4}) is a fiber bundle over the corresponding sector shown in Fig.~\ref{pyramid} (right), where the fiber over a polynomial (\ref{trun}) is equal to the interval between two its highest critical values. Clearly, all three of these domains are contractible. The set of polynomials (\ref{con7}) belongs to the boundaries of all these domains.

Thus, we have the following statement.

\begin{proposition}
\label{pro19}
There are natural one-to-one correspondences between the following three-element sets:
\begin{enumerate}
\item
The set of connected components of the space of non-dis\-cri\-mi\-nant polynomials $($\ref{vdD4}$)$ having a minimum point and three critical points with Morse index 1, only one of whose critical values is positive;
\item
The set of critical points with Morse index 1 of an arbitrary polynomial $($\ref{vdD4}$)$ having a minimum point;
\item
The set of pairs of lines of the zero-level set of the polynomial $X^3-3X Y^2 + Z^2$ in the plane $\{Z=0\}.$ 
\end{enumerate} 
\end{proposition}

\noindent
{\it Proof.} The sets mentioned in item 2 for all polynomials (\ref{vdD4}) having a minimum point are in a standard one-to-one correspondence with each other. This correspondence is defined by the continuity within the contractible set of all such polynomials. Thus, it is sufficient to define the correspondence $1 \leftrightarrow 2$ assuming that the polynomial mentioned in item 2 belongs to this component. In this case, we associate with the connected component the critical point with a positive critical value. 
Also, it is sufficient to define the correspondence $2 \leftrightarrow 3$ assuming that the polynomial (\ref{vdD4}) satisfies the conditions (\ref{con7}). In this case, the zero-level set of such a polynomial in the plane $\{Z=0\}$ looks as in Fig.~\ref{perd4} (right). The desired correspondence associates with each critical point the pair of lines of Fig.~\ref{perd4} (left) that are parallel to the lines in Fig.~\ref{perd4} (right) passing through this point. \hfill $\Box$
\medskip

\noindent
{\it Proof of Lemma \ref{lem1318}.} For each $l \in L$, the orthogonal slice $U(l)$ of the discriminant stratum $\{D_4^-\} \subset \Xi_2$ at the point $\{F(l)\}$ is a versal deformation of the corresponding singularity of class $D_4^-$. The inducing map from the definition of versal deformations (see \cite{AVGZ}, \S~8) identifies this slice with a neighborhood of the origin in the parameter space $\R^4$ of the deformation (\ref{vdD4}) and maps the discriminant to the discriminant. This map defines a bijection between the union of three domains mentioned in item 1 of Proposition \ref{pro19} and the union $V(l) \cup W(l)$, see Definition \ref{defUVW}. 

Consider again the case $l \equiv \frac{x-y}{\sqrt{2}}$. Let $f$ be a point of the set $V(l) \cup W(l)$ it the corresponding slice $U(l)$. We assume that this point $f$ is sufficiently close to $F(l)$, so that the positive critical value $\xi$ at its highest critical point obtained by splitting the $D_4^-$ singularity is below the smallest value $\tilde \omega_1$ of all other positive critical values at real critical points.
 Then, acting as in the proof of Lemma \ref{lem1216}, we can move the polynomial $f$ within the discriminant complement in $\Xi_2$ in such a way that the critical values at the real critical points remain fixed, and the critical values at the two non-real points converge to a real value between $\xi $ and $\tilde \omega_1$. For the same reasons as in the proof of Lemma \ref{lem1216}, the corresponding critical points meet at a real critical point of class $A_2$. Next, we slightly perturb the obtained polynomial so that this point splits into two real critical points. We again obtain a polynomial with Coxeter-Dynkin graph shown in Fig.~\ref{504}. This polynomial has four critical points with Morse index 1; three of these points arise from the splitting of the $D_4^-$ singularity. The three corresponding vertices in the middle row of this graph are the leaf vertex and some two of the vertices connected to upper vertices. Only one of these three critical points has a positive critical value. 
Depending on the component containing $f$, the vertex in the graph corresponding to this critical point is or is not the leaf-like rightmost vertex. This choice uniquely determines the virtual component of the virtual function associated with the obtained polynomial (and, consequently, also with the polynomial $f$ that belongs to the same component of the discriminant complement). Indeed, the information on this polynomial known to us determines all the elements of the virtual function 
up to the surgeries not crossing the discriminant: all these surgeries consist in collisions of the critical values at distant critical points.
Starting from an arbitrary element of this virtual component, our program listed all of its elements, in particular found its $\mbox{Card}$ invariant. It turns out that $\mbox{Card}=844$ when the rightmost vertex in Fig.~\ref{504} corresponds to the critical point of $f$ with positive critical value, and $\mbox{Card}=1318$ in the other two cases. \hfill $\Box$ \medskip

 The union of the contractible domains $V(l) \subset U(l)$ over all $l \in L$ is a connected set that is $S(3)$-invariant. All of these domains consist of polynomials with $\mbox{Card}=844$. This proves the assertion of Theorem \ref{mthmp2} for the virtual components with $\mbox{Card}=844$.
\medskip

The union of all sets $W(l)$ also is $S(3)$-invariant. It is the space of a fiber bundle over the circle $L / \Z_2$, with the fiber $W(l)$ consisting of two contractible domains.

\begin{lemma}
\label{monod}
The monodromy over the cycle $L/\Z_2 \sim S^1$ permutes the two components of the sets $W(l)$ in the orthogonal four-dimensional slices $U(l)$ of the $\{D_4^-\}$ stratum at the corresponding points $\{F(l)\}$.
\end{lemma}

\begin{figure}
\begin{center}
\unitlength 0.8 mm
\begin{picture}(70,60)
\put(5,30){\line(1,0){50}}
\put(0,15){\line(1,0){50}}
\put(0,15){\line(1,3){5}}
\bezier{40}(5,30)(7.5,37.5)(10,45)
\bezier{30}(10,45)(16,45)(21.5,45)
\bezier{30}(22.5,45)(29,45)(37,45)
\bezier{40}(38.5,45)(47,45)(54,45)
\put(55,45){\line(1,0){4.5}}
\put(50,15){\line(1,3){10}}
\put(30,30){\line(1,2){12}}
\put(30,30){\line(-1,2){12}}
\bezier{200}(18,54)(18,50)(30,50)
\bezier{200}(30,50)(42,50)(42,54)
\bezier{200}(18,54)(18,58)(30,58)
\bezier{200}(30,58)(42,58)(42,54)
\bezier{30}(18,6)(18,10)(30,10)
\bezier{30}(30,10)(42,10)(42,6)
\bezier{200}(18,6)(18,2)(30,2)
\bezier{200}(30,2)(42,2)(42,6)
\bezier{40}(30,30)(26,22)(22,14)
\bezier{40}(30,30)(34,22)(38,14)
\put(18,6){\line(1,2){4.3}}
\put(42,6){\line(-1,2){4.3}}
\put(5,6){\line(1,0){50}}
\put(5,30){\line(0,1){24}}
\put(5,6){\line(0,1){8.5}}
\bezier{30}(5,15.5)(5,22)(5,29)
\put(55,6){\line(0,1){48}}
\put(55,54){\line(-1,0){50}}
\put(53,55){\tiny $l=0$}
\end{picture} \qquad \qquad
\begin{picture}(60,70)
\put(-1.5,10.5){\line(1,0){50}}
\put(-1.5,10.5){\line(1,3){12.7}}
\bezier{80}(48.5,10.5)(49.5,13.5)(50.5,16.5)
\bezier{50}(50,15)(55,30)(59.4,43.2)
\bezier{40}(59.4,43.2)(60.3,45.9)(61.2,48.6)
\bezier{100}(61.2,48.6)(57.2,48.6)(52.8,48.6)
\bezier{30}(52.5,48.6)(45,48.6)(39.5,48.6)
\bezier{40}(38.5,48.6)(28,48.6)(21,48.6)
\bezier{20}(20,48.6)(17,48.6)(11.7,48.6)
\put(30,30){\line(1,2){12}}
\put(30,30){\line(-1,2){12}}
\bezier{200}(18,54)(18,50)(30,50)
\bezier{200}(30,50)(42,50)(42,54)
\bezier{200}(18,54)(18,58)(30,58)
\bezier{200}(30,58)(42,58)(42,54)
\bezier{30}(18,6)(18,10)(30,10)
\bezier{30}(30,10)(42,10)(42,6)
\bezier{200}(18,6)(18,2)(30,2)
\bezier{200}(30,2)(42,2)(42,6)
\bezier{40}(30,30)(25,20)(20.25,10.5)
\bezier{100}(20.25,10.5)(19.25,8.5)(18,6)
\bezier{40}(30,30)(35,20)(39.75,10.5)
\bezier{100}(42,6)(41,8)(39.75,10.5)
\put(38.4,51){\line(2,5){2.1}}
\bezier{50}(30,30)(34,40)(38,50)
\bezier{50}(30,30)(26,20)(22,10)
\put(22,10){\line(-2,-5){2.5}}
\put(40.5,56.25){\line(3,-2){20.5}}
\put(40.5,56.25){\line(-3,2){20.5}}
\put(30,30){\line(3,-2){20.5}}
\put(30,30){\line(-3,2){20.5}}
\put(61,42.6){\line(-2,-5){18.2}}
\put(20,69.9){\line(-2,-5){21}}
\put(19.5,3.75){\line(3,-2){10}}
\put(19.5,3.75){\line(-3,2){9.5}}
\bezier{16}(8.8,10.75)(2.8,14.75)(0.1,16.65)
\bezier{30}(0.1,16.65)(-0.65,17.15)(-1.25,17.55)
\put(27,66){\tiny $l=0$}
\end{picture}
\end{center}
\caption{Plane sections of the surface $x^3+y^3+z^3-6x y z=0$}
\label{sec1}
\end{figure}
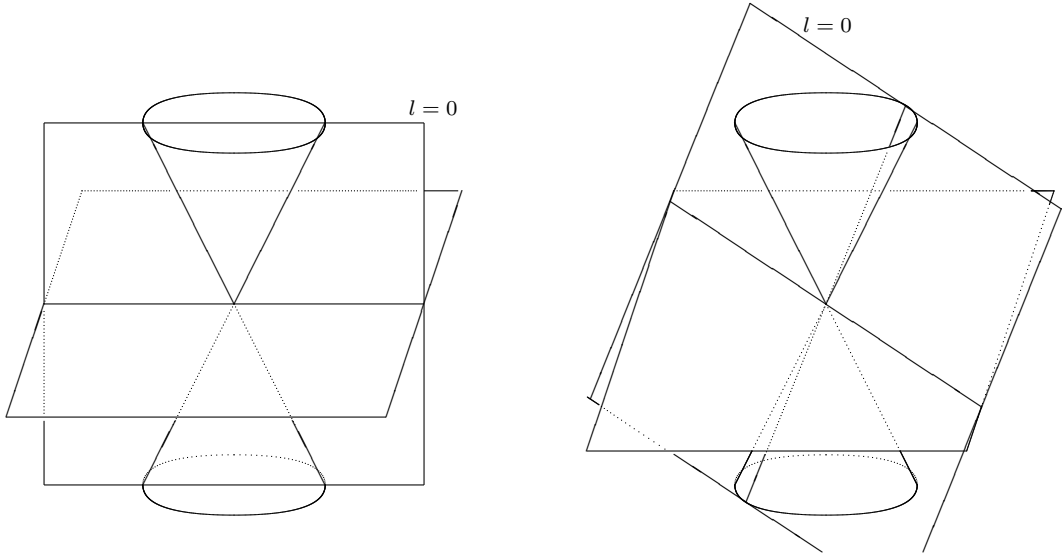

\noindent
{\it Proof.} For any $l \in L$, the subset in $\R^3$ defined by equations
\begin{equation}
l=0, F(l)=0
\label{3lines}
\end{equation}
consists of three lines, see Fig.~\ref{perd4}. These lines have different nature. Two of them belong to the component of the set where
\begin{equation}
\label{equ6}
x^3+y^3+z^3 - 6 x y z =0
\end{equation}
that is homeomorphic to the quadratic cone in $\R^3$. The third line belongs to the component that is homeomorphic to the plane, see Fig.~\ref{sec1} (left) (in which we schematically depicted the set (\ref{equ6}) by the union of the quadratic cone and the horizontal plane). 

Proposition \ref{pro19} reduces the monodromy action on the set of components of the space $V(l) \cup W(l) \subset U(l)$ to the action on the set of pairs of these lines. For any $\varkappa =(\varkappa_1, \dots, \varkappa_4)$, denote by $g_\varkappa$ the corresponding polynomial (\ref{vdD4}). According to the definition of versal deformations, there is a neighborhood $\tilde U$ in the space $\R^4$ of parameters $\varkappa$, a diffeomorphism $\theta: \tilde U \to U(l)$, and a family of local diffeomorphisms $\Gamma_\varkappa$ from a neighborhood of the origin in the space of arguments $(X, Y, Z)$ of polynomials (\ref{vdD4}) to a neighborhood of the origin in the space of arguments $(x,y,z) $ of polynomials of class $\Xi_2$, such that the diffeomorphisms $\Gamma_\varkappa$ depend smoothly on the parameters $\varkappa \in \tilde U$, and for any $\varkappa \in \tilde U$ we have
\begin{equation}
g_\varkappa(X,Y,Z) \equiv f_{\theta(\varkappa)} \circ \Gamma_\varkappa (X, Y, Z).
\end{equation}

The diffeomorphism $\Gamma_0$ can be considered as a change of local coordinates that reduces the polynomial $F(l)$ to the normal form $X^3-3X Y^2+ Z^2$. The kernel of the second differential of a polynomial is independent of the coordinates, therefore, this diffeomorphism maps the plane $\{Z=0\}$ to a surface that is tangent to the plane $\{l=0\}$. It also maps the three lines distinguished by the equation $g_0 = 0$ in the plane $\{Z=0\}$ to three curves that are tangent to three lines forming the set (\ref{3lines}). The map $\theta$ defines a bijection between the set of polynomials $g_\varkappa \subset \tilde U$ having a minimum point and $\mbox{Ind}=-1$, on the one hand, and the union of the subsets $V(l) $ and $W(l)$ of the slice $U(l)$, on the other hand. According to Proposition \ref{pro19}, we obtain a one-to-one correspondence between the pairs of lines forming the set (\ref{3lines}) and the set of connected components of the space $V(l) \cup W(l)$. 

This correspondence does not depend on the choice of maps $\Gamma_\varkappa$ and $\theta$ reducing the versal deformations (\ref{vdD4}) and $U(l)$ to one another. Indeed, let $\bar f \in U(l)$ be an arbitrary polynomial very close to $F(l)$, which has a minimum point and three critical points $a, b, $ and $c$ with Morse index 1 and critical value 0. (Such points form a single curve in $U(l)$ which is the image of the curve (\ref{con7}) under the map $\theta$.) To each of the points $a, b, c$, \ a connected component of the set $V(l) \cup W(l)$ corresponds: it contains all polynomials neighboring to $\bar f$ and such that the critical values of their critical points neighboring to the selected critical point of $\bar f$ is positive, and the critical values at the other two critical points with Morse index 1 are negative. Also, a pair of lines forming the set (\ref{3lines}) corresponds to each point $a, b, $ and $c$. These two lines are nearly parallel to the two edges of the triangle $\triangle_{a b c} \subset \R^3$ that contain this point. The composition of these two correspondences defines an invariant description of the correspondence between the components of $V(l) \cup W(l)$ and the pairs of lines in set (\ref{3lines}). 

Thus, the monodromy action of the circle $L/{\Z_2}$ on the set of pairs of lines in (\ref{3lines}) coincides with its action on the set of corresponding connected components of the set $V(l) \cup W(l)$.

As $l$ moves along the circle $L /\Z_2$, these lines in $\ker l$ move in such a way that one pair of them returns to itself, while the other two pairs permute. Thus, the monodromy over this circle acts in the same way on the set of three related domains in $U(l)$. We know one of these domains which is moved by this monodromy to itself: it is the domain $V(l)$ considered in Lemma \ref{lem1318}. Therefore, the permuted domains belong to the set $W(l)$. \hfill $\Box$ \medskip

Theorem \ref{mthmp2} is thus proved for the virtual component with $\mbox{Card} = 1318$. 
\medskip

\noindent
{\bf Component with $\mbox{Card}=156$.}

\begin{lemma}
\label{lem51}
Let $l: \R^3 \to \R$ be an arbitrary non-zero linear function, whose kernel is tangent to the component of the surface $($\ref{equ6}$)$ that is homeomorphic to the quadratic cone, see Fig.~\ref{sec1} $($right$)$. 
The polynomial 
\begin{equation}
\label{156}
x^3 + y^3 + z^3 -6xyz + l^2 
\end{equation} 
then has a singularity of class $D_5^-$ at the origin in $\R^3$. In particular, it has the form 
\begin{equation}
\label{equ65}
\zeta^2 + \xi^2 \eta - \eta^4
\end{equation}
 in some local coordinates at this point. The critical values of the polynomial $($\ref{156}$)$ at all its other real critical points are positive.
\end{lemma} 

\noindent
{\it Proof.} All singularities of corank 2 that can occur in the deformation (\ref{versP}) of $P_8$ singularities are of class $D_4$, $D_5$ or $E_6$, see \cite{siersma}, \cite{Jaw2}. The restriction of the polynomial $($\ref{156}$)$ to the plane $\ker l$ is a degenerate cubic form. Therefore, its singularity at the origin is more complicated than $D_4$. 
According to Corollary 10.4.9 of \cite{LL}, this singularity is not of class $E_6$. Thus, it remains to 
choose between the classes $D_5^+$ and $D_5^-$. 
In the particular case when 
\begin{equation}
\label{partic}
l \equiv 2x-5y+11z, 
\end{equation}
the substitutions 
$\tilde l = l + \left( \frac{3}{16} (5y-11z)^2 -\frac{3}{2} y z\right)$ and $\tilde y=y-3z$
reduce the polynomial (\ref{156}) to a polynomial whose lower weighted homogeneous part is $ \tilde l^2 + \frac{63}{2}\tilde y^2z - \frac{9}{4}z^4$. A dilation of coordinates reduces it to the form (\ref{equ65}), and technique of \cite{AVGZ}, \S 12, removes all higher weighted homogeneous terms of the Taylor expansion.

The set of all polynomials (\ref{156}) such that the plane $\ker l$ is tangent to the same component of the zero-level set of (\ref{equ6}) is connected. Therefore, the singularity will be of the same class for all other functions (\ref{156}). The positivity of other critical values follows from the same arguments as in the proof of Lemma \ref{lemccd}. \hfill $\Box$ \medskip

\begin{lemma}
\label{lele}
For any linear function $l$ as in Lemma \ref{lem51}, the set of polynomials with $\mbox{Ind}=-3$ has exactly one local connected component near
the corresponding polynomial $($\ref{156}$)$ in the space $\Xi_2$. 
\end{lemma}

\noindent
{\it Proof.} By the versality of the deformation (\ref{versP}), the discriminant is equisingular along its smooth stratum $\{D_5^-\} \subset \Xi_2$. Therefore, it is sufficient to consider only a small five-dimensional transversal slice of this stratum in $\Xi_2$ at the point (\ref{156}). Any such slice realizes a versal deformation of the corresponding singularity. 
The components of the discriminant complements of simple singularities in the spaces of their versal deformations are listed explicitly in \cite{isr}. For $D_5^-$ singularities, only one of these components consists of functions with $\mbox{Ind}=-3$. 
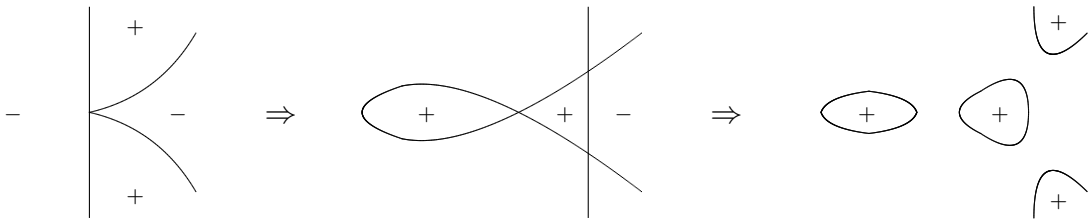
\begin{figure}
\unitlength 0.65mm
\begin{center}
\begin{picture}(44,40)
\put(14,0){\line(0,1){40}}
\bezier{150}(34,5)(27,17)(14,20)
\bezier{150}(34,35)(27,23)(14,20)
\put(0,18.5){\tiny $-$}
\put(29,18.5){\tiny $-$}
\put(21,35){\tiny $+$}
\put(21,3){\tiny $+$}
\put(47,18){$\Rightarrow$}
\end{picture} \qquad
\begin{picture}(70,40)
\put(50,0){\line(0,1){40}}
\bezier{300}(60,35)(30,12)(15,15)
\bezier{300}(60,5)(30,28)(15,25)
\bezier{250}(15,25)(0,20)(15,15)
\put(55,18.5){\tiny $-$}
\put(18,18.5){\tiny $+$}
\put(44,18.5){\tiny $+$}
\put(73,18){$\Rightarrow$}
\end{picture} \qquad
\begin{picture}(60,40)
\bezier{250}(50,0)(50,15)(60,5)
\bezier{250}(50,40)(50,25)(60,35)
\bezier{150}(49,20)(49,30)(40,24)
\bezier{150}(49,20)(49,10)(40,16)
\bezier{150}(40,24)(32,20)(40,16)
\bezier{150}(28,20)(27,23)(19,24)
\bezier{150}(19,24)(11,23)(10,20)
\bezier{150}(19,16)(11,17)(10,20)
\bezier{150}(28,20)(27,17)(19,16)
\put(17,18.5){\tiny $+$}
\put(42,18.5){\tiny $+$}
\put(53,36){\tiny $+$}
\put(53,2){\tiny $+$}
\end{picture}
\caption{Perturbation of a $D_5^-$ singularity with $\mbox{Ind}=-3$}
\label{perd5}
\end{center}
\end{figure}
Specifically, in the case of any singularity in two variables that is equivalent to $\xi^2 \eta -\eta^4$ it is the component containing the perturbations whose zero-level sets are as shown in Fig.~\ref{perd5} (right). \hfill $\Box$

\begin{lemma}
\label{lel156}
All the polynomials in the set considered in Lemma \ref{lele} have $\mbox{Card}=156.$
\end{lemma}

\begin{proposition}
\label{pro156}
Let the linear function $l$ be given by $($\ref{partic}$)$. Then the set considered in Lemmas \ref{lele} and \ref{lel156} 
contains a generic polynomial with the following properties:
\begin{itemize}
\item it has exactly six real critical points: three with Morse index 1 and negative critical values, and three with Morse index 2 and positive critical values;
\item the intersection indices of vanishing cycles corresponding to the five lowest real critical values are defined by the matrix
\begin{equation}
 \left|
\begin{array}{ccccc}
$-2$ \ & 0 & 0 & 1 & 0 \\
0 & $-2$ \ & 0 & 1 & 0 \\
0 & 0 & $-2$ \ & 1 & 1 \\
1 & 1 & 1 & $-2$ \ & 0 \\
0 & 0 & 1 & 0 & $-2$ 
\end{array}
\right| \ . \label{mat156}
\end{equation} 
\end{itemize}
\end{proposition}

\noindent
{\it Proof.} According to the Gusein-Zade---A'Campo method \cite{GZ}, the intersection matrix of five cycles vanishing in the splitting of $D_5^-$ singularity shown in Fig.~\ref{perd5} is given by the Coxeter--Dynkin diagram \ \
\begin{picture}(29,9)
\unitlength 0.7mm
\put(19,4){\circle*{1.5}}
\put(0,4){\circle*{1.5}}
\put(9,4){\circle{1.5}}
\put(27,0){\circle{1.5}}
\put(27,8){\circle{1.5}}
\put(8.3,4){\line(-1,0){9}}
\put(9.7,4){\line(1,0){9}}
\put(26.2,0.5){\line(-2,1){6.4}}
\put(26.2,7.5){\line(-2,-1){6.4}}
\end{picture}, 
where white (respectively, black) circles correspond to the critical points with Morse index 1 and negative critical value (respectively, with Morse index 2 and positive critical value). Using the Lyashko--Looijenga map, we can achieve that their critical values are ordered as follows: \ \
\begin{picture}(30,9)
\unitlength 0.7mm
\put(19,4){\circle*{1.5}}
\put(0,4){\circle*{1.5}}
\put(9,4){\circle{1.5}}
\put(27,0){\circle{1.5}}
\put(27,8){\circle{1.5}}
\put(8.3,4){\line(-1,0){9}}
\put(9.7,4){\line(1,0){9}}
\put(26.2,0.5){\line(-2,1){6.4}}
\put(26.2,7.5){\line(-2,-1){6.4}}
\put(0,0){\tiny 5}
\put(8,0){\tiny 3}
\put(18,0){\tiny 4}
\put(28,-1){\tiny 2}
\put(28,6){\tiny 1}
\end{picture}, which corresponds to the matrix (\ref{mat156}). All these perturbations can be performed in an arbitrarily small neighborhood of the polynomial (\ref{156}), so that they will not create real critical points in addition to six critical points of this polynomial. 
\hfill $\Box$
\medskip

\noindent
{\it Proof of Lemma \ref{lel156}.} All these sets corresponding to different planes $\ker l$ belong to the same connected component of the discriminant complement. Therefore, we can assume that $l$ is given by (\ref{partic}), and any neighborhood of the polynomial (\ref{156}) contains a polynomial satisfying the conditions listed in Proposition \ref{pro156}.
Our program, upon a special request, listed all virtual functions of $P_8^2$ class satisfying these conditions. There are exactly six of them, and all six belong to the virtual component with $\mbox{Card}=156.$ \hfill $\Box$
\medskip

The set of all polynomials (\ref{156}) with $\ker l$ tangent to the cone is connected and $S(3)$-invariant. According to Lemma \ref{lele}, the same is true for the union of all local domains with $\mbox{Card}=156$ in the transversal slices to the $\{D_5^-\}$ stratum at the points corresponding to these polynomials. Indeed, this union is the space of a fiber bundle with connected base and fibers. This implies the statement of Theorem \ref{mthmp2} for the virtual component with $\mbox{Card}=156$. \hfill $\Box$ 

\begin{conjecture}
The connected components with $\mbox{\rm Card}$ equal to 1216, 1318, 844, or 156 have non-trivial first homology groups.
\end{conjecture}

\section{Petrovskii lacunas}
\label{lacu} 
The above enumeration of the connected components of the discriminant complements at parabolic singularities provides a complete list of 
{\em local Petrovskii lacunas} at the corresponding singular points of wavefronts
of hyperbolic PDEs. In this theory, {\em generating families} of functions play a crucial role. 
Wave functions are defined on the parameter spaces of these families via integrals
of special differential forms along certain cycles in the zero-level sets of the corresponding functions. {\em Wavefronts} are the discriminant sets of these families.
Outside of the wavefront, the integration cycles and integrand forms depend regularly on the parameters, so that the wave function is regular as well.
However, over the discriminant sets of parameters, the zero-level sets degenerate and the wave function may become irregular. 
The type of irregularity of a wave function near a wavefront point is characterized by the asymptotics of the wave function as the non-discriminant parameter values approach the point. The type of this asymptotic behavior primarily depends on the singularity class of the corresponding discriminant function and on the homology class of the integration cycle. These homology classes depend on the connected components of the wavefront complement from which we approach the wavefront point. 
Such a component is called a {\em local lacuna} if the wave function can be continued from it to a regular function on entire neighborhood of the wavefront point. 

If the wavefront point corresponds to a function with isolated singularity, then whether a neighboring component of the wavefront complement is a local lacuna is determined by the topological data contained in the virtual function associated with an arbitrary non-discriminant function from that component (see, for example, \cite{APLT}, Chapter V).

Versal deformations of stably equivalent function singularities can be naturally identified with each other. This identification establishes a one-to-one correspondence between the components of the discriminant complements of these deformations. The set of local lacunas at a discriminant point depends on the stable equivalence class of the corresponding singularity, on the parity of the number $N$ of variables, and on the parity of the positive inertia index $i_+$ of the quadratic part of the singularity. For example, the functions $f(x_1, \dots, x_n)$ and $\tilde f(x_1, \dots, x_{n+2}) \equiv f(x_1, \dots, x_n) + x_{n+1}^2 + x_{n+2}^2$ have the same set of local lacunas. On contrary, the function $\check f(x_1, \dots, x_{n+2}) \equiv f(x_1, \dots, x_n) + x_{n+1}^2 - x_{n+2}^2$ or a stably equivalent function in $n+2k-1$ variables generally has a different set.

\begin{theorem}
The numbers of local lacunas near parabolic singularities of wavefronts in versal generating families are given in Table \ref{t13}.
\end{theorem}

\begin{table}
\begin{center}
\caption{Numbers of local lacunas at parabolic singularities} \label{t13}
\begin{tabular}{|l|c|c|c|c|}
\hline Singularity & $N$ even & $N$ even & $N$ odd & $N$ odd \cr 
class & $i_+$ even & $i_+$ odd & $i_+$ even & $i_+$ odd \cr 
\hline $P_8^1$ & $0$ & $0$ & $2$ & $0$ \cr 
\hline $P_8^2$ & $0$ & $0$ & $3$ & 0 \cr 
\hline $ X_9^+$ & $1$ & 0 & $2$ & 0 \cr 
\hline $ X_9^-$ & $1$ & 0 & 0 & $2$ \cr 
\hline $X_9^1$ & 0 & 0 & 0 & 0 \cr 
\hline $X_9^2$ & 0 & $ 4$ & 0 & 0 \cr 
\hline $J_{10}^3$ & $0$ & $1$ & 0 & 0 \cr 
\hline $J_{10}^1$ & $0$ & $0$ & 0 & 0 \cr 
\hline
\end{tabular}
\end{center}
\end{table}

Almost all of these local lacunas were described in \cite{kashin}, and it was conjectured that there are no additional local lacunas. The above explicit enumeration of connected components of the discriminant complements proves this conjecture for all singularity classes except for $P_8^2$. For $P_8^2$, it gives us one new local lacuna, thus replacing the number 2 conjectured in \cite{kashin} by 3.

In our current terms, the unique lacuna of the $X_9^{+}$ singularity for even $N$ and $i_+$ is component $(a)$ of Fig.~\ref{x9+}. The two lacunas of the $X_9^+$ singularity for odd $N$ and even $i_+$ are the components $(a)$ and $(b)$ of the same figure. The lacunas for the $X_9^-$ singularities are obtained from these of $X_9^+$ singularities by multiplying of all functions by $-1$.
The four lacunas at the $X_9^2$ singularity with even $N$ and odd $i_+$ are shown in picture $(i)$ of Fig.~\ref{X92}, or in its rotation by the angle $\frac{\pi}{2}$, or in the pictures obtained by applying the operation (\ref{ze}) to these two. The unique lacuna at the $J_{10}^3$ singularity with even $N$ and odd $i_+$ is shown in picture $(m)$ of Fig.~\ref{J103}. The two local lacunas of the $P_8^1$ singularity with odd $N$ and even $i_+$ are the components with $\mbox{Ind} = 1$. Three local lacunas of the $P_8^2$ singularity with odd $N$ and even $i_+$ are also the components with $\mbox{Ind} = 1$. The only new local lacuna is the component with $\mbox{Card}=262$ for $P_8^2$ singularity.

}
\end{document}